\documentclass[11pt,a4paper]{article}
\usepackage[T1]{fontenc}
\usepackage{lmodern}
\usepackage[margin=27mm]{geometry}
\usepackage{amsmath,amssymb,amsthm,mathtools}
\usepackage{enumitem}
\usepackage{graphicx}
\usepackage{xcolor}
\usepackage{tikz}
\usetikzlibrary{arrows.meta,decorations.markings}
\usepackage[hidelinks]{hyperref}
\allowdisplaybreaks[2]
\newtheorem{theorem}{Theorem}[section]
\newtheorem{lemma}[theorem]{Lemma}
\newtheorem{proposition}[theorem]{Proposition}
\newtheorem{corollary}[theorem]{Corollary}
\theoremstyle{definition}
\newtheorem{definition}[theorem]{Definition}
\theoremstyle{remark}

\newcommand{\PP}{\mathbb P}
\newcommand{\EE}{\mathbb E}
\newcommand{\one}{\mathbf 1}
\newcommand{\NN}{\mathcal N}
\newcommand{\DD}{\mathbb D}
\newcommand{\ii}{\mathrm i}
\newcommand{\bk}{\varkappa}
\hypersetup{pdftitle={Double-dimer loops in the biased Aztec diamond converge to a diffeomorphic image of nested CLE4},
 pdfauthor={Zhongyang Li},
 pdfsubject={Global curve convergence to a diffeomorphic image of nested CLE(4)}}
\DeclareMathOperator{\diam}{diam}
\DeclareMathOperator{\CLE}{CLE}
\DeclareMathOperator{\Lam}{Lam}
\DeclareMathOperator{\Tr}{Tr}
\DeclareMathOperator{\tr}{tr}
\DeclareMathOperator{\dn}{dn}
\DeclareMathOperator{\nest}{nest}
\DeclareMathOperator{\CP}{CP}
\title{Double-dimer loops in the biased Aztec diamond\\
converge to a diffeomorphic image of nested $\CLE_4$}
\author{Zhongyang Li\\
  \small Department of Mathematics, University of Connecticut\\
  \small\href{mailto:lizhongyangmath@gmail.com}{\nolinkurl{lizhongyangmath@gmail.com}}}
\date{}
\begin{document}
\maketitle

\begin{abstract}
We superimpose two independent domino tilings of the Aztec diamond of
order $N$, with vertical weight $a>0$ and horizontal weight $1$.
Represent each domino by the edge joining the centers of its two squares.
The edges used by exactly one tiling form loops. After rescaling by $N^{-1}$, the full loop collection converges in law
to the image of the nested conformal loop ensemble $\CLE_4$ under the
diffeomorphism
\[
 w\longmapsto
 \frac{2(\Re w+\ii a\,\Im w)}{\sqrt{1+a^2}\,(1+|w|^2)} 
\]
from the unit disk onto the liquid region of the rescaled diamond.
Thus the bias rescales the horizontal and vertical coordinates of the
unbiased limit. The topology of convergence records loop multiplicities
and the order of traversal around each loop, up to reversal.
Transfer estimates in the Krawtchouk basis control slab traversals and
opposite crossings along entire sections, including where they meet
the arctic ellipse. These estimates give tightness and, together with
finite-puncture observables, identify the limit and rule out backtracking.
\end{abstract}

\section{Introduction}
\label{sec:model}

\begin{figure}[htbp]
 \centering
 \begin{tikzpicture}[x=1cm,y=1cm]
 \begin{scope}[xshift=-3.35cm,scale=0.40]
 \foreach \i in {-5,...,4} {
  \foreach \j in {-5,...,4} {
   \pgfmathtruncatemacro{\ok}{(abs(\i+0.5)+abs(\j+0.5)<=5) ? 1 : 0}
   \pgfmathtruncatemacro{\pty}{mod(\i+\j+20,2)}
   \ifnum\ok=1
    \ifnum\pty=0 \fill[black!14] (\i,\j) rectangle ++(1,1);\fi
    \draw[black!45,line width=0.3pt] (\i,\j) rectangle ++(1,1);
   \fi
  }
 }
 \foreach \i in {-5,...,4} {
  \foreach \j in {-5,...,4} {
   \pgfmathtruncatemacro{\ok}{(abs(\i+0.5)+abs(\j+0.5)<=5) ? 1 : 0}
   \pgfmathtruncatemacro{\okr}{(abs(\i+1.5)+abs(\j+0.5)<=5) ? 1 : 0}
   \pgfmathtruncatemacro{\oku}{(abs(\i+0.5)+abs(\j+1.5)<=5) ? 1 : 0}
   \ifnum\ok=1
    \ifnum\okr=1 \draw[line width=0.5pt] (\i+0.5,\j+0.5)--(\i+1.5,\j+0.5);\fi
    \ifnum\oku=1 \draw[line width=1.5pt] (\i+0.5,\j+0.5)--(\i+0.5,\j+1.5);\fi
   \fi
  }
 }
 \foreach \i in {-5,...,4} {
  \foreach \j in {-5,...,4} {
   \pgfmathtruncatemacro{\ok}{(abs(\i+0.5)+abs(\j+0.5)<=5) ? 1 : 0}
   \pgfmathtruncatemacro{\pty}{mod(\i+\j+20,2)}
   \ifnum\ok=1
    \ifnum\pty=0 \fill (\i+0.5,\j+0.5) circle (2.6pt);
    \else \draw[fill=white,line width=0.5pt] (\i+0.5,\j+0.5) circle (2.6pt);\fi
   \fi
  }
 }
 \end{scope}
 \begin{scope}[xshift=3.05cm,scale=2]
  \fill[black!6] (1,0)--(0,1)--(-1,0)--(0,-1)--cycle;
  \fill[white] (0,0) ellipse [x radius=0.5547,y radius=0.83205];
  \draw[black!70] (1,0)--(0,1)--(-1,0)--(0,-1)--cycle;
  \draw[line width=0.7pt] (0,0) ellipse [x radius=0.5547,y radius=0.83205];
  \draw[densely dashed,black!65] (-0.8,0.2)--(0.2,-0.8);
  \node at (-0.07,0.23) {$L$};
  \node[anchor=south west,inner sep=1pt] at (0.28,0.72) {$\partial\Diamond$};
 \end{scope}
 \end{tikzpicture}
 \caption{In the order-five Aztec diamond (left), edges join
 neighboring square centers. Thick edges represent vertical dominoes
 of weight $a$. For $a=3/2$ (right), the arctic ellipse bounds the
 liquid region $L$ inside the rescaled diamond with boundary
 $\partial\Diamond$. The dashed line crosses the liquid and frozen regions.}
 \label{fig:aztec-interior-dual}
\end{figure}
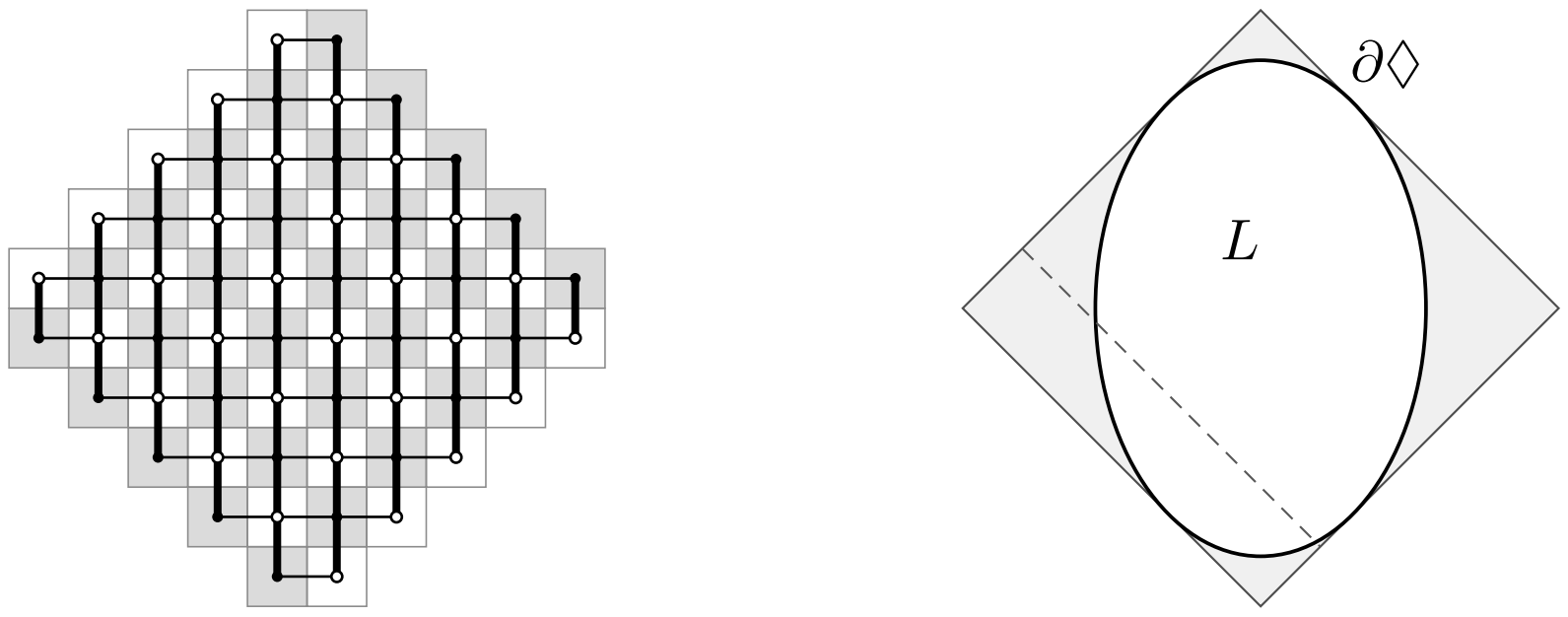

For a positive integer $N$, the Aztec diamond of order $N$ is the
union of the unit lattice squares whose centers belong to the set
\begin{equation}\label{eq:aztec-vertices}
 V_N=\bigl\{(X,Y)\in(\mathbb Z+\tfrac12)^2:|X|+|Y|\le N\bigr\}.
\end{equation}
A domino consists of two such squares sharing an edge. A domino
tiling partitions the diamond into horizontal $2\times1$ and vertical
$1\times2$ dominoes. Give vertical dominoes weight $a>0$
and horizontal dominoes weight $1$ (Figure~\ref{fig:aztec-interior-dual},
left). A tiling with $v$ vertical dominoes
then has probability proportional to $a^v$, with $a=1$ giving the
uniform measure.

\begin{figure}[tbp]
 \centering
 \begin{minipage}[t]{0.495\linewidth}
  \centering
  \includegraphics[width=\linewidth]{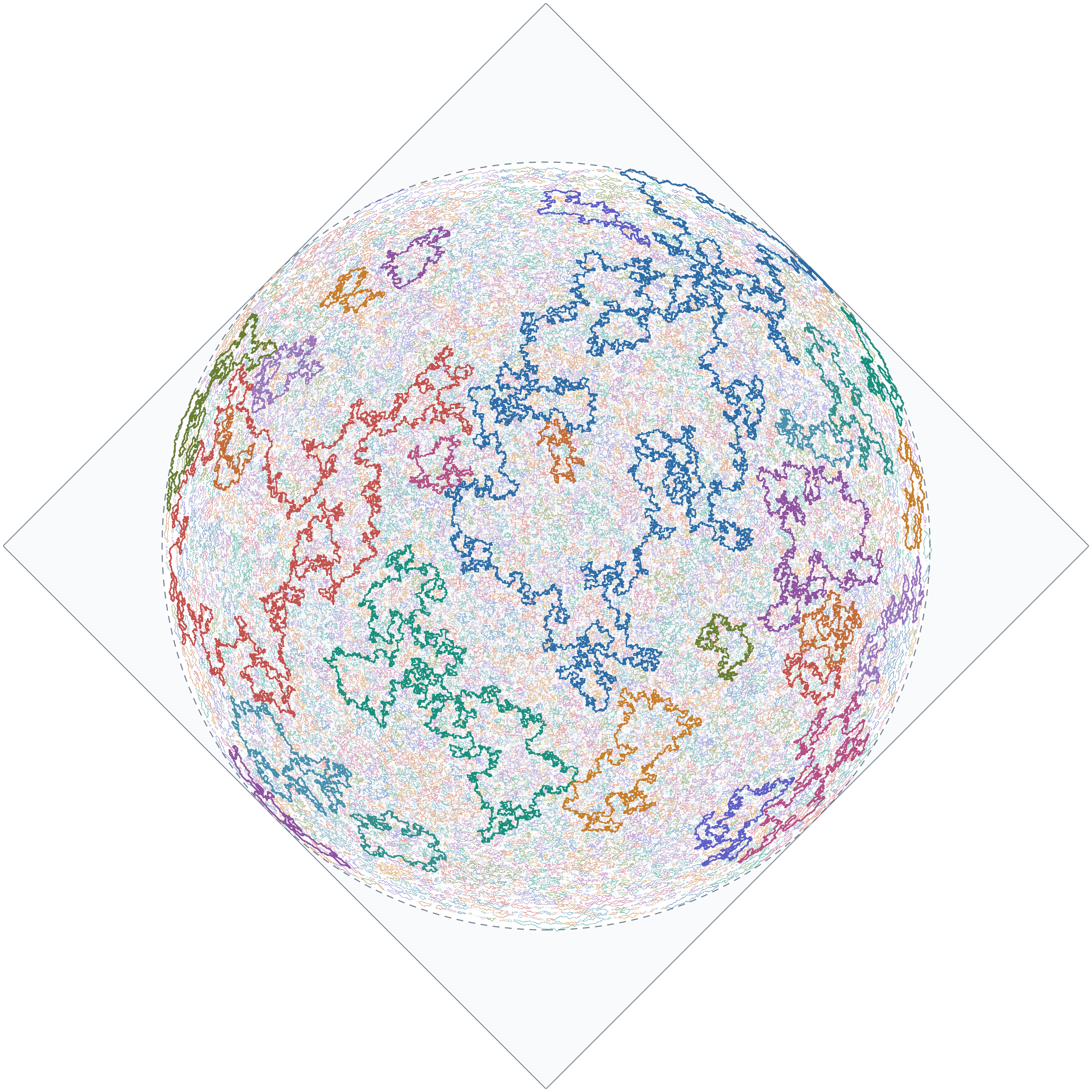}\\
  {\small(a) $a=1$}
 \end{minipage}\hfill
 \begin{minipage}[t]{0.495\linewidth}
  \centering
  \includegraphics[width=\linewidth]{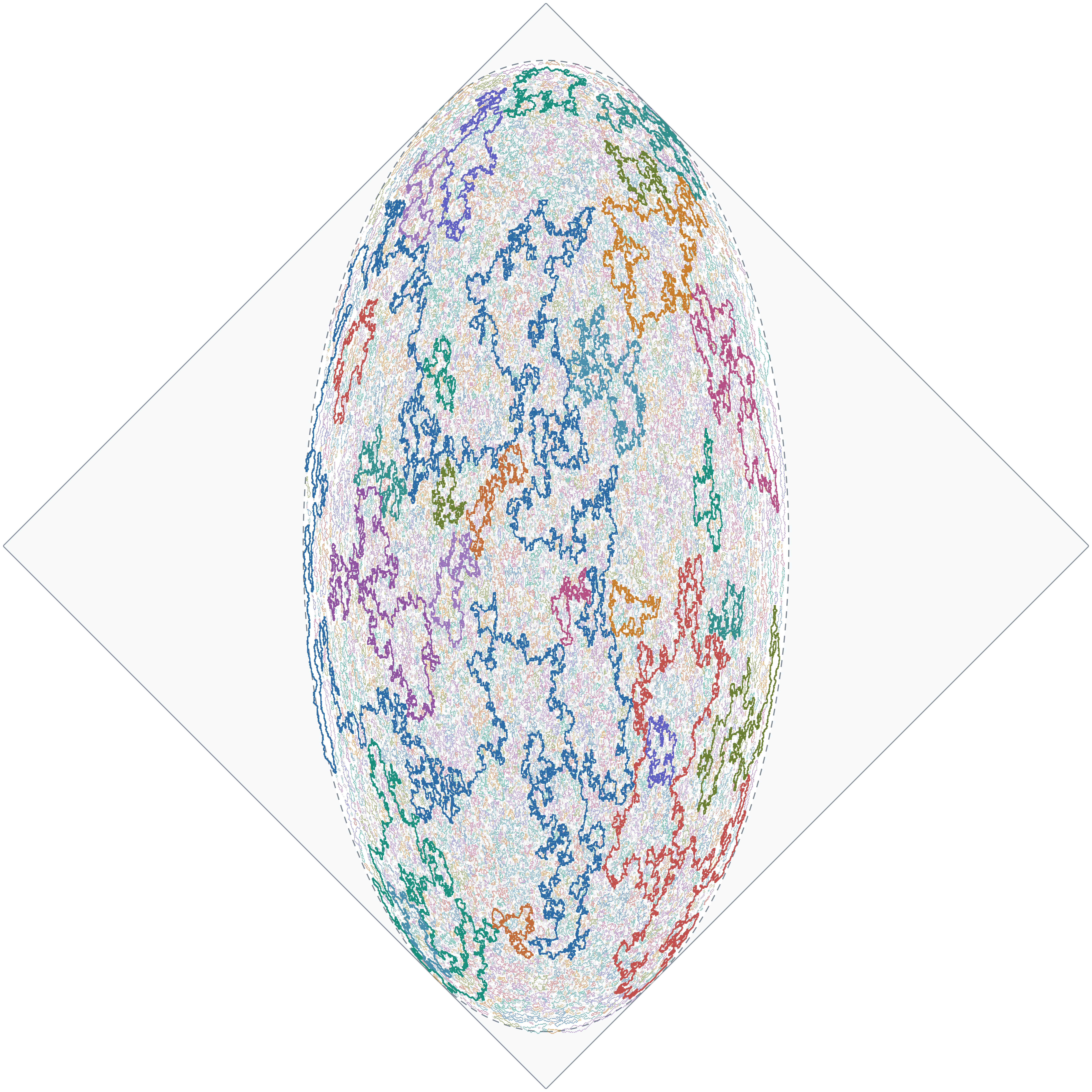}\\
  {\small(b) $a=2$}
 \end{minipage}
 \caption{Double-dimer loops on the Aztec diamond of order
 $N=1024$, in the scaled coordinates
 of~\eqref{eq:limiting-aztec-diamond}: (a)~the uniform measure,
 $69{,}503$ nondoubled loops; (b)~vertical weight $a=2$, that is
 $\bk=4$, $47{,}636$ nondoubled loops. All nondoubled loops are drawn
 and the longest ones are emphasized; the colors only distinguish the
 loops. In both panels the two tilings were sampled exactly by domino
 shuffling~\cite{EKLP,JPS,CEP,Propp}, where each new $2\times2$ block
 receives two horizontal dominoes with probability $1/(1+\bk)$ and two
 vertical ones otherwise. The dashed curves are the arctic boundaries
 $\partial L$ of~\eqref{eq:physical-liquid}: the circle of radius
 $1/\sqrt2$ in~(a) and the ellipse with semi-axes $1/\sqrt5$ and
 $2/\sqrt5$ in~(b). By Theorem~\ref{thm:global-curves}, the limit
 in~(b) is the image of the limit in~(a) under the linear map
 $\mathsf A_\bk$ of~\eqref{eq:linear-map}.}
 \label{fig:samples}
\end{figure}

Draw each domino as the segment joining the centers of its two
squares and superimpose two independent tilings sampled with these
weights. Removing the segments used by both tilings (the doubled
edges) leaves disjoint polygonal loops, called double-dimer loops.
Figure~\ref{fig:samples} shows these loops for $a=1$ and $a=2$.

The liquid region is the part of the rescaled diamond where the
local tiling remains random as $N\to\infty$. Its boundary is the
arctic ellipse~\cite{CEP} (Figure~\ref{fig:aztec-interior-dual}, right). The conformal loop ensemble with parameter~$4$,
denoted by $\CLE_4$, is a random collection of disjoint simple loops.
Sampling an independent copy inside each loop and iterating gives
nested $\CLE_4$. We prove that, after rescaling by
$N^{-1}$, the double-dimer loops converge jointly in law as curves to
the image of nested $\CLE_4$ under an explicit diffeomorphism from the
unit disk onto the liquid region. Changing $a$ rescales the horizontal
and vertical coordinates of the limiting loops.

Sheffield~\cite{Sheffield} introduced the conformal loop ensembles
$\CLE_\kappa$ using branching variants of Schramm--Loewner evolution.
For $8/3<\kappa\le4$, the loops are simple and pairwise disjoint.
Sheffield and Werner~\cite{SW} characterized the non-nested ensembles
in this range by conformal invariance and a spatial Markov property.
In particular, the outer boundaries of the outermost clusters of a
Brownian loop-soup of intensity~$1$ have the law of non-nested
$\CLE_4$. Kenyon predicted a $\CLE_4$ scaling limit for double-dimer
loops~\cite{Kenyon,Dubedat,BC}.

Kenyon's monodromy observables record the homotopy classes of
double-dimer loops relative to finitely many marked faces
\cite{Kenyon}. For the square lattice in the upper half-plane,
Dub\'edat proved that the correlators built from these observables
converge, for locally unipotent monodromies close to the identity, to
the corresponding correlators of $\CLE_4$~\cite{Dubedat}. Basok and
Chelkak showed that such correlators determine the probabilities of
cylindrical events~\cite{BC}. Combined with the crossing estimate of
Bai and Wan~\cite{BW}, this leads to the convergence of these
probabilities to those of $\CLE_4$, recorded for the half-plane
in~\cite[Theorem~4.2]{BI}. These results concern the homotopy
classes of the loops relative to finitely many points. Convergence of
the loops as curves requires tightness in addition, and no
Russo--Seymour--Welsh theory is available for double
dimers~\cite{BC,BI}.
A loop can make a long excursion and return beside its outgoing arc
without changing its class relative to any of the marked faces.
To control such excursions, we bound the probability that two arcs
of the same loop cross a short interval of a diagonal section with opposite
net signs. In the diamond the bound holds along the entire section,
including where it meets the arctic boundary.

\subsection{The biased double-dimer model}
\label{subsec:aztec-model}

Join nearest neighbors in $V_N$ from~\eqref{eq:aztec-vertices}
to obtain the Aztec diamond graph $G_N$. A tiling corresponds to a
perfect matching $M$ of $G_N$, a set of edges covering each vertex
exactly once~\cite{EKLP}. Writing $v(M)$ for the number of vertical
edges in $M$, the tiling measure becomes
\begin{equation}\label{eq:biased-measure-a}
 \PP\{M\}=a^{v(M)}/Z_N,\qquad Z_N=\sum_M a^{v(M)} .
\end{equation}
Let $M_N^{(1)},M_N^{(2)}$ be the two independent perfect matchings and
$\mathcal L_N$ the collection of their double-dimer loops scaled by
$N^{-1}$. Set $\bk=a^2$.

Color $(X,Y)$ white when $X+Y+N$ is odd and black otherwise. There
are $N(N+1)$ vertices of each color. A vertical edge is a \emph{down
edge} if its black endpoint is uppermost; let $\dn(M)$ count the
down edges of $M$. Above the line $Y=j$ for any $j\in\mathbb Z$, the two colors
occur equally often, since each row has even length. The matched
edges crossing that height therefore have equally many upper
endpoints of each color. Summing over $j$ gives
\begin{equation}\label{eq:row-balance}
 v(M)=2\dn(M).
\end{equation}
The transfer matrices use the equivalent weights
\begin{equation}\label{eq:biased-measure}
 \PP\{M\}=\bk^{\dn(M)}/Z_N,
\end{equation}
with weight $\bk$ on down edges and $1$ elsewhere.

\subsection{The liquid ellipse}
\label{subsec:liquid-coordinates}

A lattice vertex $(X,Y)\in V_N$ has rescaled physical coordinate
$z=(X+\ii Y)/N$. In the $z$-plane the limiting diamond is
\begin{equation}\label{eq:limiting-aztec-diamond}
 \Diamond=\{z\in\mathbb C:|\Re z|+|\Im z|\le1\}.
\end{equation}
By the arctic ellipse theorem~\cite[Theorem~26]{CEP},
the liquid region is
\begin{equation}\label{eq:physical-liquid}
 L=\Bigl\{z\in\mathbb C:
 (1+\bk)\bigl((\Re z)^2+(\Im z)^2/\bk\bigr)<1\Bigr\}.
\end{equation}
The arctic ellipse $\partial L$ touches the four sides at
$(\pm(1+\bk)^{-1},\pm\bk(1+\bk)^{-1})$. For uniform tilings
it is the circle of radius $1/\sqrt2$ \cite{JPS}.

Write $\DD=\{w:|w|<1\}$ and $L_1=\{z:|z|<1/\sqrt2\}$. The bijection
$\Psi_1(w)=\sqrt2\,w/(1+|w|^2)$ from $\DD$ to $L_1$ gives the
uniform model's intrinsic coordinate $\Phi_1=\Psi_1^{-1}$~\cite[Eq.~(1.1)]{CR}.
The real-linear map
\begin{equation}\label{eq:linear-map}
 \mathsf A_\bk(z)
    =\sqrt{\frac{2}{1+\bk}}\,\bigl(\Re z+\ii\sqrt\bk\,\Im z\bigr)
\end{equation}
sends $L_1$ onto $L$. Set $\Phi=\Phi_1\circ\mathsf A_\bk^{-1}$ and
$\Psi=\mathsf A_\bk\circ\Psi_1$, so that
\begin{equation}\label{eq:intrinsic-map}
 \begin{aligned}
 \Phi(z)&=\frac{\sqrt{1+\bk}\,(\Re z+\ii\,\Im z/\sqrt\bk)}
   {1+\sqrt{1-(1+\bk)\bigl((\Re z)^2+(\Im z)^2/\bk\bigr)}},&
 \Psi(w)&=\frac{2(\Re w+\ii\sqrt\bk\,\Im w)}
                 {\sqrt{1+\bk}\,(1+|w|^2)} .
 \end{aligned}
\end{equation}
These maps extend to inverse homeomorphisms of $\overline L$ and
$\overline{\DD}$. See~\cite[Eq.~(6.3)]{CJY} and~\cite[Prop.~7.2]{BL}.

\subsection{The curve topology and main theorem}
\label{subsec:curve-topology}

\begin{definition}[Loop collections]\label{def:loop-collection-space}
Let $\mathbb T=\mathbb R/\mathbb Z$. For continuous loops
$\gamma,\eta:\mathbb T\to\Diamond$, define the Fr\'echet pseudometric
\begin{equation}\label{eq:curve-frechet-distance}
 d_{\rm curv}(\gamma,\eta)
 =\inf_{\varphi\in\operatorname{Homeo}(\mathbb T)}
   \sup_{t\in\mathbb T}|\gamma(t)-\eta(\varphi(t))|.
\end{equation}
Here $\varphi$ may preserve or reverse orientation
(cf.\ \cite[Section~2.1]{AB}).
Identifying nonconstant loops at zero $d_{\rm curv}$-distance gives
the metric space $\mathcal C_\Diamond$.

Let $\mathcal X_\Diamond$ consist of the unordered, at most countable
collections $\Gamma$ of loops in $\mathcal C_\Diamond$ satisfying
\begin{equation}\label{eq:loop-collection-local-finiteness}
 \#\{\gamma\in\Gamma:\diam\gamma\ge\varepsilon\}<\infty
 \qquad(\varepsilon>0),
\end{equation}
where loops are counted with multiplicity.

A partial pairing $\pi$ between $\Gamma$ and $\Gamma'$ is a set of
pairs $(\gamma,\gamma')$, one loop copy from each collection, with
no copy used twice. Let $E(\pi)$ be
the least $r\ge0$ such that $d_{\rm curv}(\gamma,\gamma')\le r$
for every $(\gamma,\gamma')\in\pi$ and every unpaired loop in either collection
has diameter at most $2r$.
Equip $\mathcal X_\Diamond$ with the metric
\begin{equation}\label{eq:loop-collection-distance}
 d_{\rm coll}(\Gamma,\Gamma')=\inf_\pi E(\pi)
\end{equation}
and its Borel $\sigma$-algebra.
\end{definition}

\begin{theorem}[Global curve limit]\label{thm:global-curves}
For every fixed $a>0$, as $N\to\infty$,
\begin{equation}\label{eq:global-curve-limit}
 \mathcal L_N
 \overset{\text{in law}}{\Longrightarrow}
 \Psi\bigl(\CLE_4^{\rm nested}(\DD)\bigr)
 \qquad\text{in }(\mathcal X_\Diamond,d_{\rm coll}).
\end{equation}
\end{theorem}

For each fixed diameter cutoff, with arbitrarily high probability all
discrete loops above the cutoff lie a positive distance inside $L$ for
large $N$ (Corollary~\ref{cor:arctic-containment}).

\subsection{Diagonal sections, traversals and tested arcs}
\label{subsec:seams}

\paragraph{White sections and states.}
For $1\le r\le N$, the \emph{white section} $r$ consists of the white
vertices
\begin{equation}\label{eq:white-sections}
 w_{r,x}=\bigl(r-x-\tfrac12,\ r+x-N-\tfrac12\bigr),
 \qquad x\in E_N=\{0,1,\ldots,N\}.
\end{equation}

\begin{definition}[Section states]\label{def:section-states}
For a perfect matching $M$ and $1\le r\le N$, the \emph{state} of
$M$ at section $r$ is the set $S_r(M)\subset E_N$ of the indices $x$
for which $w_{r,x}$ is matched in $M$ to the black vertex immediately
to its right or immediately above it. We write $S_r$ when $M$ is
understood. The \emph{particle number} of a set $S\subseteq E_N$ is
its cardinality $|S|$.
\end{definition}

Counting the vertices on the lower left side of the section shows that
$|S_r(M)|=r$ for every $M$. At the formal boundary sections $0$ and $N+1$, set
$S_0=\varnothing$ and $S_{N+1}=E_N$.

\paragraph{Diagonal coordinates.}
For each $N$, define the \emph{diagonal coordinates} $(s,y)$ of a
point $z$ in the rescaled physical plane by
\begin{equation}\label{eq:seam-physical-coordinates}
 z=(s-y)+\ii(s+y-1)-\frac{1+\ii}{2N}.
\end{equation}
Then $N^{-1}w_{r,x}$ has coordinates $(s,y)=(r/N,x/N)$. Thus the white
section $r$ lies on the line $s=r/N$, the index $x$ gives the
\emph{transverse coordinate} $y=x/N$, and both coordinates change by
$1/(2N)$ along every edge of the scaled graph. The last term
in~\eqref{eq:seam-physical-coordinates} is a translation that tends to
zero. In the limit the diagonal coordinates and the physical
coordinate are related by
\[
  z=(s-y)+\ii(s+y-1),
\]
and this affine map sends
\begin{equation}\label{eq:liquid}
 \Omega=\bigl\{(s,y)\in(0,1)^2:
   \bigl((1+\bk)y-(1-s)-\bk s\bigr)^2<4\bk s(1-s)\bigr\}
\end{equation}
onto $L$. Indeed, substituting $\Re z=s-y$ and $\Im z=s+y-1$
in~\eqref{eq:physical-liquid} gives the inequality defining $\Omega$. The set $\Omega$ is an open
elliptical region, symmetric under $(s,y)\mapsto(y,s)$; for
$\bk=1$ it is the disk $(2s-1)^2+(2y-1)^2<1$.

\paragraph{Traversals.}
For $0<u<1$ let $r_u\in\{1,\ldots,N\}$ be the index for which $r_u/N$
is nearest to $u$, the smaller one in case of a tie.

\begin{definition}[Traversals]\label{def:traversals}
Let $1\le r<r'\le N$. A \emph{traversal} between the white sections
$r$ and $r'$ is a subarc of a loop of the superposition
$M_N^{(1)}\cup M_N^{(2)}$ that has one endpoint on each of the two
sections and whose interior lies strictly between the two section
lines. Different subarcs of the same loop are counted separately. For
$0<u<v<1$, let $\NN_N(u,v)$ be the number of traversals between the
sections $r_u$ and $r_v$, and put $\NN_N(u,v)=0$ if $r_u=r_v$.
\end{definition}

\paragraph{Tested arcs.}
\begin{definition}[Deep regions, tested arcs, index]
\label{def:tested-arcs}
Fix $c\in(0,1)$ and positive numbers $\rho_-,\rho_+$ with
$0<c-\rho_-<c+\rho_+<1$.
\begin{enumerate}[label=\textnormal{(\roman*)},leftmargin=*]
\item The \emph{left} and \emph{right deep regions} are
$\{s\le c-\rho_-\}$ and $\{s\ge c+\rho_+\}$.
\item Choose an orientation of every loop. A \emph{tested arc} is a
subarc of a loop, directed along the loop, with one endpoint in each
deep region.
\item Let $I$ be an interval of the transverse coordinate,
identified with the segment $\{(c',y):y\in I\}$ of a line
$s=c'$ that lies between the two deep regions. Assume that the
endpoints of the segment are not on the graph. For an oriented arc
$\alpha$, the \emph{index} $\iota_I(\alpha)$ is the number of
crossings of $\alpha$ through the segment from $s<c'$ to $s>c'$ minus
the number of crossings in the reverse direction.
\end{enumerate}
\end{definition}

A loop $g$ is a simple closed curve, so $\iota_I(g)\in\{-1,0,1\}$: the
index is the difference of the winding numbers of $g$ around the two
endpoints of the segment (Lemma~\ref{lem:winding-index}).

\subsection{Reflection of sections}
\label{subsec:reflection}

Reflection across the vertical axis sends the lattice vertex
$(X,Y)$ to $(-X,Y)$. In the rescaled physical plane it acts by
\begin{equation}\label{eq:reflection}
 \varrho(z)=-\overline z.
\end{equation}
It preserves the horizontal and vertical domino weights and exchanges the vertex colors.
If a vertex has diagonal coordinates $(s,y)$, then its reflection
across the imaginary axis has diagonal coordinates
$(y+1/(2N),\,s-1/(2N))$ by~\eqref{eq:seam-physical-coordinates}.
Reflection takes the white sections to sections in the other diagonal
direction and preserves crossing probabilities for every $a>0$.

\subsection{Proof structure}
\label{subsec:outline}

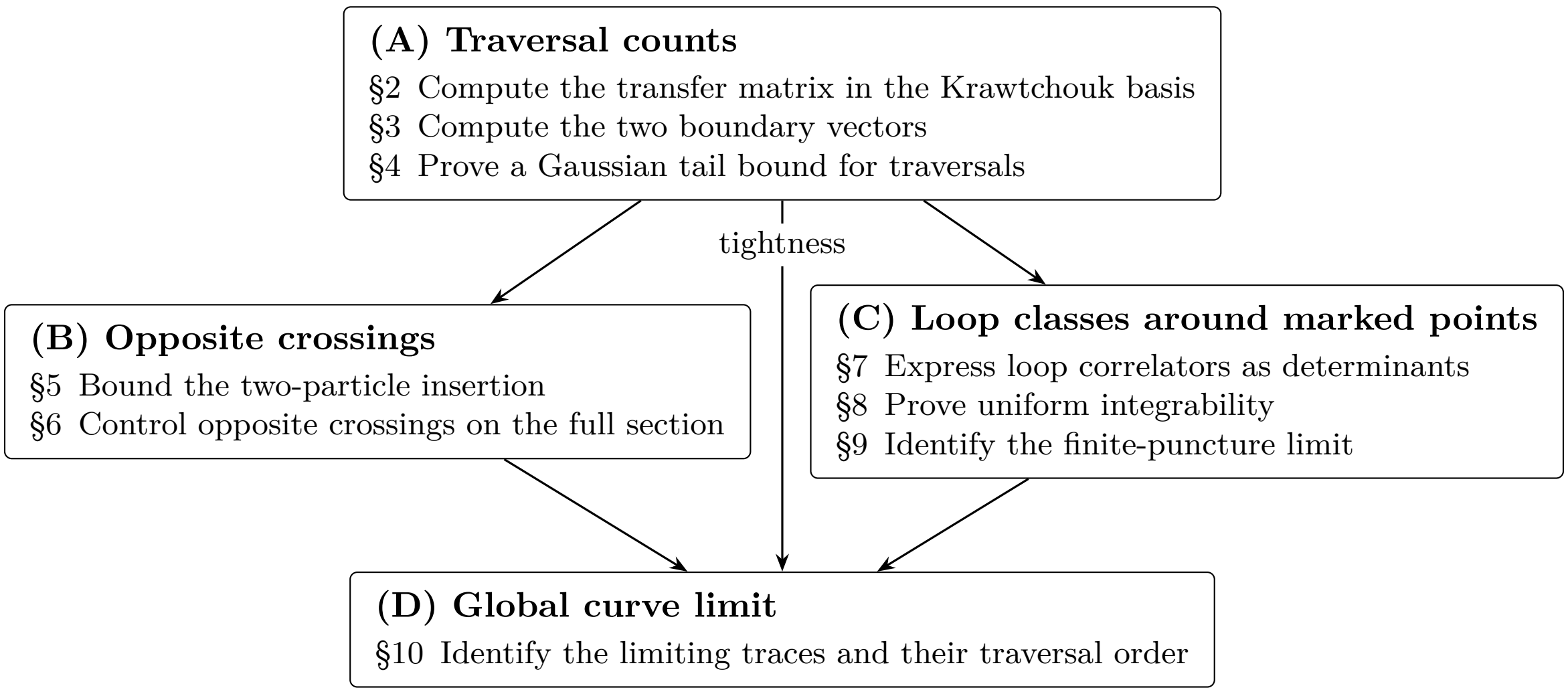
\begin{figure}[htbp]
  \centering
  \begin{tikzpicture}[
    stage/.style={draw, rounded corners=2pt, align=left, font=\footnotesize,
                  inner xsep=7pt, inner ysep=5pt},
    arrow/.style={-{Stealth[length=5pt]}, line width=0.6pt}]
    \node[stage] (A) at (0,2.75) {
      {\small\bfseries (A) Traversal counts}\\[2pt]
      \S\ref{sec:transfer}\enspace Compute the transfer matrix in the Krawtchouk basis\\
      \S\ref{sec:caps}\enspace Compute the two boundary vectors\\
      \S\ref{sec:traversals}\enspace Prove a Gaussian tail bound for traversals};

    \node[stage] (B) at (-4.0,0) {
      {\small\bfseries (B) Opposite crossings}\\[2pt]
      \S\ref{sec:insertion}\enspace Bound the two-particle insertion\\
      \S\ref{sec:global-opposite}\enspace Control opposite crossings on the full section};

    \node[stage] (C) at (4.0,0) {
      {\small\bfseries (C) Loop classes around marked points}\\[2pt]
      \S\ref{sec:punctures}\enspace Express loop correlators as determinants\\
      \S\ref{sec:word-length}\enspace Prove uniform integrability\\
      \S\ref{sec:interior-boundary}\enspace Identify the finite-puncture limit};

    \node[stage] (D) at (0,-2.45) {
      {\small\bfseries (D) Global curve limit}\\[2pt]
      \S\ref{sec:global-curves}\enspace Identify the limiting traces and their traversal order};

    \draw[arrow] (A) -- (B);
    \draw[arrow] (A) -- (C);
    \draw[arrow] (A) -- node[pos=0.12, fill=white, inner sep=2pt,
                              font=\footnotesize] {tightness} (D);
    \draw[arrow] (B) -- (D);
    \draw[arrow] (C) -- (D);
  \end{tikzpicture}
  \caption{Proof of Theorem~\ref{thm:global-curves}.}
  \label{fig:proof-structure}
\end{figure}

Figure~\ref{fig:proof-structure} outlines the proof of
Theorem~\ref{thm:global-curves}.

\paragraph{(A) A bound on traversals
(Sections~\ref{sec:transfer}--\ref{sec:traversals}).}
For two sections a fixed positive distance apart in the rescaled
diamond, the probability of at least $2k$ traversals between them is
at most $Ce^{-ck^2}$, uniformly in $N$
(Proposition~\ref{prop:gaussian-traversal-tail}).
The bound counts traversals along the full sections, including their
frozen parts. To prove it, we count the ways to assign the edges of each crossing path
alternately to the two matchings in the slab and use the transfer matrix in the
Krawtchouk basis to bound the resulting sum.

\paragraph{(B) Opposite crossings
(Sections~\ref{sec:insertion}--\ref{sec:global-opposite}).}
Partition the entire section into intervals of length at most $h$.
The probability that two tested arcs of the same loop have indices of
opposite signs on some interval tends to zero as $N\to\infty$ followed
by $h\downarrow0$, uniformly over the partitions
(Theorem~\ref{thm:global-opposite}).
The proof uses planarity to construct a nonnegative quantity for each
interval that is positive whenever such a pair occurs there. We bound
the sum of its expectations in the Krawtchouk basis and apply the
traversal estimate from~\textnormal{(A)}.

\paragraph{(C) The finite-puncture limit
(Sections~\ref{sec:punctures}--\ref{sec:interior-boundary}).}
Fix finitely many marked points in $L$ and record, with multiplicity,
the homotopy classes of the loops enclosing at least two of them.
The law of this collection converges in total variation to the
corresponding law for nested $\CLE_4$ in $\DD$, with each marked
point $z$ mapped to $\Phi(z)$
(Theorem~\ref{thm:aztec-finite-puncture}).
We compute the limiting loop correlators using Kenyon's determinant
identity~\cite{Kenyon} and an exact formula for the inverse Kasteleyn
matrix. We identify the limit by comparison with square-lattice
double dimers in the upper half-plane
(Lemma~\ref{lem:canonical-reference-identification}).

\paragraph{(D) From homotopy classes to curves
(Section~\ref{sec:global-curves}).}
The traversal bound in~(A) gives tightness of $\mathcal L_N$.
For any subsequential limit, (C) pairs its loops with those of
$\Psi(\CLE_4^{\rm nested}(\DD))$ by the points they surround in a
fixed countable dense set. By~(B), paired loops have the same trace.
A winding-number argument shows that backtracking would force two discrete arcs to have opposite
crossing indices on arbitrarily short intervals, contradicting~(B)
(Lemma~\ref{lem:no-backtracking}).

\paragraph{Relation to the half-plane model.}
Our proof builds on the crossing strategy of~\cite{HP}.
Three features of the Aztec diamond argument are worth distinguishing.

\begin{enumerate}[leftmargin=*]
\item \emph{Crossing estimates on full sections.}
For every $a>0$, we compute the transfer matrix between two sections
and its singular values in the Krawtchouk basis. We use these formulas
to prove a Gaussian tail for the number of traversals between the
sections and to bound the probability of opposite crossings on short
intervals anywhere along a section, including where it meets the
arctic ellipse (Proposition~\ref{prop:gaussian-traversal-tail} and
Theorem~\ref{thm:global-opposite}).

\item \emph{Exact counting identities.}
We average over the matching assignments at the endpoints of every
arch that reaches an outer section. This gives an exact count for
transfers between particle-number sectors
(Lemma~\ref{lem:sector-count}) and expresses a quadratic function
of signed crossings as a matrix element in expectation
(Lemma~\ref{lem:operator-identity}).

\item \emph{The finite-puncture law.}
We compute the limiting correlators from the inverse Kasteleyn
matrix of the biased Aztec diamond and identify them by comparison
with the half-plane model. For any fixed set of marked points in $L$,
the resulting law of reduced loop classes converges in total variation
to that of nested $\CLE_4$ at their images under $\Phi$
(Theorem~\ref{thm:aztec-finite-puncture}).
\end{enumerate}

\subsection{Notation}
\label{subsec:notation}

Constants may depend on the fixed bias $\bk=a^2$ and are uniform in $N$
unless stated otherwise. Hats refer to the gauge in~\eqref{eq:cap-gauge}.
From Section~\ref{sec:insertion} onward, $a,b$ also occur as labels
for section levels (Table~\ref{tab:notation}).

\begin{table}[htbp]
\centering
\small
\renewcommand{\arraystretch}{1.12}
\begin{tabular}{@{}p{0.25\linewidth}p{0.50\linewidth}p{0.18\linewidth}@{}}
\hline
Symbol & Meaning & Defined in\\
\hline
$\bk=a^2$ & bias parameter; $a$ is the weight of a vertical domino
  & Sec.~\ref{subsec:aztec-model}\\
$\mathcal L_N$ & scaled loop collection
  & Sec.~\ref{subsec:aztec-model}\\
$(X,Y)$, $z$ & lattice coordinates and rescaled position $z=(X+\ii Y)/N$
  & Sec.~\ref{subsec:liquid-coordinates}\\
$L$, $\Phi$, $\Psi$ & liquid region and its coordinate maps
  & Sec.~\ref{subsec:liquid-coordinates}\\
$(s,y)$, $E_N$, $w_{r,x}$ & diagonal coordinates, $\{0,\ldots,N\}$,
  white vertices & Sec.~\ref{subsec:seams}\\
$S_r(M)$ & state of $M$ at section $r$; $|S_r|=r$
  & Def.~\ref{def:section-states}\\
$\NN_N(u,v)$, $P_N^{\rm tr}(k;r,r')$ & number of traversals; probability
  of at least $2k$ traversals
  & Def.~\ref{def:traversals}, Sec.~\ref{subsec:deep-arches}\\
$\iota_I(\alpha)$ & index of an arc $\alpha$ on an interval $I$
  & Def.~\ref{def:tested-arcs}\\
$\psi_j$, $e_A$, $\varepsilon_S$ & Krawtchouk basis of
  $V=\mathbb R^{N+1}$ and bases of $\Lambda^pV$
  & Sec.~\ref{subsec:krawtchouk-basis}, \ref{subsec:exact-transfer}\\
$\widehat T_t^{(p)}$, $q_t(j)$ & transfer through $t$ layers on
  $\Lambda^pV$; its singular factors
  & Sec.~\ref{subsec:exact-transfer}\\
$\widehat L_r$, $\widehat R_r$, $\Omega_r$, $G_{r,r'}$ & cap vectors,
  their common direction, transfer factor
  & Sec.~\ref{subsec:cap-partition}\\
$\mathcal H_r$, $\widehat F_\tau(\Theta;\mathcal C)$, $J$, $P$
  & two-matching space, colour tensors, endpoint matrices
  & Sec.~\ref{subsec:marked-caps}\\
$n(\omega)$, $d(\omega)$, $m(\omega)$ & numbers of active loops, deep
  arches, crossing connecting paths
  & Def.~\ref{def:deep-arches}\\
$r_-<r_a<r_c<r_b<r_+$ & outer sections and the three cuts
  & Sec.~\ref{sec:insertion}, \ref{subsec:opposite-setting}\\
$\mathcal O_I$, $\widehat{\mathcal T}_\pm$, $G_\pm$ & two-particle
  insertion, transfers to and from the middle cut, their transfer
  factors
  & Sec.~\ref{sec:insertion}\\
$\sigma(x)$, $t_\eta$, $t_{\rm tot}$, $X_I$ & signs on the middle cut
  and the function built from them
  & Sec.~\ref{subsec:XI}\\
$B_{N,I}$, $\mathcal F_M$ & opposite crossings on $I$; at most $M$
  deep arches per cap & Sec.~\ref{subsec:opposite-setting}\\
$F_N(\rho)$, $X_{\rm unip}$ & topological correlators, unipotent
  representations & Sec.~\ref{sec:punctures}\\
$\mathcal Z$, $z_N$, $\chi_{N,j}(z)$, $S_j$ & probes, their
  realizations, indicators, signatures
  & Sec.~\ref{subsec:probes}, \ref{subsec:signatures}\\
\hline
\end{tabular}
\caption{Notation used in more than one section.}
\label{tab:notation}
\end{table}

\section{Exact transfer in the Krawtchouk basis}
\label{sec:transfer}

We fix the states on two consecutive white sections and express the
total weight of the matchings between them as a determinant. These
weights form the transfer matrix for one step. After a diagonal
rescaling, this matrix acts as a weighted shift in the Krawtchouk
basis. Iterating the shift gives the transfer matrix for $t$ steps
and its singular values (Proposition~\ref{prop:exact-transfer}).

\subsection{Single-layer matchings}
\label{subsec:layer}

The black vertices between white sections $r$ and $r+1$ are
\[
 b_{r,i}=\bigl(r-i-\tfrac12,\ r+i-N+\tfrac12\bigr),
 \qquad 0\le r\le N,\quad 0\le i<N.
\]
For $1\le r<N$, the neighbors of $b_{r,i}$ are
$w_{r,i},w_{r,i+1},w_{r+1,i},w_{r+1,i+1}$.
The neighbors of $b_{0,i}$ all lie on white section $1$, and those
of $b_{N,i}$ all lie on white section $N$.
For $0\le r\le N$, define \emph{layer $r$} to consist of the
black vertices $b_{r,i}$, $0\le i<N$, their white neighbors,
and the edges incident to those black vertices.
Its left and right white sections are $r$ and $r+1$, respectively,
when present. For each black vertex, its left and lower neighbors
belong to the left section, and its right and upper neighbors to
the right section.

Write $S,T$ for the input and output states of a layer, so that
$|T|=|S|+1$. A matching in the layer with these boundary states covers every black
vertex of the layer. The white vertices it covers are precisely those
indexed by $S$ on the left and by $T^c=E_N\setminus T$ on the right. Let $U(T,S)$ be the sum of the
weights of these matchings. As in Section~\ref{subsec:aztec-model},
each down edge $b_{r,i}w_{r,i}$ has weight $\bk$ and all other edges
have weight $1$. At the boundary layers,
$S=\varnothing$ at $r=0$ and $T=E_N$ at $r=N$ leave the absent
section empty.

Define the matrix $H$, with rows and columns indexed by $E_N$, by
\begin{equation}\label{eq:H-definition}
 H(y,x)=\one_{\{y=x\}}+(1+\bk)\one_{\{y>x\}},
 \qquad x,y\in E_N.
\end{equation}

We compute the single-layer weight $U(T,S)$ directly.
The determinant formula~\eqref{eq:layer} also follows from Borodin
and Ferrari~\cite[Theorem~3.8]{BorodinFerrari} by summing over the
particle positions between the two sections and multiplying by
$(1+\bk)^{|S|}$. The product formula~\eqref{eq:layer-weight} agrees
with Tokuyama's weight~\cite[Eq.~(2.8)]{Tokuyama} for the two decreasing
rows $(N-t_i)_{i=1}^{p+1}$ and $(N-s_i)_{i=1}^p$, with $z_n=1$ and $t=\bk$.

\begin{lemma}\label{lem:layer-matching}
Let $S=\{s_1<\cdots<s_p\}$ and
$T=\{t_1<\cdots<t_{p+1}\}$ be subsets of $E_N$, with $0\le p\le N$.
Then
\begin{equation}\label{eq:layer}
 U(T,S)=\det
 \begin{pmatrix}
  1&H(t_1,s_1)&\cdots&H(t_1,s_p)\\
  \vdots&\vdots&&\vdots\\
  1&H(t_{p+1},s_1)&\cdots&H(t_{p+1},s_p)
 \end{pmatrix}.
\end{equation}
This number equals
\begin{equation}\label{eq:layer-weight}
 (1+\bk)^{|S\setminus T|}\,\bk^{\,e(T,S)},
 \qquad e(T,S)=\#\{1\le i\le p:\ s_i=t_i\},
\end{equation}
if
\begin{equation}\label{eq:interlacing}
 t_1\le s_1\le t_2\le\cdots\le s_p\le t_{p+1},
\end{equation}
and equals zero otherwise. For $p=0$, the interlacing condition
is understood to be vacuous.
\end{lemma}

\begin{proof}
For each $x\in E_N$, let
\[
 n_x=\one_{\{x\in S\}}+\one_{\{x\notin T\}}
\]
be the number of selected white vertices of index $x$ on the two sections. Put
\[
 c_x=\sum_{z=0}^x n_z
 =x+1+|S\cap\{0,\ldots,x\}|-|T\cap\{0,\ldots,x\}|.
\]
In any such matching, the black vertices $b_{r,0},\ldots,b_{r,x-1}$
must use white indices at most $x$, while such white vertices
can only use $b_{r,0},\ldots,b_{r,x}$. Thus
$x\le c_x\le x+1$ for $0\le x<N$.
These inequalities are equivalent to~\eqref{eq:interlacing}.
Conversely, assume these inequalities and assign $b_{r,x}$ to white
index $x$ if $c_x=x+1$, and index $x+1$ if $c_x=x$.
These choices of index are forced by the prefix counts.
Taking differences of consecutive prefix counts shows that these
assignments supply exactly $n_x$ black neighbors at each index.
When $n_x=1$, the edge is forced. It is a down edge exactly when
the selected white vertex is $w_{r,x}$ and its assigned black neighbor
is $b_{r,x}$, that is, when $x\in S\cap T$ and $c_x=x+1$.
For $x=s_i\in S\cap T$ we have
$c_x=x+1+i-|T\cap\{0,\ldots,x\}|$, and interlacing gives
$|T\cap\{0,\ldots,x\}|=i$ if $s_i=t_i$ and $i+1$ if $s_i=t_{i+1}$.
Hence a forced edge at $x=s_i$ is a down edge exactly when $s_i=t_i$;
all other forced edges have weight $1$.
When $n_x=2$, equivalently $x\in S\setminus T$, the two assigned
black vertices are $b_{r,x-1}$ and $b_{r,x}$. Of the two possible
matchings, only the one containing $b_{r,x}w_{r,x}$ uses a down edge,
so their weights sum to $1+\bk$. The choices at different indices
are independent, giving~\eqref{eq:layer-weight}.

To evaluate the determinant in~\eqref{eq:layer}, subtract each row
from the next, working from bottom to top, and expand along the
first column. The remaining matrix has entries
\[
 R_{ij}=H(t_{i+1},s_j)-H(t_i,s_j)
 =\begin{cases}
  1+\bk,&t_i<s_j<t_{i+1},\\
  \bk,&s_j=t_i,\\
  1,&s_j=t_{i+1},\\
  0,&\text{otherwise}.
 \end{cases}
\]
Only the identity permutation can contribute to its determinant.
Indeed, a nonzero term with an inversion $i<k$ and
$\sigma(i)>\sigma(k)$ would require
\[
 s_{\sigma(i)}\le t_{i+1}\le t_k\le s_{\sigma(k)},
\]
contradicting the ordering of the $s_j$.
The determinant is therefore $\prod_{i=1}^pR_{ii}$, which is
exactly the weight~\eqref{eq:layer-weight}. For $p=0$, both values are $1$.
\end{proof}

\subsection{The Krawtchouk basis}
\label{subsec:krawtchouk-basis}

Define the Krawtchouk coefficients $K_j$ and normalized functions
$\psi_j$ by
\begin{align}
 K_j(x)&=[z^j](\bk+z)^{N-x}(1-z)^x,\label{eq:krawtchouk}\\
 \psi_j(x)&=(-1)^j
       \sqrt{\frac{\binom Nx\bk^{x}}
                  {(1+\bk)^N\binom Nj\bk^{N-j}}}\,K_j(x),
       \qquad 0\le j,x\le N.\label{eq:psi}
\end{align}
Dividing the generating function by $\bk^{N-x}$ and setting
$z=\bk t$ gives $(1+t)^{N-x}(1-\bk t)^x$, the Krawtchouk generating
function with parameter $1/(1+\bk)$~\cite{Krawtchouk1929}. At $\bk=1$ it is symmetric.
Johansson~\cite[Section~2]{Johansson2005} used the extended
Krawtchouk kernel for the particle process of Aztec tilings.
The superscript ${}^{\mathsf T}$ denotes matrix transpose.
Set
\[
 V=\mathbb R^{N+1},\qquad
 D=\operatorname{diag}\left(\sqrt{\binom Nx\bk^{-x}}\right)_{x=0}^N,
 \qquad \one=(1,\ldots,1)^{\mathsf T}\in V,
\]
and, for $0\le j\le N$, define
\[
 a_j=\sqrt{\bk}\,\sqrt{\frac{N-j}{j+1}},
 \qquad
 b_j=(-1)^j\bk^{-j/2}\frac{\binom{N+1}{j+1}}{\sqrt{\binom Nj}}.
\]

\begin{lemma}\label{lem:krawtchouk-action}
The vectors $\psi_0,\ldots,\psi_N$ form a real orthonormal basis
of $V$. Moreover,
\begin{equation}\label{eq:gauged-one}
 \begin{aligned}
 DHD^{-1}\psi_j
   &=b_j\psi_0+a_j\psi_{j+1},
   &&0\le j<N,\\
 DHD^{-1}\psi_N&=(-1)^N\bk^{-N/2}\psi_0,\\
 D\one&=(1+\bk^{-1})^{N/2}\psi_0.
 \end{aligned}
\end{equation}
\end{lemma}

\begin{proof}
Write
\[
 F_x(z)=(\bk+z)^{N-x}(1-z)^x=\sum_{j=0}^N K_j(x)z^j.
\]
The binomial theorem gives
\[
 \begin{aligned}
 \sum_{x=0}^N\binom Nx\bk^x F_x(z)F_x(w)
 &=\bigl[(\bk+z)(\bk+w)+\bk(1-z)(1-w)\bigr]^N\\
 &=(1+\bk)^N(\bk+zw)^N.
 \end{aligned}
\]
Comparing coefficients of $z^jw^k$ yields
\[
 \sum_{x=0}^N\binom Nx\bk^x K_j(x)K_k(x)
 =(1+\bk)^N\binom Nj\bk^{N-j}\one_{\{j=k\}},
 \qquad 0\le j,k\le N.
\]
The normalization in~\eqref{eq:psi} therefore gives
$\langle\psi_j,\psi_k\rangle=\one_{\{j=k\}}$.
These $N+1$ vectors form an orthonormal basis of $V$.

To compute the action of $H$, use
\[
 (1+\bk)z\,\bk^xF_x(z)
 =(\bk+z)\bigl(\bk^xF_x(z)-\bk^{x+1}F_{x+1}(z)\bigr),
 \qquad 0\le x<N,
\]
which holds because $(\bk+z)-\bk(1-z)=(1+\bk)z$.
Summing over $x=0,\ldots,y-1$ gives, for $y\in E_N$,
\[
 \begin{aligned}
 z\left(\bk^yF_y(z)+(1+\bk)\sum_{x=0}^{y-1}\bk^xF_x(z)\right)
 &=z\bk^yF_y(z)+(\bk+z)\bigl(F_0(z)-\bk^yF_y(z)\bigr)\\
 &=(\bk+z)^{N+1}-\bk^{y+1}F_y(z).
 \end{aligned}
\]
Comparing coefficients of $z^{j+1}$ yields
\begin{equation}\label{eq:polynomial-sum}
 \bk^yK_j(y)+(1+\bk)\sum_{x=0}^{y-1}\bk^xK_j(x)
 =\binom{N+1}{j+1}\bk^{N-j}-\bk^{y+1}K_{j+1}(y),
 \qquad 0\le j\le N.
\end{equation}
Here $K_{N+1}(y)=0$ for $y\in E_N$, since $F_y$ has degree $N$.
By~\eqref{eq:psi} and the definition of $D$,
\[
 (D^{-1}\psi_j)(x)
 =\frac{(-1)^j\bk^xK_j(x)}
        {\sqrt{(1+\bk)^N\binom Nj\bk^{N-j}}},
\]
and it follows that
\[
 (DHD^{-1}\psi_j)(y)
 =\frac{(-1)^j\sqrt{\binom Ny\bk^{-y}}}
       {\sqrt{(1+\bk)^N\binom Nj\bk^{N-j}}}
   \left[\binom{N+1}{j+1}\bk^{N-j}-\bk^{y+1}K_{j+1}(y)\right].
\]
Since $K_0(y)=\bk^{N-y}$, we have
$\psi_0(y)=(1+\bk^{-1})^{-N/2}\sqrt{\binom Ny\bk^{-y}}$, which is
the last identity in~\eqref{eq:gauged-one}.
Substituting this and the definitions of $a_j,b_j,\psi_j$ gives
the first two identities in~\eqref{eq:gauged-one}.
\end{proof}

\subsection{Exact transfer and singular values}
\label{subsec:exact-transfer}

We now express the layer transfer in exterior powers of $V$.
Let $\varepsilon_0,\ldots,\varepsilon_N$ be the standard coordinate
basis. For $0\le p\le N+1$, let $\Lambda^pV$ denote its $p$th
exterior power, and set
\[
 \mathcal S_p=\{S\subseteq E_N:|S|=p\}.
\]
For $S=\{s_1<\cdots<s_p\}\in\mathcal S_p$ and
$A=\{j_1<\cdots<j_p\}\in\mathcal S_p$, write
\[
 \varepsilon_S=\varepsilon_{s_1}\wedge\cdots\wedge\varepsilon_{s_p},
 \qquad
 e_A=\psi_{j_1}\wedge\cdots\wedge\psi_{j_p}.
\]
We identify $\Lambda^pV$ with $\mathbb R^{\mathcal S_p}$ using
the basis $(\varepsilon_S)_{S\in\mathcal S_p}$, and equip it
with the inner product for which this basis is orthonormal.
Equivalently, for $u_1,\ldots,u_p,v_1,\ldots,v_p\in V$,
\[
 \left\langle
 u_1\wedge\cdots\wedge u_p,\,
 v_1\wedge\cdots\wedge v_p
 \right\rangle
 =\det\bigl[\langle u_i,v_j\rangle\bigr]_{i,j=1}^{p}.
\]
The basis $(e_A)_{A\in\mathcal S_p}$ is orthonormal by
Lemma~\ref{lem:krawtchouk-action}. We take $\Lambda^0V=\mathbb R$ and
$\varepsilon_\varnothing=e_\varnothing=1$. All operator norms use
these inner products. The induced gauge $\Lambda^pD$ applies $D$ to
each wedge factor, with $\Lambda^0D=\operatorname{id}_{\mathbb R}$.
For $0\le p\le N$, let
\[
 U^{(p)}:\mathbb R^{\mathcal S_p}\longrightarrow
          \mathbb R^{\mathcal S_{p+1}},
 \qquad (U^{(p)})_{T,S}=U(T,S),
\]
and define
\[
 \widehat U^{(p)}=(\Lambda^{p+1}D)U^{(p)}(\Lambda^pD)^{-1}.
\]
For integers $p,t\ge0$ with $p+t\le N+1$, define
\[
 \widehat T_t^{(p)}
 =\widehat U^{(p+t-1)}\cdots\widehat U^{(p)}
 :\Lambda^pV\longrightarrow\Lambda^{p+t}V,
 \qquad t\ge1,
\]
and set $\widehat T_0^{(p)}=\operatorname{id}_{\Lambda^pV}$.
We omit the superscript $(p)$ from $\widehat U^{(p)}$ and
$\widehat T_t^{(p)}$ when the input degree is clear.

For $0\le t\le N+1$ and $j\in E_N$, set
\begin{equation}\label{eq:q-definition}
 q_t(j)=
 \begin{cases}
  \displaystyle\sqrt{\frac{\binom N{j+t}}{\binom Nj}},&j+t\le N,\\[2mm]
  0,&j+t>N,
 \end{cases}
\end{equation}
and, for an input degree $p\ge0$ with $p+t\le N+1$,
\[
 B_{N,t}^{(p)}
 =(1+\bk^{-1})^{Nt/2}\bk^{pt/2}\prod_{m=0}^{t-2}a_m^{\,t-1-m}
 =(1+\bk^{-1})^{Nt/2}\bk^{pt/2+t(t-1)/4}
   \left(\prod_{k=0}^{t-1}\binom Nk\right)^{1/2}.
\]
For $\bk=1$ this constant does not depend on $p$. In general it
depends on $p$ only through the factor $\bk^{pt/2}$; in particular
\begin{equation}\label{eq:B-pair}
 B_{N,t}^{(p_1)}B_{N,t}^{(p_2)}
 \quad\text{depends on $(p_1,p_2)$ only through $p_1+p_2$.}
\end{equation}
Write $A+t=\{j+t:j\in A\}$.
Empty products equal $1$, and $\{0,\ldots,t-1\}$ is empty for $t=0$.

\begin{proposition}[Exact transfer and singular values]\label{prop:exact-transfer}
For $0\le p\le N$ and $A\in\mathcal S_p$,
\begin{equation}\label{eq:one-transfer}
 \widehat U^{(p)}e_A=
 \begin{cases}
  (1+\bk^{-1})^{N/2}\displaystyle\left(\prod_{j\in A}a_j\right)
  e_{\{0\}\cup(A+1)},&N\notin A,\\[2mm]
  0,&N\in A.
 \end{cases}
\end{equation}
More generally, for $p,t\ge0$ with $p+t\le N+1$ and
$A\in\mathcal S_p$,
\begin{equation}\label{eq:t-transfer}
 \widehat T_t^{(p)}e_A=
 \begin{cases}
  B_{N,t}^{(p)}\displaystyle\left(\prod_{j\in A}q_t(j)\right)
  e_{\{0,\ldots,t-1\}\cup(A+t)},&A+t\subseteq E_N,\\[2mm]
  0,&A+t\not\subseteq E_N.
 \end{cases}
\end{equation}
The nonzero singular values of $\widehat T_t^{(p)}$, counted with
multiplicity, are
\[
 B_{N,t}^{(p)}\prod_{j\in A}q_t(j),
 \qquad A\in\mathcal S_p,\quad A+t\subseteq E_N.
\]
In particular, its operator norm is the largest of these coefficients.
\end{proposition}

\begin{proof}
The determinant formula in Lemma~\ref{lem:layer-matching} gives
\begin{equation}\label{eq:transfer-exterior}
 U^{(p)}(v_1\wedge\cdots\wedge v_p)
 =\one\wedge Hv_1\wedge\cdots\wedge Hv_p,
 \qquad v_1,\ldots,v_p\in V.
\end{equation}
Indeed, for $v_i=\varepsilon_{s_i}$, the coefficient of
$\varepsilon_T$ on the right is the determinant in~\eqref{eq:layer};
the general identity follows by linearity.
Thus, for $A=\{j_1<\cdots<j_p\}$,
\[
 \widehat U^{(p)}e_A
 =D\one\wedge(DHD^{-1}\psi_{j_1})\wedge\cdots\wedge
   (DHD^{-1}\psi_{j_p}).
\]
By Lemma~\ref{lem:krawtchouk-action}, the first factor is
$(1+\bk^{-1})^{N/2}\psi_0$. Antisymmetry kills every term with a second factor $\psi_0$. If $N\notin A$, only the $a_j\psi_{j+1}$ terms remain,
and their indices are already increasing after the initial $0$.
If $N\in A$, the factor
$DHD^{-1}\psi_N=(-1)^N\bk^{-N/2}\psi_0$ makes
the whole exterior product zero. This proves~\eqref{eq:one-transfer}.

At each step, every existing index increases by one and a new
index $0$ is inserted. An initial index $j$ contributes
\[
 \prod_{\ell=0}^{t-1}a_{j+\ell}
 =\bk^{t/2}\sqrt{\frac{\binom N{j+t}}{\binom Nj}}=\bk^{t/2}q_t(j),
 \qquad j+t\le N.
\]
The $p$ initial indices therefore contribute
$\bk^{pt/2}\prod_{j\in A}q_t(j)$.
The index inserted at step $\ell$ contributes
$a_0\cdots a_{t-\ell-1}$ in the remaining steps.
Thus the inserted indices contribute $a_m$ exactly $t-1-m$ times.
Together with $(1+\bk^{-1})^{N/2}$ from each step, this gives the
defining product for $B_{N,t}^{(p)}$. Its binomial expression follows from
\[
 \prod_{m=0}^{t-2}a_m^{\,t-1-m}
 =\prod_{r=1}^{t-1}\left(\prod_{m=0}^{r-1}a_m\right)
 =\bk^{t(t-1)/4}\left(\prod_{r=0}^{t-1}\binom Nr\right)^{1/2}.
\]
If an initial index would exceed $N$, the iteration vanishes
at that step by~\eqref{eq:one-transfer}.
This proves~\eqref{eq:t-transfer}, including $t=0$ by the stated
empty-product conventions.

On the nonzero states, the map
$A\mapsto\{0,\ldots,t-1\}\cup(A+t)$ is injective.
The corresponding output vectors are therefore orthonormal, and
\[
 (\widehat T_t^{(p)})^{\mathsf T}\widehat T_t^{(p)}e_A
 =\left(B_{N,t}^{(p)}\prod_{j\in A}q_t(j)\right)^2e_A.
\]
Taking the positive square roots of the nonzero eigenvalues
gives the asserted singular values.
\end{proof}

\section{Normalization of the caps}
\label{sec:caps}

At section $r$, $0\le r\le N+1$, the left cap consists of layers
$0,\ldots,r-1$ and the right cap of layers $r,\ldots,N$. At an
actual white section both caps include the section vertices. Their
partition functions normalize the one-copy transfer
in~\eqref{eq:normalization}. For two copies, the bound
in Lemma~\ref{lem:aztec-marked-cap} also allows prescribed weights
for the endpoint colours of an arch.

\subsection{Cap partition functions and normalization}
\label{subsec:cap-partition}

For a state $S\in\mathcal S_r$ on an actual section, let $L_r(S)$
and $R_r(S)$ be the partition functions of the compatible left and
right cap completions. A completion covers every vertex of its cap
off the section and exactly the sites $S^c$ (left) or $S$ (right)
on the section. A completion $M$ has weight $\bk^{\dn(M)}$.
At the formal boundaries, use
$L_0(\varnothing)=R_{N+1}(E_N)=1$.
Regard these functions as vectors by setting
\[
 L_r=\sum_{S\in\mathcal S_r}L_r(S)\varepsilon_S,
 \qquad
 R_r=\sum_{S\in\mathcal S_r}R_r(S)\varepsilon_S.
\]

Set
\[
 \Omega_r=\psi_0\wedge\cdots\wedge\psi_{r-1},
 \qquad
 \ell_r=B_{N,r}^{(0)}=(1+\bk^{-1})^{Nr/2}\bk^{r(r-1)/4}
 \left(\prod_{j=0}^{r-1}\binom Nj\right)^{1/2},
\]
where $\Omega_0=\ell_0=1$, and define the gauged cap vectors
\begin{equation}\label{eq:cap-gauge}
 \widehat L_r=(\Lambda^rD)L_r,
 \qquad
 \widehat R_r=(\Lambda^rD)^{-1}R_r.
\end{equation}

\begin{proposition}[Cap partition functions and normalization]
\label{prop:cap-normalization}
The partition function $Z_N=\sum_M\bk^{\dn(M)}$, the sum being over
all perfect matchings of $G_N$, is
\[
 Z_N=(1+\bk)^{N(N+1)/2},
\]
and, for $0\le r\le N+1$, the gauged cap vectors satisfy
\begin{equation}\label{eq:cap-state}
 \widehat L_r=\ell_r\Omega_r,
 \qquad
 \widehat R_r=\frac{Z_N}{\ell_r}\Omega_r.
\end{equation}
All vector norms are Euclidean norms with respect to the
coordinate bases $(\varepsilon_S)_{S\in\mathcal S_r}$.
For integers $0\le r\le r'\le N+1$, put
\begin{equation}\label{eq:cap-transfer-factor}
 G_{r,r'}=\frac{\ell_{r'}}{\ell_r}.
\end{equation}
Then
\[
 \widehat T_{r'-r}^{(r)}\Omega_r=G_{r,r'}\Omega_{r'},
 \qquad
 Z_N=\bigl\langle\widehat R_{r'},
          \widehat T_{r'-r}^{(r)}\widehat L_r\bigr\rangle,
\]
and
\begin{equation}\label{eq:normalization}
 \frac{\|\widehat L_r\|\,
       \|\widehat R_{r'}\|\,G_{r,r'}}{Z_N}=1.
\end{equation}
\end{proposition}

\begin{proof}
For $0\le p\le N$ and $T\in\mathcal S_{p+1}$, classifying
left-cap completions according to their state $Q$ at section $p$
gives
\[
 L_{p+1}(T)
 =\sum_{Q\in\mathcal S_p}U(T,Q)L_p(Q).
\]
The decomposition into a left-cap completion and a matching of the
last layer is unique. Starting from $L_0=1$ therefore gives
$L_r=U^{(r-1)}\cdots U^{(0)}1$.

By Proposition~\ref{prop:exact-transfer},
\[
 \widehat L_r=\widehat T_r^{(0)}1=\ell_r\Omega_r.
\]
At the final section, $L_{N+1}=Z_N\varepsilon_{E_N}$ and
$\Lambda^{N+1}D$ acts by multiplication by
\[
 \det D=\bk^{-N(N+1)/4}\prod_{x=0}^N\sqrt{\binom Nx}.
\]
Taking norms in the preceding identity at $r=N+1$, using
$\|\Omega_{N+1}\|=1$, gives
\[
 (\det D)Z_N=\ell_{N+1}
 =(1+\bk^{-1})^{N(N+1)/2}\bk^{N(N+1)/4}
   \prod_{x=0}^N\sqrt{\binom Nx}
 =(1+\bk)^{N(N+1)/2}\det D,
\]
which proves the stated value of $Z_N$.

For each $S\in\mathcal S_r$, $R_r(S)$ is the
$\varepsilon_{E_N}$ coordinate of
$U^{(N)}\cdots U^{(r)}\varepsilon_S$, since the right cap
consists of layers $r,\ldots,N$ and the final state is $E_N$.
Hence
\[
 R_r
 =\bigl(U^{(N)}\cdots U^{(r)}\bigr)^{\mathsf T}
   \varepsilon_{E_N}.
\]
Since $\Lambda^pD$ is diagonal in the coordinate basis,
the gauge definitions give
\[
 \widehat R_r
 =\bigl(\widehat T_{N+1-r}^{(r)}\bigr)^{\mathsf T}
   (\Lambda^{N+1}D)^{-1}\varepsilon_{E_N}.
\]
By~\eqref{eq:t-transfer}, the transfer from section $r$ to $N+1$ annihilates
every $e_A$ except $\Omega_r$, since survival requires
$A\subseteq\{0,\ldots,r-1\}$ and $|A|=r$. Its transpose thus has
range in $\mathbb R\Omega_r$. Gluing the caps at section $r$ gives
\begin{align}
 Z_N=\sum_{S\in\mathcal S_r}L_r(S)R_r(S)
     =\langle\widehat R_r,\widehat L_r\rangle.\label{znin}
\end{align}
Equation~\eqref{znin} and $\widehat L_r=\ell_r\Omega_r$ now give
$\widehat R_r=(Z_N/\ell_r)\Omega_r$.

Now let $0\le r\le r'\le N+1$.
By the definition of the multilayer transfers,
\begin{equation}\label{eq:cap-composition}
 \widehat T_{r'-r}^{(r)}\widehat L_r
 =\widehat T_{r'-r}^{(r)}\widehat T_r^{(0)}1
 =\widehat T_{r'}^{(0)}1
 =\widehat L_{r'}.
\end{equation}
Division by $\ell_r$ gives the asserted action on $\Omega_r$.
Pairing~\eqref{eq:cap-composition} with $\widehat R_{r'}$ gives $Z_N$
by~\eqref{znin}. Finally, $\|\widehat L_r\|=\ell_r$ and
$\|\widehat R_{r'}\|=Z_N/\ell_{r'}$ give~\eqref{eq:normalization}.
\end{proof}

For two independent perfect matchings, the cap vectors are
$\widehat L_r\otimes\widehat L_r$ and
$\widehat R_{r'}\otimes\widehat R_{r'}$, and the partition
function is $Z_N^2$.
We use the inner product for which the vectors
$\varepsilon_{S_1}\otimes\varepsilon_{S_2}$ are orthonormal.
Squaring~\eqref{eq:normalization} gives
\begin{equation}\label{eq:normalization-squared}
 \frac{
 \|\widehat L_r\otimes\widehat L_r\|\,
 \|\widehat R_{r'}\otimes\widehat R_{r'}\|\,
 G_{r,r'}^{\,2}
 }{Z_N^2}=1.
\end{equation}

\subsection{Endpoint colours at a cut}
\label{subsec:marked-caps}

Cutting a double-dimer configuration at a section leaves arches in
each cap. Their endpoint colours record which perfect matching uses each cut
edge. We shall change the colour rule on selected arches while keeping
the cap norm bounded by its value for two independent perfect matchings.

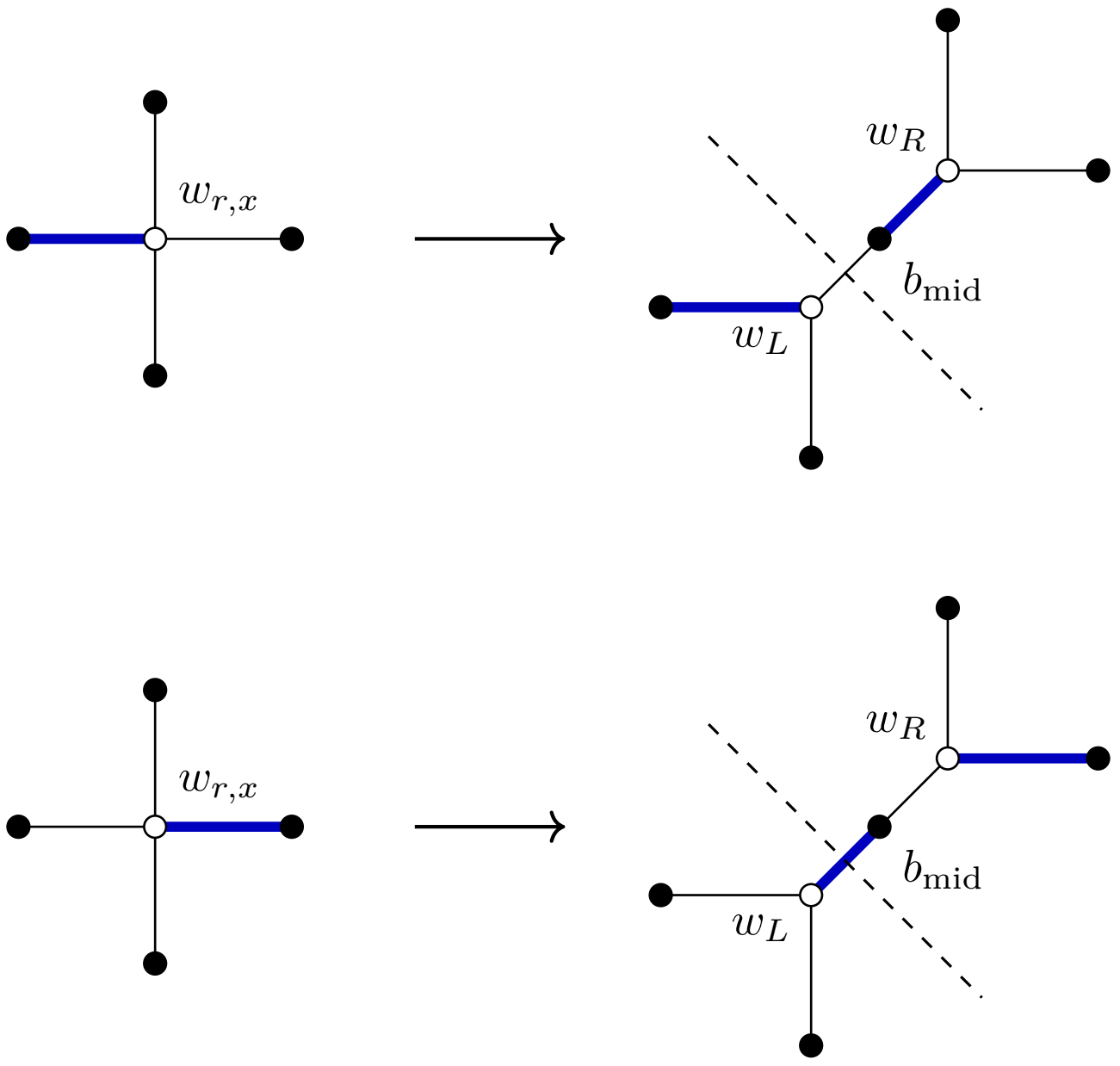
\begin{figure}[htbp]
 \centering
 \begin{tikzpicture}[scale=1.05,
   wv/.style={circle,draw,fill=white,inner sep=1.7pt,line width=0.5pt},
   bv/.style={circle,draw,fill=black,inner sep=1.7pt,line width=0.5pt},
   ed/.style={line width=0.5pt},
   md/.style={line width=2.2pt,blue!75!black,line cap=round}]
 \foreach \yy/\isleft in {0/1,-4.3/0} {
  \begin{scope}[yshift=\yy cm]
   \draw[ed] (-1,0)--(1,0);
   \draw[ed] (0,-1)--(0,1);
   \ifnum\isleft=1 \draw[md] (-1,0)--(0,0);\else \draw[md] (0,0)--(1,0);\fi
   \node[bv] at (-1,0) {};
   \node[bv] at (1,0) {};
   \node[bv] at (0,1) {};
   \node[bv] at (0,-1) {};
   \node[wv,label={[label distance=-1pt]above right:{\small $w_{r,x}$}}] at (0,0) {};
   \draw[->,line width=0.8pt] (1.9,0)--(3.0,0);
   \begin{scope}[xshift=5.3cm]
    \draw[ed] (-1.6,-0.5)--(-0.5,-0.5)--(-0.5,-1.6);
    \draw[ed] (-0.5,-0.5)--(0,0)--(0.5,0.5);
    \draw[ed] (1.6,0.5)--(0.5,0.5)--(0.5,1.6);
    \ifnum\isleft=1
     \draw[md] (-1.6,-0.5)--(-0.5,-0.5);
     \draw[md] (0,0)--(0.5,0.5);
    \else
     \draw[md] (0.5,0.5)--(1.6,0.5);
     \draw[md] (-0.5,-0.5)--(0,0);
    \fi
    \draw[dashed,line width=0.6pt] (-1.25,0.75)--(0.75,-1.25);
    \node[bv] at (-1.6,-0.5) {};
    \node[bv] at (-0.5,-1.6) {};
    \node[bv] at (1.6,0.5) {};
    \node[bv] at (0.5,1.6) {};
    \node[wv,label={[label distance=-2pt]below left:{\small $w_L$}}] at (-0.5,-0.5) {};
    \node[bv,label={[label distance=-1pt]below right:{\small $b_{\rm mid}$}}] at (0,0) {};
    \node[wv,label={[label distance=-2pt]above left:{\small $w_R$}}] at (0.5,0.5) {};
   \end{scope}
  \end{scope}
 }
 \end{tikzpicture}
 \caption{A white vertex before and after splitting. Thick edges
 belong to the perfect matching, and the dashed line is the cut. Matching the
 original vertex on the left (top row) leaves the cut edge empty.
 Matching it on the right (bottom row) occupies the cut edge, so
 $x\in S_r(M)$. The same rule applies to the downward and upward
 neighbors, respectively.}
 \label{fig:vertex-splitting}
\end{figure}

\paragraph{The split graph and the cut.}
Fix a white section $r$, with $1\le r\le N$.
Construct an auxiliary graph by replacing each white vertex
$w_{r,x}$, $x\in E_N$, by a path
$w_L-b_{\rm mid}-w_R$, where $w_L,w_R$ are new white
vertices and $b_{\rm mid}$ is a new black vertex.
Both internal edges have weight $1$. Attach the original edges
extending left or down to $w_L$ and those extending right or up to
$w_R$ (Figure~\ref{fig:vertex-splitting}). A perfect matching extends uniquely
by matching $b_{\rm mid}$ to whichever of $w_L,w_R$ remains free.

If two perfect matchings use the same side of a split vertex, the other
internal edge is doubled. If they use opposite sides, the loop uses
both internal edges. Splitting thus preserves all loops and moves
their traces by $O(N^{-1})$.
In the $(s,y)$ coordinates, make each replacement inside the disk
of radius $1/(4N)$ centered at $(r/N,x/N)$, leaving the embedding
fixed outside these disks. Embed the new edges so that $w_Lb_{\rm mid}$ crosses the line
$s=r/N$ at the point $(r/N,x/N)$, with $w_L$ on the left and
$b_{\rm mid},w_R$ on the right.

We call the line $s=r/N$ the \emph{cut}, the edge $w_Lb_{\rm mid}$ of
the split vertex $w_{r,x}$ the \emph{cut edge at site $x$}, and its
endpoints $w_L$ and $b_{\rm mid}$ the \emph{cap vertices} of the site.
The cut edge at $x$ belongs to the extension of a perfect matching $M$ if and
only if $x\in S_r(M)$. From now on the \emph{left cap} at section $r$
is the part of the split graph to the left of the cut, which contains
the vertices $w_L$, and the \emph{right cap} is the part to its right,
which contains the vertices $b_{\rm mid}$ and $w_R$. In either cap a completion with state $S$ covers every vertex except
the cap vertices at sites of $S$. Splitting gives the completions of
Section~\ref{subsec:cap-partition} with their original weights.

\begin{definition}[Superpositions, arches, colours]
\label{def:cap-superposition}
Fix the left or the right cap at section $r$.
\begin{enumerate}[label=\textnormal{(\roman*)},leftmargin=*]
\item A \emph{superposition} $C$ of the cap is the union, with
multiplicities, of two completions $M_1,M_2$ of the cap, with the
matching labels erased. Its number of down edges,
counted with multiplicity, is $\dn(C)=\dn(M_1)+\dn(M_2)$.
\item The \emph{cut multiplicity} $m_x(C)\in\{0,1,2\}$ of a site
$x\in E_N$ is the number of the two completions whose state contains
$x$. The site is \emph{singly occupied} if $m_x(C)=1$. Write $X(C)$ for
the set of singly occupied sites.
\item In $C$ every vertex has degree two, except that the cap vertex of
a site $x$ has degree $2-m_x(C)$. Hence the edges of multiplicity one
form disjoint cycles, the \emph{internal loops} of $C$, and disjoint
paths joining the cap vertices of two singly occupied sites, the
\emph{arches} of $C$. Write $\ell(C)$ for the number of internal loops
and $\mathcal A(C)$ for the set of arches. Every singly occupied site
is an end of exactly one arch.
\item A pair of states $(Q_1,Q_2)$, $Q_1,Q_2\subseteq E_N$, is
\emph{compatible} with $C$ if
$\one_{\{x\in Q_1\}}+\one_{\{x\in Q_2\}}=m_x(C)$ for all $x$.
A compatible pair is the same thing as a \emph{colouring}
$\xi:X(C)\to\{+,-\}$: put $\xi(x)=+$ if $x\in Q_1$ and $\xi(x)=-$ if
$x\in Q_2$.
\end{enumerate}
\end{definition}

Each original completion has $r$ occupied sites, hence
$\sum_xm_x(C)=2r$. Compatibility preserves $|Q_1|+|Q_2|=2r$ and
allows the two individual cardinalities to vary.

\begin{lemma}[Colours of two completions]\label{lem:ordinary-colours}
Let $C$ be a superposition of a cap and $(Q_1,Q_2)$ a compatible pair
with colouring $\xi$. The number of ordered pairs of completions
$(M_1,M_2)$ with states $(Q_1,Q_2)$ and superposition $C$ is
$2^{\ell(C)}$ if $\xi$ takes different values at the two ends of every
arch, and $0$ otherwise.
\end{lemma}

\begin{proof}
In such a pair the edges of multiplicity one alternate between $M_1$
and $M_2$ along every internal loop and every arch. An arch joins two
cap vertices of the same bipartite colour, both $w_L$ or both
$b_{\rm mid}$. Its length is even, so its first and last edges belong
to different completions. Each endpoint lies in the state of the
completion missing its incident arch edge. The two endpoint colours
therefore differ. Conversely, if $\xi$
takes different values at the two ends of every arch, the colours
determine how each arch is divided between $M_1$ and $M_2$, and each
internal loop can be divided in two ways.
\end{proof}

Let
\[
 \mathcal H_r=
 \bigoplus_{\substack{s\in\mathbb Z\\0\le r+s,\,r-s\le N+1}}
       \bigl(\Lambda^{r+s}V\otimes\Lambda^{r-s}V\bigr),
\]
with the inner product for which the vectors
$\varepsilon_{Q_1}\otimes\varepsilon_{Q_2}$, $|Q_1|+|Q_2|=2r$, are
orthonormal. The summand with index $s$ is the \emph{sector} $s$; in it
the particle numbers of the two states are $r+s$ and $r-s$. Let
$d_x=D_{xx}=\binom Nx^{1/2}\bk^{-x/2}$ be the one-site gauge
of~\eqref{eq:cap-gauge}, and put
\[
 g_L(C)=\prod_{x\in E_N}d_x^{\,m_x(C)},\qquad
 g_R(C)=\prod_{x\in E_N}d_x^{-m_x(C)} .
\]
For every pair $(Q_1,Q_2)$ compatible with $C$, the map
$\Lambda^{|Q_1|}D\otimes\Lambda^{|Q_2|}D$ multiplies
$\varepsilon_{Q_1}\otimes\varepsilon_{Q_2}$ by $g_L(C)$.

\begin{definition}[Endpoint rules and colour tensors]
\label{def:colour-tensor}
An \emph{endpoint rule} $\Theta$ assigns to every superposition $C$ of
a cap and every arch $A\in\mathcal A(C)$ a real symmetric
$2\times2$ matrix $\Theta_A(C)$ with rows and columns indexed by
$\{+,-\}$. Let $\tau\in\{L,R\}$ denote the side of the cap and let
$\mathcal C$ be a set of superpositions of the cap. The
\emph{colour tensor} $\widehat F_\tau(\Theta;\mathcal C)\in\mathcal H_r$
is the vector with coordinates
\begin{equation}\label{eq:auxiliary-cap-vectors}
 \bigl(\widehat F_\tau(\Theta;\mathcal C)\bigr)_{Q_1,Q_2}
 =\sum_{\substack{C\in\mathcal C\\
           (Q_1,Q_2)\ \mathrm{compatible}}}
   g_\tau(C)\,\bk^{\dn(C)}\,2^{\ell(C)}
   \prod_{A\in\mathcal A(C)}
          \bigl(\Theta_A(C)\bigr)_{\xi(x_A),\xi(y_A)} ,
\end{equation}
where $x_A,y_A$ are the two ends of $A$ and $\xi$ is the colouring
given by $(Q_1,Q_2)$. Each superposition occurs once in the sum. We
write $\widehat F_\tau(\Theta)$ when $\mathcal C$ is the set of all
superpositions of the cap.
\end{definition}

The two endpoint rules used below have matrices
\[
 J=\begin{pmatrix}0&1\\1&0\end{pmatrix},\qquad
 P=\frac12\begin{pmatrix}1&1\\1&1\end{pmatrix}.
\]
Both have operator norm $1$. With $\Theta_A=J$ the two ends of $A$ must
have different colours, as for two completions. Taking $\Theta_A=P$ distributes the same total weight $2$
uniformly among all four endpoint pairs. If $\Theta_A(C)=J$ for all $C$ and $A$, then
Lemma~\ref{lem:ordinary-colours} and the definition of the gauge give
\begin{equation}\label{eq:ordinary-colour-tensor}
 \widehat F_L(J)=\widehat L_r\otimes\widehat L_r,\qquad
 \widehat F_R(J)=\widehat R_r\otimes\widehat R_r .
\end{equation}
These vectors lie in the sector $s=0$. A rule that uses $P$ on some
arches produces components in the other sectors as well.

\begin{lemma}[Norm bound for colour tensors]
\label{lem:aztec-marked-cap}
Let $\Theta$ be an endpoint rule with $\|\Theta_A(C)\|\le1$ in the
operator norm for every $C$ and every arch $A$. Then, for every set
$\mathcal C$ of superpositions of the left or right cap at section $r$,
\begin{equation}\label{eq:marked-cap}
 \|\widehat F_L(\Theta;\mathcal C)\|
     \le \|\widehat L_r\otimes\widehat L_r\|,\qquad
 \|\widehat F_R(\Theta;\mathcal C)\|
     \le \|\widehat R_r\otimes\widehat R_r\| .
\end{equation}
\end{lemma}

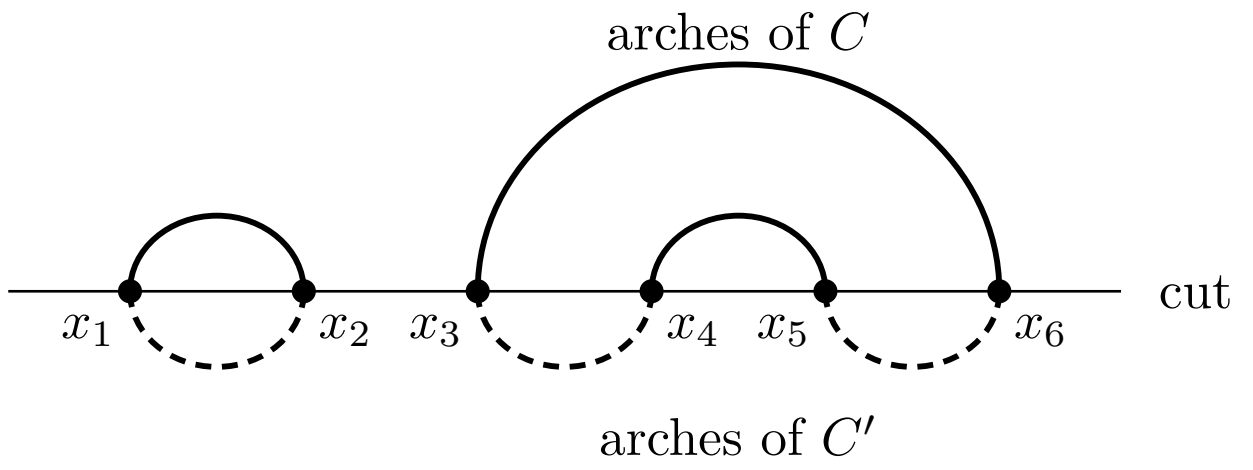
\begin{figure}[htbp]
 \centering
 \begin{tikzpicture}[x=1.15cm,y=1cm,
   st/.style={circle,fill=black,inner sep=1.6pt}]
  \draw[line width=0.5pt] (0.3,0)--(6.7,0);
  \node[anchor=west] at (6.8,0) {\small cut};
  \draw[line width=1.1pt] (1,0) arc (180:0:0.5);
  \draw[line width=1.1pt] (3,0) arc (180:0:1.5);
  \draw[line width=1.1pt] (4,0) arc (180:0:0.5);
  \draw[line width=1.1pt,densely dashed] (1,0) arc (180:360:0.5);
  \draw[line width=1.1pt,densely dashed] (3,0) arc (180:360:0.5);
  \draw[line width=1.1pt,densely dashed] (5,0) arc (180:360:0.5);
  \foreach \i in {1,...,6} {
    \node[st] at (\i,0) {};
  }
  \foreach \i/\dx in {1/-0.24,2/0.24,3/-0.24,4/0.24,5/-0.24,6/0.24} {
    \node at (\i+\dx,-0.26) {\small $x_{\i}$};
  }
  \node at (4.5,1.72) {\small arches of $C$};
  \node at (4.5,-0.95) {\small arches of $C'$};
 \end{tikzpicture}
 \caption{Two superpositions $C$ and $C'$ of the same cap with the
 same singly occupied sites $x_1,\ldots,x_6$. The arches of $C$ (solid)
 and of $C'$ (dashed) are drawn on opposite sides of the cut only for
 legibility. Together they form a cycle through $x_1,x_2$ and a cycle
 through $x_3,x_6,x_5,x_4$. Each cycle contributes the trace of the
 product of its endpoint matrices.}
 \label{fig:arch-cycles}
\end{figure}

\begin{proof}
We treat the left cap; the right cap is the same with $g_R$.
Expanding the square of~\eqref{eq:auxiliary-cap-vectors} and summing
over $(Q_1,Q_2)$ gives
\begin{equation}\label{eq:norm-expansion}
 \|\widehat F_L(\Theta;\mathcal C)\|^2
 =\sum_{C,C'\in\mathcal C}
   g_L(C)g_L(C')\,\bk^{\dn(C)+\dn(C')}\,2^{\ell(C)+\ell(C')}\,
   W(C,C';\Theta),
\end{equation}
where $W(C,C';\Theta)$ is the sum, over the pairs $(Q_1,Q_2)$
compatible with both $C$ and $C'$, of
\[
 \prod_{A\in\mathcal A(C)}\bigl(\Theta_A(C)\bigr)_{\xi(x_A),\xi(y_A)}
 \prod_{A'\in\mathcal A(C')}
     \bigl(\Theta_{A'}(C')\bigr)_{\xi(x_{A'}),\xi(y_{A'})} .
\]
If the cut multiplicities of $C$ and $C'$ differ at some site, no pair
is compatible with both and $W(C,C';\Theta)=0$. Assume they agree.
Then $X(C)=X(C')=:X$, and the pairs compatible with both are the
colourings $\xi\in\{+,-\}^X$.

The arches of $C$ pair up the sites of $X$, and so do the arches of
$C'$. Together the two pairings decompose $X$ into disjoint cycles
along which arches of $C$ and of $C'$ alternate
(Figure~\ref{fig:arch-cycles}); a cycle has two sites when an arch of
$C$ and an arch of $C'$ have the same ends. The sum over $\xi$
factorizes over the cycles. For a cycle through the sites
$x_1,\ldots,x_{2k}$ with arches $A_1,\ldots,A_{2k}$, where $A_i$ joins
$x_i$ to $x_{i+1}$ and $x_{2k+1}=x_1$, the sum over the colours of
these sites is
\[
 \sum_{\xi_1,\ldots,\xi_{2k}\in\{+,-\}}\ \prod_{i=1}^{2k}
   \bigl(\Theta_{A_i}\bigr)_{\xi_i,\xi_{i+1}}
 =\tr\bigl(\Theta_{A_1}\Theta_{A_2}\cdots\Theta_{A_{2k}}\bigr).
\]
Here $\Theta_{A_i}$ stands for $\Theta_{A_i}(C)$ or
$\Theta_{A_i}(C')$, and the symmetry of the matrices is used to write
each factor in the direction of the cycle. Since each matrix has norm at most $1$ and
$|\tr B|\le2\|B\|$ in dimension two, the absolute contribution of
each cycle is at most $2$, attained by the ordinary rule since
$\tr(J^{2k})=2$. Multiplying over the cycles,
\[
 |W(C,C';\Theta)|\le W(C,C';J).
\]
All other factors in~\eqref{eq:norm-expansion} are positive and do not
depend on $\Theta$. Therefore
\[
 \|\widehat F_L(\Theta;\mathcal C)\|^2
 \le\sum_{C,C'\in\mathcal C}
   g_L(C)g_L(C')\,\bk^{\dn(C)+\dn(C')}\,2^{\ell(C)+\ell(C')}\,W(C,C';J)
 \le\|\widehat F_L(J)\|^2 ,
\]
where the last step adds the nonnegative terms with $C$ or $C'$
outside $\mathcal C$. By~\eqref{eq:ordinary-colour-tensor} the right-hand
side is $\|\widehat L_r\otimes\widehat L_r\|^2$.
\end{proof}

\section{Particle-number shifts and traversal counts}
\label{sec:traversals}

A shift from particle numbers $(r,r)$ to $(r+k,r-k)$ costs at most
$e^{-tk^2/(N+1)}$ in the product of the two transfer norms across
$t$ layers, relative to the unshifted product $G_{r,r+t}^2$.
The counting identity~\eqref{eq:sector-count} relates this cost to
the event of $2k$ traversals. Applying it recursively on nested slabs
gives the Gaussian tail~\eqref{eq:gaussian-tail}.

\subsection{The cost of a particle-number shift}
\label{subsec:shift-cost}

\begin{lemma}
\label{lem:transfer-mode-decay}
For integers $1\le t\le N$ and $0\le j<N-t$,
\begin{align}
 \frac{q_t(j+1)}{q_t(j)}
 &=\sqrt{\frac{(N-j-t)(j+1)}
               {(N-j)(j+t+1)}},
 \label{eq:q-ratio}\\
 \log\frac{q_t(j+1)}{q_t(j)}
 &\le-\frac t2
       \left(\frac1{N-j}+\frac1{j+t+1}\right)
 \le-\frac{t}{N+1}.
 \label{eq:q-gap}
\end{align}
In particular, the positive values $q_t(0),\ldots,q_t(N-t)$
are decreasing.
\end{lemma}

\begin{proof}
The definition~\eqref{eq:q-definition} of $q_t$
gives~\eqref{eq:q-ratio}. Taking logarithms and using
$\log(1-z)\le-z$ for $0<z<1$
gives the first inequality in~\eqref{eq:q-gap}. For the second, use
\[
 \frac1{N-j}+\frac1{j+t+1}
 \ge\frac4{N+t+1}
 \ge\frac2{N+1}.
\]
\end{proof}

\begin{proposition}[Particle-number shift bound]
\label{prop:particle-shift-bound}
Let $t$ be an integer with $0\le t\le N+1$.
For every integer $0\le p\le N+1-t$,
\[
 \|\widehat T_t^{(p)}\|
 =B_{N,t}^{(p)}\prod_{j=0}^{p-1}q_t(j).
\]
Let $r,k\ge0$ be integers satisfying
\[
 r+t\le N+1,
 \qquad
 k\le\min\{r,N+1-r-t\}.
\]
Then, with $G_{r,r+t}$ defined
by~\eqref{eq:cap-transfer-factor},
\begin{align}
 \frac{\|\widehat T_t^{(r+k)}\|\,
       \|\widehat T_t^{(r-k)}\|}
      {G_{r,r+t}^{\,2}}
 &=\prod_{j=0}^{k-1}
       \frac{q_t(r+j)}{q_t(r-1-j)}
 \notag\\
 &\le e^{-tk^2/(N+1)}.
 \label{eq:sector-cost}
\end{align}
\end{proposition}

\begin{proof}
The largest singular value in~\eqref{eq:t-transfer} is attained
at $A=\{0,\ldots,p-1\}$ by Lemma~\ref{lem:transfer-mode-decay}.
This also holds when $p=0$ (the empty product) or $t=0$
(all $q_0(j)=1$). Thus $\Omega_p$ realizes the stated norm.
By Proposition~\ref{prop:cap-normalization},
\[
 \|\widehat T_t^{(r)}\|
 =\|\widehat T_t^{(r)}\Omega_r\|
 =G_{r,r+t}.
\]
By~\eqref{eq:B-pair},
$B_{N,t}^{(r+k)}B_{N,t}^{(r-k)}=\bigl(B_{N,t}^{(r)}\bigr)^2$.
The bias factors $\bk^{kt/2}$ and $\bk^{-kt/2}$ cancel in the
product of shifted norms. Cancelling the common $q_t$ factors gives the equality
in~\eqref{eq:sector-cost}. If $t=0$ or $k=0$, the bound
is immediate. Otherwise, summing~\eqref{eq:q-gap} gives
\[
 \log\frac{q_t(r+j)}{q_t(r-1-j)}
 \le-\frac{t(2j+1)}{N+1},
 \qquad 0\le j<k.
\]
Now use $\sum_{j=0}^{k-1}(2j+1)=k^2$.
\end{proof}

\subsection{Deep arches and connecting paths}
\label{subsec:deep-arches}

For white sections $1\le r<r'\le N$ and an integer $k\ge0$, let
$P_N^{\mathrm{tr}}(k;r,r')$ be the probability that there are at least
$2k$ traversals between the sections $r$ and $r'$
(Definition~\ref{def:traversals}). With $r_u,r_v$ and $\NN_N(u,v)$ as
in that definition,
\[
 \PP\{\NN_N(u,v)\ge2k\}
 =P_N^{\mathrm{tr}}(k;r_u,r_v)
 \qquad\text{when }r_u<r_v.
\]

Fix integers
\[
 1\le r_-<r_a<r_b<r_+\le N .
\]
Split the graph at $r_a$ and $r_b$ as in
Section~\ref{subsec:marked-caps}. These are the left and right cuts,
with caps to their left and right, respectively. The sections $r_-$
and $r_+$ are the outer sections. Both $w_L$ and $w_R$ at a split
vertex count as vertices of its section.

A \emph{full configuration} $\omega$ is the superposition of the two
perfect matchings on the split graph, retaining edge multiplicities and
erasing the matching labels. It consists of doubled
edges and disjoint loops, and
\begin{equation}\label{eq:configuration-probability}
 \PP\{\omega\}=Z_N^{-2}\,\bk^{\dn(\omega)}\,2^{\ell(\omega)},
\end{equation}
Here $\dn(\omega)$ counts down edges with multiplicity, and each
of the $\ell(\omega)$ loops has two alternating colourings. The
arches of the restriction to each cap
(Definition~\ref{def:cap-superposition}) are the left and right
arches of $\omega$.

\begin{definition}[Deep arches, active loops, connecting paths]
\label{def:deep-arches}
Let $\omega$ be a full configuration.
\begin{enumerate}[label=\textnormal{(\roman*)},leftmargin=*]
\item A left arch is \emph{deep} if it contains a vertex of the outer
section $r_-$. A right arch is \emph{deep} if it contains a vertex of
the outer section $r_+$.
\item A loop is \emph{active} if it contains a deep arch, and
\emph{inactive} otherwise.
\item Remove from every active loop the edges of its deep arches. The
remaining pieces, called \emph{connecting paths}, begin and end with
cut edges. An active loop with $j$ deep arches has $j$ such paths.
\item A connecting path is \emph{crossing} if its two end cut edges
lie on different cuts, and \emph{returning} if they lie on the same
cut.
\item Write $n(\omega)$ for the number of active loops, $d(\omega)$ for
the number of deep arches, which is also the number of connecting
paths, and $m(\omega)$ for the number of crossing connecting paths.
\end{enumerate}
\end{definition}

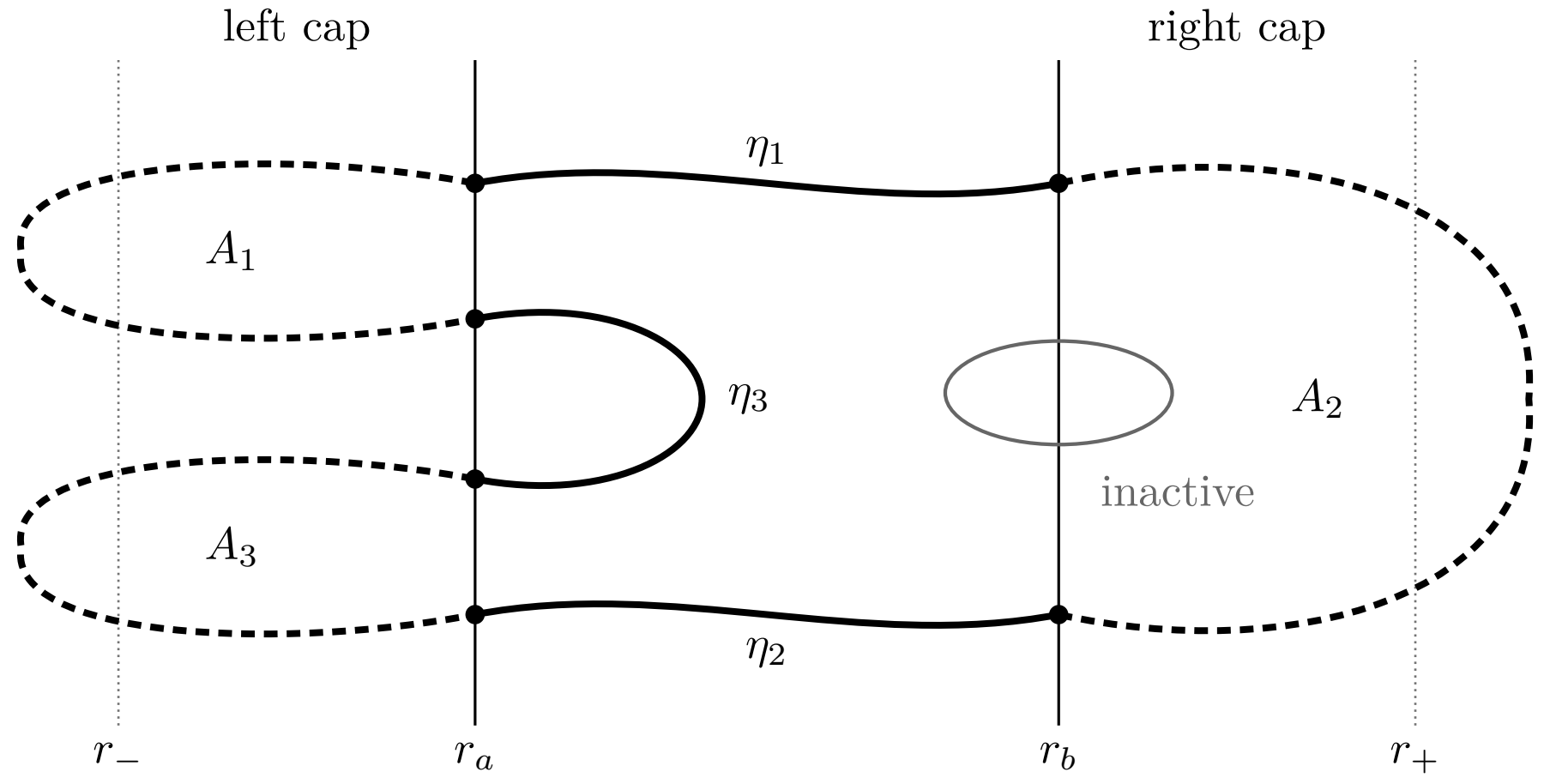
\begin{figure}[htbp]
 \centering
 \begin{tikzpicture}[x=1.25cm,y=0.95cm,
   deep/.style={line width=1.5pt,densely dashed},
   conn/.style={line width=1.5pt},
   st/.style={circle,fill=black,inner sep=1.5pt}]
  \draw[line width=0.5pt,black!55,densely dotted] (0,0.6)--(0,6.0);
  \draw[line width=0.5pt,black!55,densely dotted] (8,0.6)--(8,6.0);
  \draw[line width=0.6pt] (2.2,0.6)--(2.2,6.0);
  \draw[line width=0.6pt] (5.8,0.6)--(5.8,6.0);
  \node[below] at (0,0.6) {\small $r_-$};
  \node[below] at (2.2,0.6) {\small $r_a$};
  \node[below] at (5.8,0.6) {\small $r_b$};
  \node[below] at (8,0.6) {\small $r_+$};
  \node at (1.1,6.25) {\small left cap};
  \node at (6.9,6.25) {\small right cap};
  \draw[deep] (2.2,5.0) .. controls (1.0,5.3) and (-0.7,5.2) .. (-0.6,4.45)
              .. controls (-0.7,3.7) and (1.0,3.6) .. (2.2,3.9);
  \draw[deep] (2.2,2.6) .. controls (1.0,2.9) and (-0.7,2.8) .. (-0.6,2.05)
              .. controls (-0.7,1.3) and (1.0,1.2) .. (2.2,1.5);
  \draw[deep] (5.8,1.5) .. controls (7.2,1.1) and (8.8,1.6) .. (8.7,3.25)
              .. controls (8.8,4.9) and (7.2,5.4) .. (5.8,5.0);
  \draw[conn] (5.8,5.0) .. controls (4.6,4.7) and (3.4,5.3) .. (2.2,5.0);
  \draw[conn] (2.2,3.9) .. controls (3.0,4.1) and (3.6,3.7) .. (3.6,3.25)
              .. controls (3.6,2.8) and (3.0,2.4) .. (2.2,2.6);
  \draw[conn] (2.2,1.5) .. controls (3.4,1.8) and (4.6,1.2) .. (5.8,1.5);
  \draw[line width=0.8pt,black!60] (5.8,3.3) ellipse [x radius=0.7,y radius=0.42];
  \foreach \p in {(2.2,5.0),(2.2,3.9),(2.2,2.6),(2.2,1.5),(5.8,1.5),(5.8,5.0)} {
    \node[st] at \p {};
  }
  \node at (0.7,4.45) {\small $A_1$};
  \node at (0.7,2.05) {\small $A_3$};
  \node at (7.4,3.25) {\small $A_2$};
  \node[above] at (4.0,5.0) {\small $\eta_1$};
  \node[below] at (4.0,1.45) {\small $\eta_2$};
  \node[anchor=west] at (3.65,3.25) {\small $\eta_3$};
  \node[anchor=west,text=black!60] at (5.95,2.5) {\small inactive};
 \end{tikzpicture}
 \caption{An active loop with three deep arches $A_1,A_3$ (left) and
 $A_2$ (right), drawn dashed, and three connecting paths, drawn solid:
 $\eta_1$ and $\eta_2$ are crossing, $\eta_3$ is returning. Here
 $n=1$, $d=3$, $m=2$, and the two traversals between $r_-$ and $r_+$
 run along $\eta_1$ and $\eta_2$. The gray loop is inactive: its right
 arch does not reach $r_+$.}
 \label{fig:deep-arches}
\end{figure}

Figure~\ref{fig:deep-arches} shows an example. A connecting path may
pass through arches that are not deep.

\begin{lemma}[Deep arches and traversals]\label{lem:deep-traversals}
Let $\omega$ be a full configuration.
\begin{enumerate}[label=\textnormal{(\roman*)},leftmargin=*]
\item If the left cap has $M$ deep arches, there are at least $2M$
traversals between the sections $r_-$ and $r_a$. If the right cap has
$M$ deep arches, there are at least $2M$ traversals between $r_b$ and
$r_+$.
\item The number $m(\omega)$ of crossing connecting paths equals the
number of traversals between $r_-$ and $r_+$. It is even.
\end{enumerate}
\end{lemma}

\begin{proof}
(i) Let $A$ be a deep left arch. Its two ends are on section $r_a$,
and it visits section $r_-$. Follow $A$ from one end to its first
vertex on $r_-$, and let $u$ be the last vertex on $r_a$ before that.
The part of $A$ from $u$ to this first vertex on $r_-$ is a traversal
between $r_-$ and $r_a$. The same construction from the other end of
$A$ gives a second traversal. The two have disjoint interiors, also
when the first and the last visit of $A$ to $r_-$ coincide. Disjoint arches give distinct traversals. The same construction
applies on the right.

(ii) Every visit of a loop to the section $r_-$ lies on a deep left
arch or on a loop contained in the left cap, and a loop contained in
the left cap has no traversal between $r_-$ and $r_+$. The same holds
for $r_+$. Consider a traversal between $r_-$ and $r_+$. It starts on
a deep left arch $A$ and ends on a deep right arch $A'$, or the
reverse, and its interior does not visit $r_-$ or $r_+$. It therefore
contains no intervening deep arch. Hence $A$ and $A'$ are
consecutive deep arches of the loop, and the traversal contains the
connecting path between them, which is crossing. Conversely, let
$\eta$ be a crossing connecting path. It joins a deep left arch $A$ to
a deep right arch $A'$, and it does not visit $r_-$ or $r_+$, since
its parts in the caps are arches that are not deep. The part of the
loop from the last vertex of $A$ on $r_-$ through $\eta$ to the first
vertex of $A'$ on $r_+$ is a traversal between $r_-$ and $r_+$. These
two constructions are inverse to each other. Finally, in the cyclic
order of the deep arches of a loop, the crossing connecting paths are
the places where a left arch is followed by a right arch or a right
arch by a left arch. Their number is even.
\end{proof}

\begin{definition}[Labellings]\label{def:labellings}
Let $\omega$ be a full configuration. A \emph{labelling} of $\omega$
assigns a label $1$ or $2$ to every edge of multiplicity one that lies
between the left and the right cut, the cut edges included. The edges
inside the two caps receive no label. A labelling must satisfy two
rules.
\begin{enumerate}[label=\textnormal{(R\arabic*)},leftmargin=*]
\item At every vertex between the left and the right cut that has two
distinct edges of $\omega$, these two edges have different labels.
\item The cut edges at the two ends of an arch that is not deep have
different labels.
\end{enumerate}
\end{definition}

Label $1$ gives colour $+$ to a cut edge and label $2$ gives colour
$-$. Under (R1) the two labels specify two matchings between the cuts,
with every doubled edge in both. Rule~(R2) assigns the ordinary
endpoint matrix $J$ to each nondeep arch. Deep arches leave the two
endpoint colours free.

\begin{lemma}[Structure of labellings]\label{lem:labellings}
Let $\omega$ be a full configuration.
\begin{enumerate}[label=\textnormal{(\roman*)},leftmargin=*]
\item In a labelling the labels alternate along every connecting path
and along every inactive loop that is not contained in a cap, if the
alternation is continued through the arches that are not deep.
Conversely, choosing one of the two alternating labellings for each of
these paths and loops gives a labelling.
\item In a labelling the two end cut edges of a returning connecting
path have different labels, and those of a crossing connecting path
have the same label.
\end{enumerate}
\end{lemma}

\begin{proof}
(i) Each arch has even length (Lemma~\ref{lem:ordinary-colours}).
Hence (R2) extends alternation through each nondeep arch. With (R1),
this determines all labels on a connecting path or on an inactive loop
meeting the region between the cuts once one edge label is given.
There is no constraint across a deep arch, so these choices are
independent.

(ii) A connecting path that returns to the left cut joins two vertices
$w_L$, and one that returns to the right cut joins two vertices
$b_{\rm mid}$; such a path has an even number of edges. A crossing
connecting path joins a vertex $w_L$ of the left cut to a vertex
$b_{\rm mid}$ of the right cut and has an odd number of edges.
\end{proof}

Let $\Theta$ be the endpoint rule of Definition~\ref{def:colour-tensor}
that assigns the averaging matrix $P$ to every deep arch and the
ordinary matrix $J$ to every other arch. Let
\[
 U_L=\widehat F_L(\Theta)\in\mathcal H_{r_a},\qquad
 U_R=\widehat F_R(\Theta)\in\mathcal H_{r_b}
\]
be the colour tensors of the left cap at $r_a$ and of the right cap at
$r_b$. Put $t=r_b-r_a$. For an integer $k$ with
$0\le k\le\min\{r_a,N+1-r_b\}$ let $U_{L,k}$ and $U_{R,k}$ be the
components of $U_L$ and $U_R$ in the sector $k$, and put
\begin{equation}\label{eq:W-definition}
 \mathcal W_k=
 \frac{\bigl\langle U_{R,k},
   \bigl(\widehat T_t^{(r_a+k)}\otimes\widehat T_t^{(r_a-k)}\bigr)
   U_{L,k}\bigr\rangle}{Z_N^2}.
\end{equation}

\begin{lemma}[Counting identity]\label{lem:sector-count}
For $0\le k\le\min\{r_a,N+1-r_b\}$,
\begin{equation}\label{eq:sector-count}
 \mathcal W_k=\EE\left[2^{-n(\omega)-m(\omega)}
   \binom{m(\omega)}{m(\omega)/2+k}\right],
\end{equation}
where a binomial coefficient $\binom mj$ is zero unless $0\le j\le m$.
\end{lemma}

\begin{proof}
\emph{Step 1. Expansion.}
The gauge of $U_L$ is cancelled by the inverse gauge on the input side
of the transfer, and the gauge of $U_R$ by the gauge on its output
side. We may therefore expand~\eqref{eq:W-definition} using the
ungauged cap tensors of~\eqref{eq:auxiliary-cap-vectors} and the layer
weights $U(T,S)$. The result is that $Z_N^2\mathcal W_k$ is a sum over
the following data.
\begin{enumerate}[label=\textnormal{(\alph*)},leftmargin=*]
\item A superposition $C_L$ of the left cap and a superposition $C_R$
of the right cap.
\item A pair of states $(Q_1^L,Q_2^L)$ at the left cut compatible with
$C_L$, with $|Q_1^L|=r_a+k$ and $|Q_2^L|=r_a-k$, and a pair
$(Q_1^R,Q_2^R)$ at the right cut compatible with $C_R$.
\item For $\nu=1,2$, a matching $M_\nu$ of the layers between the two
cuts with input state $Q_\nu^L$ and output state $Q_\nu^R$.
\end{enumerate}
Its weight is
\begin{equation}\label{eq:term-weight-two-cuts}
 \bk^{\dn}\;2^{\ell(C_L)+\ell(C_R)}
 \prod_{A}\bigl(\Theta_A\bigr)_{\xi(x_A),\xi(y_A)},
\end{equation}
where $\dn$ is the total number of down edges in $C_L$, $C_R$,
$M_1$ and $M_2$, the product is over the arches $A$ of $C_L$ and
$C_R$ with ends $x_A,y_A$, and $\xi$ denotes the colours given by the
states in~(b).

\emph{Step 2. Terms and labellings.}
Take the union, with multiplicities, of $C_L$, $C_R$, $M_1$, $M_2$ and
the cut edges, where the cut edge of a site is taken as many times as
the site occurs in the pair of states at its cut. Every vertex of the
split graph has degree two in this union. Hence the union is a full
configuration $\omega$, with $\bk^{\dn}=\bk^{\dn(\omega)}$, and its
caps are $C_L$ and $C_R$. The nonzero terms with a given $\omega$
correspond bijectively to the labellings of $\omega$ that also satisfy
\begin{enumerate}[label=\textnormal{(R\arabic*)},leftmargin=*,start=3]
\item at the left cut, the number of cut edges with label $1$ minus
the number with label $2$ is $2k$.
\end{enumerate}
Rules (R1) and (R2) enforce the matchings and the nondeep arch
factors, respectively. Rule (R3) fixes the sector in (b), since
doubly occupied sites contribute equally to both states. Each deep
arch contributes $\tfrac12$. The factor $2^{\ell(C_L)+\ell(C_R)}$
already counts the two choices on each loop contained in a cap.

\emph{Step 3. Counting.}
In the imbalance (R3) the two ends of each nondeep left arch cancel
by (R2). The remaining sites are endpoints of connecting paths.
Returning paths at the left cut contribute zero by
Lemma~\ref{lem:labellings}(ii). Each crossing path contributes
$+1$ or $-1$, so exactly $m(\omega)/2+k$ must have label $1$ at their
endpoints. The remaining $d(\omega)-m(\omega)$ paths each admit both
labellings. By Lemma~\ref{lem:labellings}(i) this gives
$\binom{m(\omega)}{m(\omega)/2+k}\,2^{d(\omega)-m(\omega)}$
labellings of the connecting paths, and every inactive loop
contributes the factor $2$. With the factors $\tfrac12$ of the
$d(\omega)$ deep arches, the sum of the
weights~\eqref{eq:term-weight-two-cuts} over all terms with a given
$\omega$ is
\[
 \bk^{\dn(\omega)}\;2^{\ell(\omega)-n(\omega)}\;
 \binom{m(\omega)}{m(\omega)/2+k}\,2^{d(\omega)-m(\omega)}\;
 2^{-d(\omega)}
 =\bk^{\dn(\omega)}2^{\ell(\omega)}\cdot
   2^{-n(\omega)-m(\omega)}\binom{m(\omega)}{m(\omega)/2+k}.
\]
Dividing by $Z_N^2$ and
using~\eqref{eq:configuration-probability}
gives~\eqref{eq:sector-count}.
\end{proof}

\begin{lemma}[Recursive inequality]
\label{lem:traversal-recursion}
For integers $1\le r_-<r_a<r_b<r_+\le N$ and integers $k,M\ge1$,
\begin{equation}\label{eq:recursion}
 \begin{split}
 P_N^{\mathrm{tr}}(k;r_-,r_+)
 \le{}&P_N^{\mathrm{tr}}(M;r_-,r_a)
       +P_N^{\mathrm{tr}}(M;r_b,r_+)\\
 &+2^{4M}\exp\left(-\frac{(r_b-r_a)k^2}{N+1}\right).
 \end{split}
\end{equation}
\end{lemma}

\begin{proof}
By Lemma~\ref{lem:deep-traversals}(ii) the event on the left
of~\eqref{eq:recursion} is $\{m\ge2k\}$. Let $\mathcal F$ be the
event that each cap has fewer than $M$ deep arches. Part~(i) gives
\[
 \PP(\mathcal F^c)\le P_N^{\mathrm{tr}}(M;r_-,r_a)+P_N^{\mathrm{tr}}(M;r_b,r_+).
\]
It remains to bound $\PP(\mathcal F\cap\{m\ge2k\})$.

Different crossing paths end at different singly occupied sites of the
left cut and of the right cut. The singly occupied sites of a cut at
section $r$ form the symmetric difference of two $r$-element subsets
of $E_N$, so there are at most $2\min\{r,N+1-r\}$ of them. Hence
$m(\omega)\le2\min\{r_a,N+1-r_b\}$, and we may assume
$k\le\min\{r_a,N+1-r_b\}$.

On $\mathcal F\cap\{m\ge2k\}$, the binomial coefficient
in~\eqref{eq:sector-count} is at least $1$. Moreover
$n(\omega)\le d(\omega)<2M$, because every active loop contains a deep
arch, and $m(\omega)\le d(\omega)<2M$. The integrand
in~\eqref{eq:sector-count} is at least $2^{-4M}$ on this event, and
Lemma~\ref{lem:sector-count} gives
\[
 \PP(\mathcal F\cap\{m\ge2k\})\le2^{4M}\,\mathcal W_k .
\]
The matrices $P$ and $J$ have operator norm one, so
Lemma~\ref{lem:aztec-marked-cap} gives
$\|U_{L,k}\|\le\|U_L\|\le\|\widehat L_{r_a}\|^2$ and
$\|U_{R,k}\|\le\|\widehat R_{r_b}\|^2$. By the Cauchy--Schwarz
inequality, the normalization~\eqref{eq:normalization} and the sector
estimate~\eqref{eq:sector-cost} with $r=r_a$ and $t=r_b-r_a$,
\[
 \mathcal W_k
 \le\frac{\|\widehat L_{r_a}\|^2\|\widehat R_{r_b}\|^2}{Z_N^2}\,
    \|\widehat T_t^{(r_a+k)}\|\,\|\widehat T_t^{(r_a-k)}\|
 =\frac{\|\widehat T_t^{(r_a+k)}\|\,
          \|\widehat T_t^{(r_a-k)}\|}{G_{r_a,r_b}^{\,2}}
 \le \exp\left(-\frac{tk^2}{N+1}\right).
\]
Adding the three bounds proves~\eqref{eq:recursion}.
\end{proof}

\subsection{The Gaussian tail}
\label{subsec:gaussian-tail}

\begin{proposition}[Uniform Gaussian traversal tail]
\label{prop:gaussian-traversal-tail}
For every $0<\delta<1$, there are constants
$C_\delta,c_\delta>0$ such that, for every integer $k\ge0$,
\begin{equation}\label{eq:gaussian-traversal-integer}
 \sup_{\substack{
   N,r,r'\in\mathbb Z\\
   N\ge1,\;1\le r<r'\le N\\
   r'-r\ge\delta(N+1)
 }}
 P_N^{\mathrm{tr}}(k;r,r')
 \le C_\delta e^{-c_\delta k^2}.
\end{equation}
In particular, for every fixed $0<u<v<1$, there are constants
$C,c_0>0$ such that, for every integer $k\ge0$,
\begin{equation}\label{eq:gaussian-tail}
 \sup_{N\ge1}\PP\{\NN_N(u,v)\ge2k\}
 \le Ce^{-c_0k^2}.
\end{equation}
\end{proposition}

\begin{proof}
\emph{Capacity.}
Section $r$ has $N+1$ white vertices, and each lies on at most one
loop, with two loop edges. Different traversals between $r$ and $r'$
end at $r$ along different loop edges, so there are at most $2(N+1)$
of them. Hence
\begin{equation}\label{eq:capacity}
 P_N^{\mathrm{tr}}(k;r,r')=0\qquad\text{for }k>N+1 .
\end{equation}

\emph{One splitting step.}
Let $r<r'$ be sections of width $\Delta=r'-r\ge4$, and let $k'\ge1$.
Apply Lemma~\ref{lem:traversal-recursion} with
\[
 r_-=r,\qquad r_a=r+\lfloor\Delta/4\rfloor,\qquad
 r_b=r'-\lfloor\Delta/4\rfloor,\qquad r_+=r',\qquad M=8k' .
\]
The two outer slabs have width $\lfloor\Delta/4\rfloor$ and the middle
has width $r_b-r_a\ge\Delta/2$. Thus
\begin{equation}\label{eq:one-split}
 P_N^{\mathrm{tr}}(k';r,r')
 \le P_N^{\mathrm{tr}}(8k';r,r_a)+P_N^{\mathrm{tr}}(8k';r_b,r')
   +2^{32k'}\exp\left(-\frac{\Delta\,k'^2}{2(N+1)}\right).
\end{equation}

\emph{Iteration.}
Let $r'-r=\Delta\ge\delta(N+1)$ and $k\ge256/\delta$. Apply
\eqref{eq:one-split} to $P_N^{\mathrm{tr}}(k;r,r')$, then to the two
probabilities on its right-hand side, and so on. After $j$ steps there
are at most $2^j$ probabilities, each with threshold and width
\[
 k_j=8^jk,\qquad \Delta_j=\lfloor\Delta/4^j\rfloor .
\]
A probability with $k_j>N+1$ is zero by~\eqref{eq:capacity} and is not
split further. If $k_j\le N+1$, then
\begin{equation}\label{eq:recursive-width}
 \frac\Delta{4^j}\ge\frac{\delta(N+1)}{4^j}
              \ge\frac{\delta k_j}{4^j}=\delta\,2^j k\ge8 ,
\end{equation}
so $\Delta_j\ge4$, the splitting step applies, and
$\Delta_j\ge\Delta/(2\cdot4^j)\ge\delta(N+1)/(2\cdot4^j)$. The last
term of~\eqref{eq:one-split} at such a step is at most
\[
 \exp\left(32(\log2)k_j-\frac{\delta k_j^2}{4\cdot4^j}\right)
 \le\exp\left(-\frac{\delta k_j^2}{8\cdot4^j}\right)
 =\exp\left(-\frac\delta8\,16^j k^2\right).
\]
The inequality holds because
$32(\log2)k_j\le\delta k_j^2/(8\cdot4^j)$ is equivalent to
$256\log2\le\delta\,2^jk$. The thresholds $k_j$ increase, so the
iteration stops after finitely many steps, and
\begin{equation}\label{eq:recursion-sum}
 P_N^{\mathrm{tr}}(k;r,r')
 \le \sum_{j\ge0}2^j\exp\left(-\frac{\delta}{8}16^j k^2\right)
 \le\frac{\exp(-\delta k^2/8)}
           {1-2\exp(-15\delta k^2/8)}
 \le2\exp(-\delta k^2/8).
\end{equation}
The second inequality uses $16^j\ge1+15j$, and the third uses
$k\ge256/\delta$. For $k<256/\delta$ the
bound~\eqref{eq:gaussian-traversal-integer} holds after enlarging
$C_\delta$.

\emph{Fixed levels.}
For fixed $u<v$ there are $\delta>0$ and $N_0$ such that
$r_v-r_u\ge\delta(N+1)$ for $N\ge N_0$, and
\eqref{eq:gaussian-traversal-integer} applies. For $N<N_0$ use
\eqref{eq:capacity} and enlarge the constant. This
proves~\eqref{eq:gaussian-tail}.
\end{proof}

In the coordinates $(s,y)$, a \emph{slab} $\mathcal R$ is the region
between two lines $s=u$, $s=v$ or between two lines $y=u$, $y=v$, with
$0<u<v<1$. Write $\NN_N(\mathcal R)$ for the number of traversals
between its boundary sections. For slabs between $y=u$ and $y=v$ the
sections are those of the reflected cut construction of
Section~\ref{subsec:reflection}.

\begin{proposition}[Uniform traversal moments]
\label{prop:traversal-moments}
For every fixed $0<u<v<1$, there is $\lambda>0$ such that
\[
 \sup_N\EE e^{\lambda\NN_N(u,v)^2}<\infty.
\]
For every finite family of fixed slabs
$\mathcal R_1,\ldots,\mathcal R_J$ and every $A<\infty$,
\begin{equation}\label{eq:traversal-exponential-moment}
 \sup_N\EE\exp\left(
 A\sum_{j=1}^J\NN_N(\mathcal R_j)
 \right)<\infty.
\end{equation}
\end{proposition}

\begin{proof}
Summing the Gaussian tail of
Proposition~\ref{prop:gaussian-traversal-tail} gives the first
assertion for sufficiently small $\lambda>0$. The inequality
$Ax\le\lambda x^2+A^2/(4\lambda)$ then gives all exponential moments
of $\NN_N(u,v)$. For a slab between $y=u$ and $y=v$, apply the
reflection $\varrho$ of Section~\ref{subsec:reflection} to the two
perfect matchings and to the cut construction. The reflected pair has the same
law, and $\varrho$ exchanges the lines $s=\mathrm{const}$ and
$y=\mathrm{const}$, so the same bounds hold. H\"older's inequality
gives~\eqref{eq:traversal-exponential-moment}.
\end{proof}

\section{An operator bound for the two-particle insertion}
\label{sec:insertion}

A two-particle insertion transfers a pair of occupied sites from one
replica to the other. We place the insertion between two transfers
of macroscopic length. Proposition~\ref{prop:full-seam-insertion}
bounds the sum of the resulting operator norms over the cells of a
partition of the section by $C h\log(2/h)$ after normalization,
where $h$ is the mesh. This summability across the whole section
will control the probability of opposite crossings in
Section~\ref{sec:global-opposite}.

Throughout Sections~\ref{sec:insertion}
and~\ref{sec:global-opposite} we fix $c\in(0,1)$, positive numbers
$\rho_-,\rho_+$ with $0<c-\rho_-<c+\rho_+<1$, and two levels $a,b$ with
\[
 c-\rho_-<a<c<b<c+\rho_+.
\]
Here and in Section~\ref{sec:global-opposite}, $a$ and $b$ denote
cut levels. The model bias enters through $\bk$.
Let $r_a,r_c,r_b$ be the indices of the white sections nearest to
these levels, as in Section~\ref{subsec:seams}. The left cap is taken
at section $r_a$, the right cap at section $r_b$, and $r_c$ is the
\emph{middle section}. Put $t_-=r_c-r_a$ and $t_+=r_b-r_c$ for the
transfer lengths. Both are at least $\eta N$, for some fixed $\eta>0$
and all sufficiently large $N$.

For $S\subseteq E_N$ and $x\in E_N$, write
$n_x(S)=|\{s\in S:s<x\}|$. Define
\[
 c_x\varepsilon_S=
 \begin{cases}
  (-1)^{n_x(S)}\varepsilon_{S\setminus\{x\}}, & x\in S,\\
  0, & x\notin S,
 \end{cases}
 \qquad
 c_x^\dagger\varepsilon_S=
 \begin{cases}
  (-1)^{n_x(S)}\varepsilon_{S\cup\{x\}}, & x\notin S,\\
  0, & x\in S,
 \end{cases}
\]
and extend both operators linearly.

Let $I$ be an interval of the transverse coordinate $y$, and let
$S(I)=\{x\in E_N:x/N\in I\}$ be the set of sites of the middle section
in $I$. The \emph{two-particle insertion} is
\begin{equation}\label{eq:insertion}
 \mathcal O_I
 =\sum_{\substack{i<j\\i,j\in S(I)}}c_i c_j\otimes c_j^\dagger c_i^\dagger .
\end{equation}
The first factor acts on the state of replica~1 and the second on the
state of replica~2. A summand vanishes on
$\varepsilon_{S_1}\otimes\varepsilon_{S_2}$ unless
$i,j\in S_1\setminus S_2$. Otherwise it removes $i$ and $j$ from $S_1$
and adds them to $S_2$, up to a sign, and leaves all other sites
unchanged. Removing $i,j$ from replica~1 divides the gauge factor by
$d_id_j$, where $d_x=D_{xx}$, and adding them to replica~2 multiplies
it by the same factor. Hence $\mathcal O_I$ commutes with the gauge
transformation.

Work on the direct sum $\mathcal H_r$ of two-replica sectors with
$2r$ particles in total, defined in Section~\ref{subsec:marked-caps}.
For a transfer of $t$ layers from section $r$, write
\[
 \widehat{\mathcal T}_{r,t}
 =\bigoplus_{p_1+p_2=2r}
     \bigl(\widehat T_t^{(p_1)}\otimes\widehat T_t^{(p_2)}\bigr)
 :\mathcal H_r\longrightarrow\mathcal H_{r+t}.
\]
Only degrees $0\le p_1,p_2\le N+1$ occur in the source, and we
interpret a block as zero if $p_1+t$ or $p_2+t$ exceeds $N+1$.
Set $\widehat{\mathcal T}_-=\widehat{\mathcal T}_{r_a,t_-}$ and
$\widehat{\mathcal T}_+=\widehat{\mathcal T}_{r_c,t_+}$.
The insertion acts within $\mathcal H_{r_c}$ by transferring two
particles from the first replica to the second. For each admissible
integer $s$, the particle numbers in a nonzero block follow
\begin{equation}\label{eq:sectors}
 \begin{aligned}
 (r_a+s,r_a-s)&\xrightarrow{\widehat{\mathcal T}_-}(r_c+s,r_c-s)\\
 &\xrightarrow{\mathcal O_I}(r_c+s-2,r_c-s+2)
 \xrightarrow{\widehat{\mathcal T}_+}(r_b+s-2,r_b-s+2).
 \end{aligned}
\end{equation}
Set $G_-=G_{r_a,r_c}$ and $G_+=G_{r_c,r_b}$,
where $G_{r,r'}$ is defined in~\eqref{eq:cap-transfer-factor}.

Set $Q_{xj}=\psi_j(x)$ and $P_I=Q^{\mathsf T}\one_{S(I)}Q$.

\begin{lemma}\label{lem:insertion-weights}
Fix $\varepsilon>0$. Consider quadruples
\[
 A_-,B_-\subseteq\{0,\ldots,N-t_-\},\qquad
 A',B'\subseteq\{0,\ldots,N-t_+\},
\]
and put
\[
 A=\{0,\ldots,t_--1\}\cup(A_-+t_-),\qquad
 B=\{0,\ldots,t_--1\}\cup(B_-+t_-).
\]
Call the quadruple compatible if
\[
 \begin{gathered}
 |A_-|+|B_-|=2r_a,\\
 A'\subseteq A,\qquad B\subseteq B',\qquad
 |A\setminus A'|=|B'\setminus B|=2.
 \end{gathered}
\]
For a compatible quadruple, set $R=A\setminus A'$ and
$T=B'\setminus B$. The corresponding absolute matrix entry of
$\widehat{\mathcal T}_+\mathcal O_I\widehat{\mathcal T}_-$,
divided by $G_+^2G_-^2$, has the form
$v\,|\det(P_I)_{R,T}|$, with $v\ge0$ independent of $I$.
For all sufficiently large $N$,
\begin{equation}\label{eq:insertion-weight-sums}
 \sum_{\text{compatible quadruples}}v\le C,\qquad
 \sum_{\substack{\text{compatible quadruples}\\
                   \exists j\in R\cup T:\ |j-r_c|>\varepsilon N}}
 v\le Ce^{-c_1N}.
\end{equation}
The sums include all admissible imbalances $s$ in~\eqref{eq:sectors}
and quadruples with zero insertion minor.
The constants $C,c_1>0$ depend only on $\varepsilon$ and the fixed sections.
\end{lemma}

\begin{proof}
Each transfer takes a Krawtchouk basis vector to a multiple of
another basis vector or to zero. We first compute the insertion
between two such vectors.

\emph{The insertion minor.}
For $R,T\subseteq E_N$ with $|R|=|T|=2$,
the Cauchy--Binet formula gives
\begin{equation}\label{eq:minor}
 \sum_{\substack{i<j\\i,j\in S(I)}}
 \det Q_{\{i,j\},R}\det Q_{\{i,j\},T}
 =\det(P_I)_{R,T}.
\end{equation}
For a compatible quadruple, the two exterior-product coefficients
of the insertion in the basis of Section~\ref{subsec:exact-transfer}
are, by antisymmetry,
\[
 \langle e_{A'},c_ic_je_A\rangle
 =\sigma_1\det Q_{\{i,j\},R},
 \qquad
 \langle e_{B'},c_j^\dagger c_i^\dagger e_B\rangle
 =\sigma_2\det Q_{\{i,j\},T}.
\]
The signs $\sigma_1,\sigma_2\in\{1,-1\}$ come from the ordering
of the occupied modes and are independent of $i,j$.
Multiplying these coefficients and summing gives
\begin{equation}\label{eq:insertion-matrix-entry}
\begin{aligned}
 \bigl\langle e_{A'}\otimes e_{B'},\,
 \mathcal O_I(e_A\otimes e_B)\bigr\rangle
 &=
 \sigma\sum_{\substack{i<j\\i,j\in S(I)}}
 \det Q_{\{i,j\},R}\det Q_{\{i,j\},T}\\
 &=\sigma\det((P_I)_{R,T}),
\end{aligned}
\end{equation}
where $\sigma=\sigma_1\sigma_2$ and the last equality
follows from~\eqref{eq:minor}.

\emph{The transfer weights.}
For a transfer of $t\ge\eta N$ layers starting at section $r$, set
\[
 F_r=\{0,\ldots,r-1\},\qquad
 \mu=\sqrt{q_t(r-1)q_t(r)}.
\]
By~\eqref{eq:t-transfer}, $\widehat T_t^{(|A|)}e_A\ne0$
if and only if $A\subseteq\{0,\ldots,N-t\}$.
For each such $A$, define
\[
 w(A)=\prod_{j\in A\setminus F_r}\frac{q_t(j)}{\mu}
       \prod_{j\in F_r\setminus A}\frac{\mu}{q_t(j)}.
\]
Here $r=r(N)$ and $t=t(N)$ are the section index and transfer
length arising from the fixed macroscopic sections above.
Since $(r+t)/N$ stays strictly below $1$, we have
$r+t\le N$ for all sufficiently large $N$, and hence
$q_t(r)>0$ by~\eqref{eq:q-definition}.
By~\eqref{eq:q-gap} and the definition of $\mu$,
\begin{equation}\label{eq:weight-factor-bound}
 \begin{aligned}
 \frac{q_t(\ell)}{\mu}
 &\le e^{-\frac{t}{N+1}(\ell-r+1/2)},
 &&r\le \ell\le N-t,\\
 \frac{\mu}{q_t(\ell)}
 &\le e^{-\frac{t}{N+1}(r-\ell-1/2)},
 &&0\le \ell<r.
 \end{aligned}
\end{equation}
Expanding the products below and using
$t/(N+1)\ge\eta/2$ gives
\begin{equation}\label{eq:state-sum}
 \begin{aligned}
 \sum_{A\subseteq\{0,\ldots,N-t\}}w(A)
 &=
 \prod_{j=0}^{r-1}
 \left(1+\frac{\mu}{q_t(j)}\right)
 \prod_{j=r}^{N-t}
 \left(1+\frac{q_t(j)}{\mu}\right)\\
 &\le
 \prod_{m=0}^{\infty}
 \left(1+e^{-\frac{\eta}{2}(m+1/2)}\right)^2
 =:C_\eta<\infty.
 \end{aligned}
\end{equation}
By~\eqref{eq:t-transfer}, when $|A|+|B|=2r$,
\begin{equation}\label{eq:two-replica-transfer-weight}
 \frac{B_{N,t}^{(|A|)}B_{N,t}^{(|B|)}
       \prod_{j\in A}q_t(j)\prod_{j\in B}q_t(j)}
      {G_{r,r+t}^2}
 =\mu^{|A|+|B|-2r}w(A)w(B)
 =w(A)w(B).
\end{equation}
Here we used~\eqref{eq:B-pair} in the form
$B_{N,t}^{(|A|)}B_{N,t}^{(|B|)}=\bigl(B_{N,t}^{(r)}\bigr)^2$ and
$G_{r,r+t}=B_{N,t}^{(r)}\prod_{j=0}^{r-1}q_t(j)$.
Since $w(A)\ge0$, \eqref{eq:state-sum} gives
\[
 \sum_{\substack{
 A,B\subseteq\{0,\ldots,N-t\}\\
 |A|+|B|=2r}}
 w(A)w(B)
 \le
 \left(\sum_{A\subseteq\{0,\ldots,N-t\}}w(A)\right)^2
 \le C_\eta^2.
\]

\emph{The normalized entry.}
Let $w_-$ and $w_+$ be obtained from the definition of
$w(A)$ by setting $(r,t)=(r_a,t_-)$ and $(r_c,t_+)$,
respectively.
Let $e_{A_-}\otimes e_{B_-}$ be an input basis vector of
$\widehat{\mathcal T}_-$, the transfer from section
$r_a$ to section $r_c$.
By~\eqref{eq:t-transfer}, its nonzero output satisfies
\[
 A=\{0,\ldots,t_--1\}\cup(A_-+t_-),\qquad
 B=\{0,\ldots,t_--1\}\cup(B_-+t_-).
\]
Set
\[
 v=w_-(A_-)w_-(B_-)w_+(A')w_+(B').
\]
By~\eqref{eq:t-transfer}, the nonzero images of distinct
input basis vectors under either $\widehat{\mathcal T}_-$
or $\widehat{\mathcal T}_+$ are scalar multiples of
distinct output basis vectors.
Let $A_+,B_+$ be the output index sets of
$\widehat{\mathcal T}_+$, so that
\[
 A_+=\{0,\ldots,t_+-1\}\cup(A'+t_+),\qquad
 B_+=\{0,\ldots,t_+-1\}\cup(B'+t_+).
\]
By~\eqref{eq:t-transfer}, only the
$e_{A'}\otimes e_{B'}$ component of
$\mathcal O_I(e_A\otimes e_B)$ contributes to the matrix entry
\[
 \left\langle
 e_{A_+}\otimes e_{B_+},\,
 \widehat{\mathcal T}_+\mathcal O_I\widehat{\mathcal T}_-
 (e_{A_-}\otimes e_{B_-})
 \right\rangle.
\]
Combining~\eqref{eq:t-transfer},
\eqref{eq:two-replica-transfer-weight}, and
\eqref{eq:insertion-matrix-entry} gives
\begin{equation}\label{eq:normalized-insertion-entry}
 \frac{
 \left|
 \left\langle
 e_{A_+}\otimes e_{B_+},\,
 \widehat{\mathcal T}_+\mathcal O_I\widehat{\mathcal T}_-
 (e_{A_-}\otimes e_{B_-})
 \right\rangle
 \right|}
 {G_+^2G_-^2}
 =v\,|\det(P_I)_{R,T}|.
\end{equation}

\emph{Modes far from $r_c$.}
Suppose some $j\in R\cup T$ satisfies
$|j-r_c|>\varepsilon N$.
If $j>r_c+\varepsilon N$, at least one of $A,B'$
contains $j\notin F_{r_c}$.
If $j<r_c-\varepsilon N$, at least one of $A',B$
omits $j\in F_{r_c}$.
By~\eqref{eq:t-transfer}, in each replica the nonzero
output contains the fixed factors $\psi_0,\ldots,\psi_{t_--1}$,
together with the input factors whose indices are increased
by $t_-$.
Thus an output index $j\ge t_-$ corresponds to the input
index $\ell=j-t_-$. Since $t_-=r_c-r_a$,
\[
 \ell-r_a=j-t_--r_a=j-r_c.
\]
For the right transfer $\widehat{\mathcal T}_+$,
apply~\eqref{eq:weight-factor-bound} with
$r=r_c$, $t=t_+$, and $\ell=j$
to bound $w_+(A')$ or $w_+(B')$.
Consequently, one of the four weights in $v$ contains
an included index above $F_r$ or a missing index below
$F_r$, at distance greater than $\varepsilon N$ from $r$.
The preceding factor bounds give, for sufficiently large $N$,
\[
 \begin{aligned}
 &\min\{w_-(A_-),w_-(B_-),w_+(A'),w_+(B')\}\\
 &\qquad\le
 e^{-\frac{\eta}{2}(\varepsilon N-1/2)}
 \le e^{-\eta\varepsilon N/4}.
 \end{aligned}
\]
All four weights are at most $1$, so
\[
 v\le e^{-\eta\varepsilon N/4},
 \qquad
 v\le e^{-\eta\varepsilon N/8}\sqrt v.
\]

\emph{Summing over all sectors.}
The product calculation in~\eqref{eq:state-sum}, with
the exponents halved, also gives
\[
 \sum_{A\subseteq\{0,\ldots,N-t\}}\sqrt{w(A)}\le C.
\]
Each transfer is injective on the basis vectors it does not annihilate.
Thus the input and output indices determine $(A_-,B_-,A',B')$.
Dropping the compatibility conditions gives
\[
 \sum_{\mathrm{compatible}}\sqrt v
 \le
 \left(\sum_{A_-\subseteq\{0,\ldots,N-t_-\}}\sqrt{w_-(A_-)}\right)^2
 \left(\sum_{A'\subseteq\{0,\ldots,N-t_+\}}\sqrt{w_+(A')}\right)^2
 \le C.
\]
Since $v\le\sqrt v$, the total weight is bounded by $C$. If a mode
lies more than $\varepsilon N$ from $r_c$, we have the stronger bound
$v\le e^{-\eta\varepsilon N/8}\sqrt v$. Summing only over these
quadruples gives the exponential estimate.
This proves~\eqref{eq:insertion-weight-sums}.
\end{proof}

\subsection{An atomless spatial limit for Krawtchouk modes}

For $0\le j\le N$, define the probability measure
\[
 \nu_{N,j}=\sum_{x=0}^N\psi_j(x)^2\,\delta_{x/N}.
\]
The weak limit has inverse-square-root singularities at the two
turning points. The squared cell masses remain summable there.

\begin{lemma}
\label{lem:atomless-mode-masses}
If $j_N/N\to\alpha\in(0,1)$, then
\begin{equation}\label{eq:mode-arcsine-limit}
 \nu_{N,j_N}\Longrightarrow\nu_\alpha,
 \qquad
 \frac{d\nu_\alpha}{dy}
 =\frac{\one_{\{|y-m_\alpha|<a_\alpha\}}}
        {\pi\sqrt{a_\alpha^2-(y-m_\alpha)^2}},
\end{equation}
where
\begin{equation}\label{eq:mode-center-width}
 m_\alpha=\frac{(1-\alpha)+\bk\alpha}{1+\bk},
 \qquad
 a_\alpha=\frac{2\sqrt{\bk\alpha(1-\alpha)}}{1+\bk}.
\end{equation}
The interval $[m_\alpha-a_\alpha,m_\alpha+a_\alpha]$ is contained in
$[0,1]$, because
$(1+\bk)(m_\alpha\pm a_\alpha)
=(\sqrt{1-\alpha}\pm\sqrt{\bk\alpha})^2$ and
$(\sqrt{1-\alpha}+\sqrt{\bk\alpha})^2\le1+\bk$
by the Cauchy--Schwarz inequality.
For every compact interval $A\Subset(0,1)$ there is $C_A<\infty$
such that, for $0<h\le1$ and every family $\mathcal P_{N,h}$ of
pairwise disjoint intervals in the transverse coordinate $y$,
each of length at most $h+O(N^{-1})$,
\begin{equation}\label{eq:mode-partition-square}
 \limsup_{N\to\infty}\sup_{j/N\in A}
 \sum_{I\in\mathcal P_{N,h}}\nu_{N,j}(I)^2
 \le C_A h\log(2/h).
\end{equation}
\end{lemma}

\begin{proof}
The generating function $F_x(z)=(\bk+z)^{N-x}(1-z)^x$
in~\eqref{eq:krawtchouk} satisfies
$(\bk+z)(1-z)F_x'(z)=\bigl[N-(1+\bk)x-Nz\bigr]F_x(z)$.
Comparing coefficients of $z^j$ gives
\[
 \bigl(N-(1+\bk)x\bigr)K_j(x)
 =\bk(j+1)K_{j+1}(x)+(1-\bk)jK_j(x)+(N-j+1)K_{j-1}(x).
\]
With the normalization~\eqref{eq:psi}, this becomes
\begin{equation}\label{eq:krawtchouk-jacobi}
 \frac{x}{N}\psi_j(x)
 =\frac{(N-j)+\bk j}{(1+\bk)N}\psi_j(x)
 +\frac{\sqrt{\bk(j+1)(N-j)}}{(1+\bk)N}\psi_{j+1}(x)
 +\frac{\sqrt{\bk j(N-j+1)}}{(1+\bk)N}\psi_{j-1}(x),
\end{equation}
where $\psi_{-1}=\psi_{N+1}=0$ by convention.
Let $J_N$ be the symmetric $(N+1)\times(N+1)$ matrix
defined by
\begin{equation}\label{eq:jacobi-matrix}
 \begin{aligned}
 (J_N)_{jj}
 &=\frac{(N-j)+\bk j}{(1+\bk)N},
 &&0\le j\le N,\\
 (J_N)_{j,j+1}=(J_N)_{j+1,j}
 &=\frac{\sqrt{\bk(j+1)(N-j)}}{(1+\bk)N},
 &&0\le j<N,
 \end{aligned}
\end{equation}
with all other entries zero.
By~\eqref{eq:krawtchouk-jacobi}
and the orthogonality of $Q$,
\[
 J_N
 =Q^{\mathsf T}\operatorname{diag}
   \left(0,\frac1N,\ldots,1\right)Q.
\]
Consequently, for every integer $m\ge0$,
\[
 \int y^m\,d\nu_{N,j}(y)
 =\sum_{x=0}^N
   \left(\frac{x}{N}\right)^m\psi_j(x)^2
 =(J_N^m)_{jj}.
\]

Fix an integer $m\ge0$. Expanding $(J_N^m)_{j_Nj_N}$
gives a finite sum over length-$m$ paths with steps in
$\{-1,0,1\}$, starting and ending at $j_N$.
These paths stay within distance $m$ of $j_N$. Along all of them
the diagonal coefficients tend to $m_\alpha$,
and both off-diagonal coefficients tend to $a_\alpha/2$.
We may take the limit term by term.
Averaging
$(m_\alpha+\tfrac{a_\alpha}{2}e^{\ii\theta}
+\tfrac{a_\alpha}{2}e^{-\ii\theta})^m$
over $\theta\in[0,2\pi]$ retains exactly the terms with
equally many $+1$ and $-1$ steps. Hence
\[
 \begin{aligned}
 \lim_{N\to\infty}\int y^m\,d\nu_{N,j_N}(y)
 &=\lim_{N\to\infty}(J_N^m)_{j_Nj_N}\\
 &=\frac1{2\pi}\int_0^{2\pi}
   \bigl(m_\alpha+a_\alpha\cos\theta\bigr)^m\,d\theta.
 \end{aligned}
\]
The limit moments are those of $m_\alpha+a_\alpha\cos\vartheta$, with
$\vartheta$ uniform on $[0,2\pi]$. This law has the density
in~\eqref{eq:mode-arcsine-limit}. Since all the measures are supported
on $[0,1]$, polynomial approximation gives weak convergence.

For $\alpha\in A$, $a_\alpha$ is bounded below by a positive
constant. Integrating the arcsine density at the two endpoints gives
\begin{equation}\label{eq:arcsine-interval-mass}
 \sup_{\alpha\in A}\sup_{|I|\le h}\nu_\alpha(I)
 \le C_A\sqrt h.
\end{equation}
For a fixed $h>0$, disjointness of the cells gives
\[
 \sum_{I\in\mathcal P_{N,h}}\nu_{N,j}(I)^2
 = (\nu_{N,j}\otimes\nu_{N,j})
   \left(\bigcup_{I\in\mathcal P_{N,h}} I\times I\right)
 \le (\nu_{N,j}\otimes\nu_{N,j})
       \{|y_1-y_2|\le h+O(N^{-1})\}.
\]
Along a subsequence realizing the limsup in~\eqref{eq:mode-partition-square},
pick maximizing indices $j_N$. Passing to a further subsequence gives
$j_N/N\to\alpha\in A$. Product weak convergence and the Portmanteau
bound on the closed sets $\{|y_1-y_2|\le h+\delta\}$, followed by
$\delta\downarrow0$, bound that limsup by
$(\nu_\alpha\otimes\nu_\alpha)\{|y_1-y_2|\le h\}$.

Write $d(y)$ for the distance from $y$ to the nearer endpoint of the
support of $\nu_\alpha$. By~\eqref{eq:arcsine-interval-mass}, both
$\nu_\alpha\{d\le2h\}$ and $\sup_y\nu_\alpha([y-h,y+h])$ are
$O_A(\sqrt h)$. The endpoint layers thus contribute $O_A(h)$ to
the product probability. On $\{d(y)>2h\}$ the arcsine density at $y$
is bounded by $C_A/\sqrt{d(y)}$. For $|y-z|\le h$ we also have
$d(z)\ge d(y)/2$, so $\nu_\alpha([y-h,y+h])\le C_Ah/\sqrt{d(y)}$.
For $2h<a_\alpha$ this gives
\[
 \begin{aligned}
 (\nu_\alpha\otimes\nu_\alpha)\{|y_1-y_2|\le h\}
 &=\int \nu_\alpha([y-h,y+h])\,d\nu_\alpha(y)\\
 &\le C_Ah+2C_Ah\int_{2h}^{a_\alpha}\frac{du}{u}
 \le C_A' h\log(2/h).
 \end{aligned}
\]
Since $\inf_{\alpha\in A}a_\alpha>0$, increasing $C_A'$ covers the
remaining $h\le1$. This proves~\eqref{eq:mode-partition-square}.
\end{proof}

\begin{lemma}
\label{lem:partition-minors}
Let $\mathcal P$ be a family of pairwise disjoint intervals
in the transverse coordinate $y$, and let
$P_I=Q^{\mathsf T}\one_{S(I)}Q$ as in~\eqref{eq:minor}.
For $R=\{a,b\}$ and $T=\{c,d\}$ with
$0\le a<b\le N$ and $0\le c<d\le N$,
\begin{equation}\label{eq:partition-minors-universal}
 \sum_{I\in\mathcal P}|\det(P_I)_{R,T}|\le1.
\end{equation}
Moreover,
\begin{equation}\label{eq:partition-minors-square}
 \sum_{I\in\mathcal P}|\det(P_I)_{R,T}|
 \le\prod_{j\in(a,b,c,d)}
        \left(\sum_{I\in\mathcal P}\nu_{N,j}(I)^2\right)^{1/4}.
\end{equation}
The product retains repeated indices if $R$ and $T$ overlap.
\end{lemma}

\begin{proof}
The exterior-product Cauchy--Schwarz inequality and Hadamard's
inequality give
\[
 |\det(P_I)_{R,T}|
 \le\sqrt{\det(P_I)_{R,R}\det(P_I)_{T,T}}
 \le\sqrt{\nu_{N,a}(I)\nu_{N,b}(I)
                  \nu_{N,c}(I)\nu_{N,d}(I)}.
\]
Here $(P_I)_{jj}=\nu_{N,j}(I)$.
Applying H\"older's inequality with four exponents equal to $4$
to this bound proves~\eqref{eq:partition-minors-square}.
Disjointness gives
\[
 \sum_{I\in\mathcal P}\nu_{N,j}(I)^2
 \le\sum_{I\in\mathcal P}\nu_{N,j}(I)\le1.
\]
Substituting this into~\eqref{eq:partition-minors-square} proves
\eqref{eq:partition-minors-universal}.
\end{proof}

\begin{proposition}
\label{prop:full-seam-insertion}
Fix $0<a<c<b<1$, and choose white-section indices
$r_a,r_c,r_b$ such that
\[
 r_a/N\to a,\qquad r_c/N\to c,\qquad r_b/N\to b.
\]
Use the direct-sum transfers of~\eqref{eq:sectors}.
For $0<h\le1$, let $\mathcal P_{N,h}$ be any partition of $[0,1]$
in the transverse coordinate $y$ into disjoint intervals
of length at most $h+O(N^{-1})$.
Then
\begin{equation}\label{eq:full-seam-insertion}
 \limsup_{N\to\infty}
 \sum_{I\in\mathcal P_{N,h}}
 \frac{\|\widehat{\mathcal T}_+\mathcal O_I
                     \widehat{\mathcal T}_-\|}{G_+^2G_-^2}
 \le C h\log(2/h).
\end{equation}
For fixed bias, $C$ depends only on the three cut levels. It is uniform
in $h$ and in the partition family $(\mathcal P_{N,h})_N$.
\end{proposition}

\begin{proof}
Choose $\varepsilon>0$ and a compact interval $A\Subset(0,1)$
such that, for all sufficiently large $N$,
\[
 |j-r_c|\le\varepsilon N\quad\Longrightarrow\quad j/N\in A.
\]
Set
\[
 v_N(h)=\sup_{j/N\in A}
          \sum_{I\in\mathcal P_{N,h}}\nu_{N,j}(I)^2.
\]

For each compatible quadruple $(A_-,B_-,A',B')$ in
Lemma~\ref{lem:insertion-weights}, use the corresponding
$v,R,T$ defined there.
By~\eqref{eq:partition-minors-square}
and~\eqref{eq:partition-minors-universal},
\[
 \sum_{I\in\mathcal P_{N,h}}|\det(P_I)_{R,T}|
 \le
 \begin{cases}
  v_N(h),
  & |j-r_c|\le\varepsilon N
    \text{ for every }j\in R\cup T,\\
  1,&\text{otherwise}.
 \end{cases}
\]

In the full direct sum, bound the operator norm by the sum of the
absolute values of all its matrix entries. Each input pair determines
its imbalance $s=|A_-|-r_a$, so each compatible quadruple occurs
exactly once across all sectors. Equation~\eqref{eq:normalized-insertion-entry} and the
weight bounds~\eqref{eq:insertion-weight-sums} give
\[
 \begin{aligned}
 &\sum_{I\in\mathcal P_{N,h}}
 \frac{\|\widehat{\mathcal T}_+\mathcal O_I
               \widehat{\mathcal T}_-\|}{G_+^2G_-^2}\\
 &\qquad\le
 \sum_{\substack{
   (A_-,B_-,A',B')\\
   \text{compatible}
 }}
 v\sum_{I\in\mathcal P_{N,h}}|\det(P_I)_{R,T}|\\
 &\qquad\le Cv_N(h)+Ce^{-c_1N}.
 \end{aligned}
\]
Taking the upper limit as $N\to\infty$ and applying
\eqref{eq:mode-partition-square} proves
\eqref{eq:full-seam-insertion}.
\end{proof}

\section{Opposite crossings on the entire section}
\label{sec:global-opposite}

To estimate opposite crossings on the middle section, fix one
interval $I$. The insertion operator of Section~\ref{sec:insertion}
gives the expectation of a signed count $X_I$ of pairs of sites on
different connecting paths. Lemma~\ref{lem:cut-signs} shows that
$X_I$ is nonnegative on every configuration. If two tested arcs of
one loop have opposite indices on $I$, Lemma~\ref{lem:detection}
supplies a positive lower bound for $X_I$ after truncating the number
of deep arches. The transfer estimate then allows us to sum over the
cells of an arbitrary partition in Theorem~\ref{thm:global-opposite}.

\subsection{Setting and the estimate for one interval}
\label{subsec:opposite-setting}

Keep $c$, $\rho_\pm$, the levels $a<c<b$ and the sections
$r_a<r_c<r_b$ of Section~\ref{sec:insertion}. Define the
\emph{outer sections}
\[
 r_-=\lceil N(c-\rho_-)\rceil,\qquad
 r_+=\lfloor N(c+\rho_+)\rfloor .
\]
Since $c-\rho_-<a$ and $b<c+\rho_+$, there is $\delta>0$ such that,
for all sufficiently large $N$,
\begin{equation}\label{eq:five-sections}
 1\le r_-<r_a<r_c<r_b<r_+\le N,\qquad
 r_a-r_-\ge\delta(N+1),\qquad r_+-r_b\ge\delta(N+1).
\end{equation}
We assume~\eqref{eq:five-sections} throughout this section.

\paragraph{Cuts.}
Split the graph at the three sections $r_a,r_c,r_b$ as in
Section~\ref{subsec:marked-caps}, and call the resulting cuts the
\emph{left cut}, \emph{middle cut} and \emph{right cut}. The left cap is the
part of the split graph to the left of the left cut, and the right cap
is the part to the right of the right cut. Full configurations $\omega$
of this split graph are defined as in
Section~\ref{subsec:deep-arches}, and
\eqref{eq:configuration-probability} holds. The loops of $\omega$ are
disjoint simple closed curves.

Deep arches, active loops and connecting paths are those of
Definition~\ref{def:deep-arches}, for the left and right cuts and the
outer sections $r_-,r_+$. As there, $n(\omega)$ is the number of
active loops and $d(\omega)$ the number of deep arches. Write
$\CP(\omega)$ for the set of connecting paths; it has $d(\omega)$
elements. Deep arches lie in the two caps, so they do not meet the
middle cut.

\paragraph{Deep regions and outer sections.}
In the diagonal coordinates of Section~\ref{subsec:seams} the white
vertices of section $r$ lie on the line $s=r/N$, and $s$ changes by
$1/(2N)$ along every edge. Hence an edge that meets
the left deep region $\{s\le c-\rho_-\}$ has its white endpoint on a
section of index at most $r_-$, and an edge that meets the right deep
region $\{s\ge c+\rho_+\}$ has its white endpoint on a section of index
at least $r_+$. Choose the splitting neighborhoods disjoint from the
two deep regions. We use the following consequence.
\begin{equation}\label{eq:deep-point}
 \parbox{0.86\linewidth}{If a loop crosses the left cut and contains a
 point $p$ of the left deep region, then $p$ lies on a deep left arch.
 The same holds on the right.}
\end{equation}
Indeed, the part of the loop in the left cap that contains $p$ is an
arch. Its ends are on section $r_a$, it contains a white vertex of a
section of index at most $r_-$, and the index of the section changes by
at most one between consecutive white vertices of a path. So the arch
visits section $r_-$.

\paragraph{Intervals on the middle cut.}
In the embedding of the split graph the middle cut is the line
$s=r_c/N$. It meets the graph only at the points $(r_c/N,x/N)$,
$x\in E_N$, where it crosses the cut edges. In this section an interval
$I$ of the transverse coordinate always stands for the segment
$\{(r_c/N,y):y\in I\}$ of the middle cut. We write $p_-$ and $p_+$ for
its lower and upper endpoint and assume that they are not on the
graph, that is, their transverse coordinates are not of the form
$x/N$ with $x\in E_N$.
As in Section~\ref{sec:insertion}, $S(I)=\{x\in E_N:x/N\in I\}$.
An oriented arc $\alpha$ crosses $I$ at the sites $x\in S(I)$ whose cut
edge it contains; the crossing is positive if $\alpha$ runs through
this edge from $w_L$ to $b_{\rm mid}$, that is, in the direction of
increasing $s$. The index $\iota_I(\alpha)$ of
Definition~\ref{def:tested-arcs} is the number of positive crossings
minus the number of negative ones.

\paragraph{The two events.}
For an interval $I$ put
\begin{equation}\label{eq:single-interval-bad-event}
 B_{N,I}=\bigl\{\text{some loop has two tested arcs }\alpha,\beta
 \text{ with }\iota_I(\alpha)\,\iota_I(\beta)<0\bigr\}.
\end{equation}
Both arcs inherit their direction from the same orientation of the
loop. They may overlap. For an integer $M\ge1$ let
\[
 \mathcal F_M=\{\text{each of the two caps has at most $M$ deep arches}\}.
\]

\begin{lemma}[One interval]\label{lem:single-interval-crossing}
For every integer $M\ge1$, every interval $I$, and all $N$
satisfying~\eqref{eq:five-sections},
\begin{equation}\label{eq:single-interval-crossing}
 \PP(B_{N,I}\cap \mathcal F_M)
 \le 4^{M+1}\,
 \frac{\|\widehat{\mathcal T}_+\mathcal O_I
                  \widehat{\mathcal T}_-\|}{G_+^2G_-^2}.
\end{equation}
\end{lemma}

The proof is given in Section~\ref{subsec:proof-single-interval}.

\subsection{The function \texorpdfstring{$X_I$}{X\_I} and its expectation}
\label{subsec:XI}

Fix a full configuration $\omega$. A site $x\in E_N$ of the middle cut
is \emph{singly occupied} if its cut edge belongs to exactly one of the
two perfect matchings. These are the sites at which the loops of $\omega$ cross
the middle cut. List \emph{all} singly occupied sites of the middle cut
in increasing order, $x_1<x_2<\cdots$, and put
\begin{equation}\label{eq:rank-sign}
 \operatorname{rk}(x_k)=k,\qquad \sigma(x)=(-1)^{\operatorname{rk}(x)} .
\end{equation}
The rank counts the sites outside $I$ and the sites on inactive loops
as well. For a loop or a path $\eta$ of $\omega$, let
\begin{equation}\label{eq:connector-sign-sum}
 t_\eta=\sum_{x}\sigma(x),
\end{equation}
the sum being over the singly occupied sites $x\in S(I)$ whose cut edge
lies on $\eta$. Put
\begin{equation}\label{eq:X-definition}
 t_{\rm tot}=\sum_{\eta\in\CP(\omega)}t_\eta,\qquad
 X_I(\omega)=2^{-n(\omega)-3}
   \Bigl(\sum_{\eta\in\CP(\omega)}t_\eta^2-t_{\rm tot}^2\Bigr).
\end{equation}

Let $\Theta$ be the endpoint rule of Definition~\ref{def:colour-tensor}
that assigns the averaging matrix $P$ to every deep arch and the
ordinary matrix $J$ to every other arch, and put
\[
 V_L=\widehat F_L(\Theta)\in\mathcal H_{r_a},\qquad
 V_R=\widehat F_R(\Theta)\in\mathcal H_{r_b},
\]
the colour tensors of the left cap at $r_a$ and of the right cap at
$r_b$.

\begin{lemma}[Operator identity]\label{lem:operator-identity}
For every interval $I$,
\begin{equation}\label{eq:contraction}
 \EE X_I=Z_N^{-2}\bigl\langle V_R,
   \widehat{\mathcal T}_+\mathcal O_I\widehat{\mathcal T}_-V_L
   \bigr\rangle .
\end{equation}
\end{lemma}

\begin{proof}
\emph{Expansion at the three cuts.}
At the left cut the gauge on $V_L$ cancels the inverse input gauge
of $\widehat{\mathcal T}_-$. At the right cut the inverse gauge on
$V_R$ cancels the output gauge of $\widehat{\mathcal T}_+$.
The middle gauges cancel because $\mathcal O_I$ commutes with the
two-replica gauge. We can thus expand~\eqref{eq:contraction}
using the ungauged cap tensors and layer weights, as in
Lemma~\ref{lem:sector-count}. A term consists of $i<j$ in
$S(I)$, superpositions $C_L,C_R$ of the two caps with compatible pairs
of states at the left and the right cut, a pair of states $(S_1,S_2)$
at the middle cut with $i,j\in S_1\setminus S_2$, and, for $\nu=1,2$,
a matching of the layers between the left and the middle cut with
output state $S_\nu$ and a matching of the layers between the middle
and the right cut with input state $S_\nu'$, where
\[
 S_1'=S_1\setminus\{i,j\},\qquad S_2'=S_2\cup\{i,j\}.
\]
The sites $i,j$ occur once in $(S_1,S_2)$ and once in $(S_1',S_2')$.
So the union of a term is again a full configuration $\omega$, in
which $i$ and $j$ are singly occupied sites of the middle cut. The
terms with a given $\omega$ and $(i,j)$ are the \emph{labellings of
$\omega$ with insertions at $i$ and $j$}. These are defined as in
Definition~\ref{def:labellings}, with one change: the cut edges of the
middle cut at $i$ and $j$ have label $1$ on their left half, attached
to $w_L$, because $i,j\in S_1$, and label $2$ on their right half,
attached to $b_{\rm mid}$, because $i,j\in S_2'$; rule~(R1) is applied
to these half-edges. A term contributes
\begin{equation}\label{eq:term-weight}
 \epsilon\;\bk^{\dn(\omega)}\;2^{\ell(C_L)+\ell(C_R)}\;2^{-d(\omega)},
\end{equation}
where $\epsilon=\pm1$ is the sign of the matrix entry of
$c_ic_j\otimes c_j^\dagger c_i^\dagger$ between
$\varepsilon_{S_1}\otimes\varepsilon_{S_2}$ and
$\varepsilon_{S_1'}\otimes\varepsilon_{S_2'}$, and $2^{-d(\omega)}$ is
the product of the factors $\tfrac12$ of the deep arches.

\emph{Sign of the insertion.}
From the definition of $c_x$ and $c_x^\dagger$, for $i<j$ in
$S_1\setminus S_2$,
\[
 c_ic_j\,\varepsilon_{S_1}
 =(-1)^{1+\#\{s\in S_1:\,i<s<j\}}\varepsilon_{S_1'},\qquad
 c_j^\dagger c_i^\dagger\,\varepsilon_{S_2}
 =(-1)^{1+\#\{s\in S_2:\,i<s<j\}}\varepsilon_{S_2'} .
\]
Hence $\epsilon=(-1)^{\mu}$, where $\mu$ is the number of sites strictly
between $i$ and $j$ counted with their multiplicity in $(S_1,S_2)$.
Doubly occupied sites contribute an even number, so $\mu$ has the
parity of the number of singly occupied sites strictly between $i$ and
$j$, which is $\operatorname{rk}(j)-\operatorname{rk}(i)-1$. Therefore
\begin{equation}\label{eq:middle-cut-sign}
 \epsilon=-\sigma(i)\sigma(j),
\end{equation}
and $\epsilon$ depends only on $\omega$ and $(i,j)$.

\emph{Counting the labellings.}
Direct every labelled edge and half-edge from its white to its black
endpoint if its label is $1$, and from black to white if its label is
$2$. By~(R1) and~(R2) the arrows keep one direction along a connecting
path or an inactive loop, except at $i$ and $j$. There the left half
of the cut edge has label $1$ and starts at the white vertex $w_L$,
the right half has label $2$ and starts at the black vertex
$b_{\rm mid}$, so both halves point towards the cut
(Figure~\ref{fig:insertion-site}). No arc of a path or loop can be
directed towards both of its ends. Hence an inactive loop contains
neither $i$ nor $j$, a connecting path contains at most one of them,
and a connecting path that contains one of them has exactly one
labelling, in which both of its parts point towards that site. By
Lemma~\ref{lem:labellings}(i), each of the other connecting paths has
two labellings, and each inactive loop contributes the factor $2$, as
in the proof of Lemma~\ref{lem:sector-count}. So the sum of the
contributions~\eqref{eq:term-weight} of all terms with a given
$\omega$ and $(i,j)$ is zero unless $i$ and $j$ lie on two distinct
connecting paths, and in that case it equals
\begin{equation}\label{eq:fixed-configuration-count}
 \epsilon\,\bk^{\dn(\omega)}\;
 2^{\ell(\omega)-n(\omega)}\cdot2^{d(\omega)-2}\cdot2^{-d(\omega)}
 =-\sigma(i)\sigma(j)\;2^{-n(\omega)-2}\;
   \bk^{\dn(\omega)}2^{\ell(\omega)} .
\end{equation}
Here $2^{\ell(\omega)-n(\omega)}$ counts the two choices on each
inactive loop, including the loops already counted in the caps.
The $d(\omega)-2$ connecting paths free of insertions give
$2^{d(\omega)-2}$, and the $d(\omega)$ deep arches each contribute
the averaging weight $1/2$.

\begin{figure}[htbp]
 \centering
 \begin{tikzpicture}[x=1.45cm,y=1cm,
   wv/.style={circle,draw,fill=white,inner sep=1.8pt,line width=0.5pt},
   bv/.style={circle,draw,fill=black,inner sep=1.8pt,line width=0.5pt},
   one/.style={line width=1.5pt,blue!75!black},
   two/.style={line width=1.5pt,red!70!black,densely dashed},
   mid/.style={postaction={decorate},
     decoration={markings,mark=at position 0.62 with
       {\arrow{Stealth[length=6pt,width=6pt]}}}}]
  \begin{scope}[yshift=0cm]
   \draw[two,mid] (0,0)--(1,0);
   \draw[one,mid] (1,0)--(2.4,0);
   \draw[two,mid] (2.4,0)--(3.4,0);
   \draw[one,mid] (3.4,0)--(4.4,0);
   \node[bv] at (0,0) {};
   \node[wv,label=below:{\small $w_L$}] at (1,0) {};
   \node[bv,label=below:{\small $b_{\rm mid}$}] at (2.4,0) {};
   \node[wv,label=below:{\small $w_R$}] at (3.4,0) {};
   \node[bv] at (4.4,0) {};
   \node[above] at (0.5,0.05) {\small 2};
   \node[above] at (1.35,0.05) {\small 1};
   \node[above] at (2.05,0.05) {\small 1};
   \node[above] at (2.9,0.05) {\small 2};
   \node[above] at (3.9,0.05) {\small 1};
   \node[anchor=east] at (-0.35,0) {\small ordinary site};
  \end{scope}
  \begin{scope}[yshift=-1.9cm]
   \draw[two,mid] (0,0)--(1,0);
   \draw[one,mid] (1,0)--(1.7,0);
   \draw[two,mid] (2.4,0)--(1.7,0);
   \draw[one,mid] (3.4,0)--(2.4,0);
   \draw[two,mid] (4.4,0)--(3.4,0);
   \node[bv] at (0,0) {};
   \node[wv,label=below:{\small $w_L$}] at (1,0) {};
   \node[bv,label=below:{\small $b_{\rm mid}$}] at (2.4,0) {};
   \node[wv,label=below:{\small $w_R$}] at (3.4,0) {};
   \node[bv] at (4.4,0) {};
   \node[above] at (0.5,0.05) {\small 2};
   \node[above] at (1.35,0.05) {\small 1};
   \node[above] at (2.05,0.05) {\small 2};
   \node[above] at (2.9,0.05) {\small 1};
   \node[above] at (3.9,0.05) {\small 2};
   \node[anchor=east] at (-0.35,0) {\small insertion site $i$ or $j$};
  \end{scope}
  \draw[dashed,line width=0.6pt] (1.7,0.75)--(1.7,-2.65);
  \node[above] at (1.7,0.75) {\small middle cut};
 \end{tikzpicture}
 \caption{A loop passing through a singly occupied site of the middle
 cut. Labels $1$ and $2$ name the matching that contains the edge.
 Edges with label $1$ are directed from white to black, edges with
 label $2$ from black to white. At an ordinary site (top) the whole cut
 edge has one label and the arrows run through. At an insertion site
 (bottom) the two halves of the cut edge have labels $1$ and $2$, and
 the arrows on both sides point towards the cut.}
 \label{fig:insertion-site}
\end{figure}
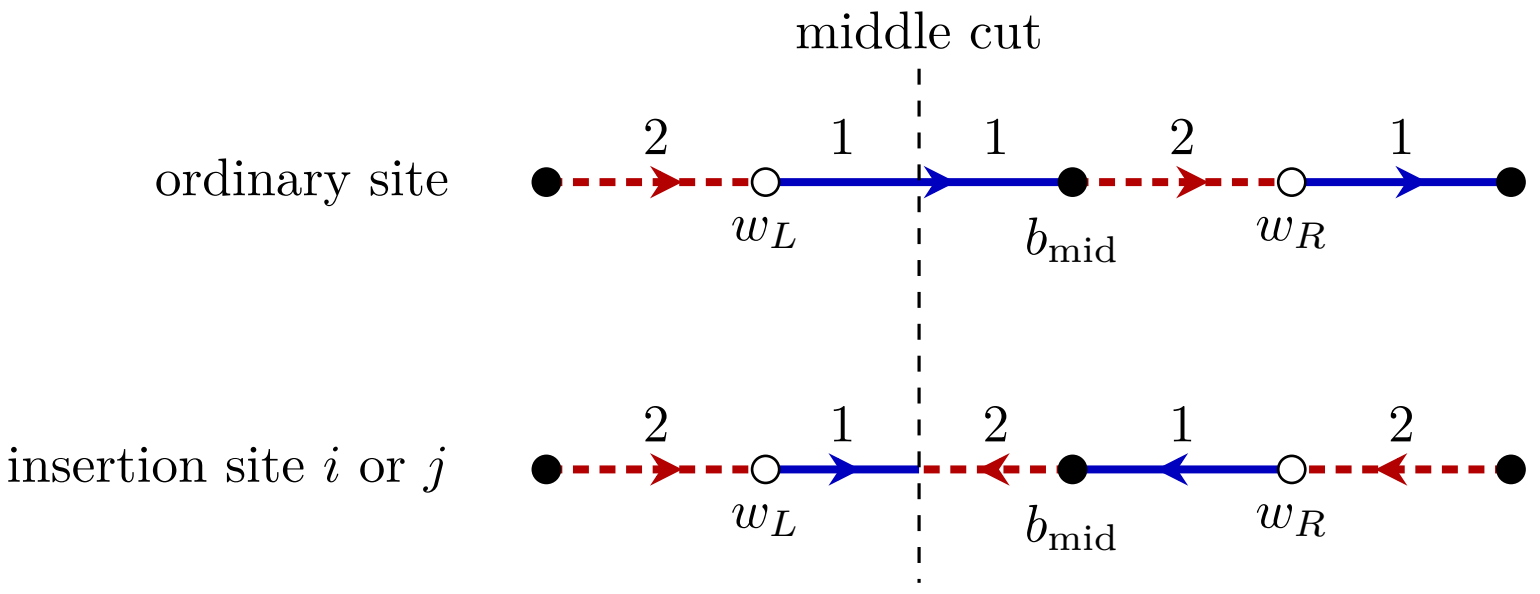

\emph{Summation over the two sites.}
Sum~\eqref{eq:fixed-configuration-count} over the pairs $i<j$ in
$S(I)$ that lie on distinct connecting paths. By the
definition~\eqref{eq:connector-sign-sum} of $t_\eta$,
\[
 -\sum_{\{\eta,\eta'\}}t_\eta t_{\eta'}
 =\frac12\Bigl(\sum_{\eta\in\CP(\omega)}t_\eta^2-t_{\rm tot}^2\Bigr),
\]
where the sum on the left is over unordered pairs of distinct
connecting paths. With~\eqref{eq:configuration-probability} this gives
\[
 Z_N^{-2}\bigl\langle V_R,
   \widehat{\mathcal T}_+\mathcal O_I\widehat{\mathcal T}_-V_L
   \bigr\rangle
 =\sum_\omega\PP\{\omega\}\,2^{-n(\omega)-3}
   \Bigl(\sum_{\eta\in\CP(\omega)}t_\eta^2-t_{\rm tot}^2\Bigr)
 =\EE X_I .\qedhere
\]
\end{proof}

\subsection{Signs along the middle cut}
\label{subsec:cut-signs}

For a loop $g$ of $\omega$ write $\operatorname{Int}(g)$ for the
bounded component of its complement and $\nest(g)$ for the number of
loops of $\omega$ that enclose $g$, that is, the number of loops
$g'\ne g$ with $g\subset\operatorname{Int}(g')$. Move along the middle
cut in the direction of increasing $y$. At a site of $g$ the cut either
\emph{enters} $g$ (passes from the outside of $g$ to
$\operatorname{Int}(g)$) or \emph{leaves} $g$.

\begin{lemma}[Signs along the cut]\label{lem:cut-signs}
Let $\omega$ be a full configuration and $I$ an interval with
endpoints $p_-,p_+$.
\begin{enumerate}[label=\textnormal{(\roman*)},leftmargin=*]
\item If $x$ is a site of a loop $g$, then
$\sigma(x)=-(-1)^{\nest(g)}$ if the cut enters $g$ at $x$, and
$\sigma(x)=(-1)^{\nest(g)}$ if it leaves $g$ at $x$.
\item For every loop $g$,
\begin{equation}\label{eq:loop-rank-sign}
 t_g=-(-1)^{\nest(g)}
 \bigl(\one_{\{p_+\in\operatorname{Int}(g)\}}
      -\one_{\{p_-\in\operatorname{Int}(g)\}}\bigr).
\end{equation}
\item $t_{\rm tot}\in\{-1,0,1\}$.
\item The number
$\sum_{\eta\in\CP(\omega)}t_\eta^2-t_{\rm tot}^2$ is a nonnegative even
integer. Consequently $X_I(\omega)\ge0$, and
\begin{equation}\label{eq:positivity}
 X_I(\omega)\ge2^{-n(\omega)-2}
 \qquad\text{whenever}\quad\sum_{\eta\in\CP(\omega)}t_\eta^2\ge2 .
\end{equation}
\end{enumerate}
\end{lemma}

Figure~\ref{fig:cut-signs} shows an example.

\begin{figure}[htbp]
 \centering
 \begin{tikzpicture}[x=0.8cm,y=1cm,
   lp/.style={line width=0.9pt},
   st/.style={circle,fill=black,inner sep=1.3pt}]
  \draw[lp] (0,3) ellipse [x radius=4.6,y radius=3.4];
  \draw[lp] (0,1.0) ellipse [x radius=1.5,y radius=0.6];
  \draw[lp] (0,2.6) ellipse [x radius=1.0,y radius=0.45];
  \draw[lp] (0,4.9) ellipse [x radius=2.6,y radius=1.1];
  \draw[lp] (0,5.0) ellipse [x radius=1.5,y radius=0.55];
  \draw[line width=0.5pt,-{Stealth[length=5pt]}] (0,-1.0)--(0,7.2)
     node[above] {\small $y$};
  \draw[line width=2.6pt,black!40] (0,1.2)--(0,4.8);
  \draw[line width=0.7pt] (-0.16,1.2)--(0.16,1.2);
  \draw[line width=0.7pt] (-0.16,4.8)--(0.16,4.8);
  \node[anchor=east] at (-0.18,1.2) {\small $p_-$};
  \node[anchor=east] at (-0.18,4.8) {\small $p_+$};
  \node[anchor=east,text=black!55] at (-0.18,3.3) {\small $I$};
  \foreach \yy/\sg/\col in {-0.4/-/black!50,0.4/+/black!50,1.6/-/black,
     2.15/+/black,3.05/-/black,3.8/+/black,4.45/-/black,
     5.55/+/black!50,6.0/-/black!50,6.4/+/black!50} {
    \node[st] at (0,\yy) {};
    \node[anchor=west,text=\col] at (0.14,\yy) {\small $\sg$};
  }
  \node[anchor=west] at (4.65,3.0) {\small $G$};
  \node[anchor=west] at (1.55,1.0) {\small $g_-$};
  \node[anchor=west] at (1.05,2.6) {\small $h$};
  \node[anchor=west] at (2.65,4.9) {\small $g_1$};
  \node[anchor=west] at (1.55,5.0) {\small $g_2$};
  \node[anchor=west,align=left] at (6.3,3.3) {\small
   $\begin{array}{@{}l@{\quad}c@{\quad}r@{}}
     \text{loop} & \nest & t_g\\[2pt]
     G & 0 & 0\\
     g_- & 1 & -1\\
     h & 1 & 0\\
     g_1 & 1 & +1\\
     g_2 & 2 & -1
    \end{array}$};
 \end{tikzpicture}
 \caption{Five loops crossing the middle cut, the signs $\sigma(x)$ at
 their ten sites (ranks increase upwards), and the interval $I$ in gray.
 The values $t_g$ add the signs at the sites of $g$ inside $I$. The
 loops with $t_g\ne0$ are those that separate $p_-$ from $p_+$: the
 chain $g_1\supset g_2$ around $p_+$, with alternating values
 $+1,-1$, and the chain $g_-$ around $p_-$, whose outermost value $-1$
 is opposite to that of $g_1$.}
 \label{fig:cut-signs}
\end{figure}

\begin{proof}
(i) Each loop contributes an odd number of crossings below a point
of the cut exactly when it contains that point. Just below a crossing
$x$ at which the cut enters $g$, the surrounding loops are precisely
those enclosing $g$. Just below an exit they are these loops and $g$
itself. Thus
$\operatorname{rk}(x)-1$ has the parity of $\nest(g)$ at an entry
and of $\nest(g)+1$ at an exit.

(ii) Between $p_-$ and $p_+$, the number of entries into $g$ minus the
number of exits equals the change in its interior indicator. Part~(i)
therefore gives~\eqref{eq:loop-rank-sign}.

(iii) Every site of an active loop lies on exactly one connecting
path, because deep arches do not meet the middle cut. Hence
\begin{equation}\label{eq:active-sign-total}
 t_{\rm tot}=\sum_{g\ \mathrm{active}}t_g .
\end{equation}
By~(ii), $t_g\ne0$ only if $g$ belongs to one of the two families
\[
 \mathcal G_+=\{g:\ p_+\in\operatorname{Int}(g),\
                   p_-\notin\operatorname{Int}(g)\},\qquad
 \mathcal G_-=\{g:\ p_-\in\operatorname{Int}(g),\
                   p_+\notin\operatorname{Int}(g)\},
\]
and then $t_g=-(-1)^{\nest(g)}$ for $g\in\mathcal G_+$ and
$t_g=(-1)^{\nest(g)}$ for $g\in\mathcal G_-$.
Each family is a chain under inclusion of loop interiors. For
consecutive loops in either chain, the nesting depths differ by one:
every intermediate loop would separate the endpoints in the same
way. Thus the nonzero $t_g$ alternate in sign along each chain.
If both chains are nonempty, their outermost loops have the same
nesting depth, since the loops enclosing either one are precisely
those containing both $p_+$ and $p_-$. The two initial signs are
therefore opposite.

In each chain the active loops form an initial segment when ordered
from the outside in. Indeed, suppose $g$ encloses an active loop in
its chain with a deep left arch. The inner loop reaches section
$r_-$, so the minimum of $s$ on the closed disk bounded by $g$
is at most $r_-/N$ and is attained on $g$. On the other hand,
$g$ crosses $I$ because it separates $p_-$ from $p_+$. Hence it
crosses the left cut and reaches an edge whose white endpoint has
section index at most $r_-$. The argument for~\eqref{eq:deep-point}
gives a deep left arch of $g$. The same argument with the maximum
of $s$ treats a deep right arch.

The active part of each alternating chain sums to zero or to its
initial sign. The initial signs are opposite when both chains are
present, so~\eqref{eq:active-sign-total} gives~(iii).

(iv) The $t_\eta$ are integers. If $t_{\rm tot}=\pm1$, some $t_\eta$ is
nonzero, so $\sum_\eta t_\eta^2\ge1=t_{\rm tot}^2$. Since
$k^2\equiv k\pmod2$ for integers $k$,
$\sum_\eta t_\eta^2\equiv t_{\rm tot}\equiv t_{\rm tot}^2\pmod2$.
If $\sum_\eta t_\eta^2\ge2$, the difference is at least $1$ and even,
hence at least $2$.
\end{proof}

\subsection{Detection of opposite crossings}
\label{subsec:detection}

\begin{lemma}[Detection]\label{lem:detection}
If $\omega\in B_{N,I}$, then
$\sum_{\eta\in\CP(\omega)}t_\eta^2\ge2$.
\end{lemma}

\begin{proof}
Let $\alpha,\beta$ be tested arcs of a loop $g$, directed along the
same orientation of $g$, with $\iota_I(\alpha)\iota_I(\beta)<0$. Put
$e=\iota_I(g)\in\{-1,0,1\}$. One of the two arcs, say $\alpha$, has
index $m=\iota_I(\alpha)\notin\{0,e\}$, since the two indices are
nonzero and different. Let $\alpha'$ be the complementary arc of
$\alpha$ in $g$, with the same orientation. Then
$\iota_I(\alpha')=e-m$, and for an integer $m\notin\{0,e\}$ with
$|e|\le1$ the numbers $m$ and $e-m$ are nonzero and have opposite
signs.

\emph{Separating the connecting paths.}
The loop $g$ contains points of both deep regions, so it crosses the
left and the right cut. By~\eqref{eq:deep-point} the endpoint of
$\alpha$ in the left deep region lies on a deep left arch $A_1$, and
its other endpoint lies on a deep right arch $A_2$. Removing $A_1$ and
$A_2$ from $g$ leaves two arcs $g_1,g_2$, with
\[
 g_1\subset\alpha\subset g_1\cup A_1\cup A_2,\qquad
 g_2\subset\alpha'\subset g_2\cup A_1\cup A_2 .
\]
Deep arches miss the middle cut. Hence $\alpha$ and $g_1$ meet the
cut at the same sites, as do $\alpha'$ and $g_2$. Each of $g_1,g_2$ is a
union of connecting paths and deep arches; let $\CP_1$ and $\CP_2$ be
the sets of connecting paths contained in $g_1$ and in $g_2$. They are
disjoint.

\emph{Crossing directions and the signs $\sigma$.}
Orient $g$ counterclockwise in the $(s,y)$ plane. If it crosses the cut at $x$ in
the direction of increasing $s$, then $\operatorname{Int}(g)$, which
is on the left of $g$, lies above $x$ on the cut, so the cut enters $g$
at $x$. If $g$ crosses in the direction of decreasing $s$, the cut
leaves $g$ at $x$. For the clockwise orientation the two cases are
exchanged. With Lemma~\ref{lem:cut-signs}(i) this shows that there is
$\varepsilon_g\in\{-1,1\}$, depending only on $g$ and its orientation,
such that the crossing of $g$ at $x$ is positive if and only if
$\sigma(x)=\varepsilon_g$. Consequently
\begin{equation}\label{eq:bad-connector-sums}
 \sum_{\eta\in\CP_1}t_\eta=\varepsilon_g\,m,\qquad
 \sum_{\eta\in\CP_2}t_\eta=\varepsilon_g\,(e-m).
\end{equation}

The two sums in~\eqref{eq:bad-connector-sums} are nonzero and have
opposite signs. So some $t_\eta$ with $\eta\in\CP_1$ and some $t_\eta$
with $\eta\in\CP_2$ are nonzero, and $\sum_\eta t_\eta^2\ge2$.
\end{proof}

\subsection{Proof of Lemma~\ref{lem:single-interval-crossing}}
\label{subsec:proof-single-interval}

On $\mathcal F_M$ there are at most $2M$ deep arches, and every active loop
contains one, so $n(\omega)\le2M$. By Lemma~\ref{lem:detection}
and~\eqref{eq:positivity}, $X_I\ge2^{-2M-2}=4^{-M-1}$ on
$B_{N,I}\cap \mathcal F_M$, and $X_I\ge0$ everywhere by
Lemma~\ref{lem:cut-signs}(iv). Therefore, by
Lemma~\ref{lem:operator-identity} and the Cauchy--Schwarz inequality,
\[
 4^{-M-1}\,\PP(B_{N,I}\cap \mathcal F_M)\le\EE X_I
 \le\frac{\|V_R\|\,\|V_L\|}{Z_N^2}\,
    \|\widehat{\mathcal T}_+\mathcal O_I\widehat{\mathcal T}_-\| .
\]
The matrices $P$ and $J$ have operator norm one, so
Lemma~\ref{lem:aztec-marked-cap} gives
\begin{equation}\label{eq:opposite-cap-bounds}
 \|V_L\|\le\|\widehat L_{r_a}\|^2,\qquad
 \|V_R\|\le\|\widehat R_{r_b}\|^2 .
\end{equation}
By the normalization~\eqref{eq:normalization},
$\|\widehat L_{r_a}\|^2\|\widehat R_{r_b}\|^2/Z_N^2=G_{r_a,r_b}^{-2}$,
and $G_{r_a,r_b}=G_-G_+$. This
proves~\eqref{eq:single-interval-crossing}. \qed

\subsection{All intervals of a partition}
\label{subsec:all-intervals}

Fix a bounded interval $J_0$ of the transverse coordinate that contains
$[0,1]$ in its interior, and a constant $C_0\ge0$. For $0<h\le1$ let
$\mathfrak P_{N,h}$ be the set of all finite partitions of $J_0$ into
intervals of length at most $h+C_0/N$ whose endpoints $y$ satisfy
$Ny\notin E_N$. For $\mathcal P\in\mathfrak P_{N,h}$ put
\[
 B^{\rm all}_{N}(\mathcal P)=\bigcup_{I\in\mathcal P}B_{N,I}.
\]

\begin{theorem}[Opposite crossings on the entire section]
\label{thm:global-opposite}
There are constants $C,c_0>0$, depending on $c,\rho_\pm,a,b,J_0,C_0$
and the bias, such that for $0<h\le1$ and every integer $M\ge1$,
\begin{equation}\label{eq:global-opposite-truncated}
 \limsup_{N\to\infty}\sup_{\mathcal P\in\mathfrak P_{N,h}}
 \PP\bigl(B^{\rm all}_{N}(\mathcal P)\bigr)
 \le Ce^{-c_0 M^2}+C\,4^M h\log(2/h).
\end{equation}
Consequently, for every $0<\alpha<1$ there is $C_\alpha$ such that for
all sufficiently small $h$,
\begin{equation}\label{eq:global-opposite}
 \limsup_{N\to\infty}\sup_{\mathcal P\in\mathfrak P_{N,h}}
 \PP\bigl(B^{\rm all}_{N}(\mathcal P)\bigr)
 \le C_\alpha h^\alpha .
\end{equation}
The same estimates hold for the reflected cuts of
Section~\ref{subsec:reflection}.
\end{theorem}

\begin{proof}
\emph{Many deep arches are unlikely.}
By Lemma~\ref{lem:deep-traversals}(i), a cap with more than $M$ deep
arches produces at least $2(M+1)$ traversals of its outer slab, between
the sections $r_-,r_a$ or between $r_b,r_+$.
By~\eqref{eq:five-sections} both slabs have width at least
$\delta(N+1)$, so~\eqref{eq:gaussian-traversal-integer} gives
\begin{equation}\label{eq:many-deep-arches}
 \PP(\mathcal F_M^c)\le P_N^{\mathrm{tr}}(M+1;r_-,r_a)
             +P_N^{\mathrm{tr}}(M+1;r_b,r_+)
 \le Ce^{-c_0M^2}.
\end{equation}

\emph{Summation over the intervals.}
Fix $M$ and $h$. For each $N$ choose
$\mathcal P_N\in\mathfrak P_{N,h}$ such that
$\PP(B^{\rm all}_{N}(\mathcal P_N)\cap \mathcal F_M)$ is within $N^{-1}$ of its
supremum over $\mathfrak P_{N,h}$. The operator $\mathcal O_I$ depends
on $I$ only through $S(I)$, and $S(I\cap[0,1])=S(I)$. The nonempty
sets $I\cap[0,1]$, $I\in\mathcal P_N$, form a partition of $[0,1]$
into intervals of length at most $h+C_0/N$. Summing
\eqref{eq:single-interval-crossing} over $I\in\mathcal P_N$ and
applying Proposition~\ref{prop:full-seam-insertion} gives
\[
 \limsup_{N\to\infty}\sup_{\mathcal P\in\mathfrak P_{N,h}}
 \PP\bigl(B^{\rm all}_{N}(\mathcal P)\cap \mathcal F_M\bigr)
 \le 4^{M+1}\,C\,h\log(2/h).
\]
Adding~\eqref{eq:many-deep-arches}, once for the whole partition,
proves~\eqref{eq:global-opposite-truncated}.

\emph{Choice of $M$.}
Choose $A$ with $c_0A^2\ge1$ and set
$M=\lceil A\sqrt{\log(1/h)}\,\rceil$. Then $e^{-c_0M^2}\le h$.
For each $\epsilon>0$,
\[
 4^M\le 4\exp\!\bigl(A(\log 4)\sqrt{\log(1/h)}\bigr)
      \le C_\epsilon h^{-\epsilon},
 \qquad
 \log(2/h)\le C_\epsilon h^{-\epsilon}.
\]
Taking $\epsilon=(1-\alpha)/2$ proves~\eqref{eq:global-opposite}.

\emph{Reflected cuts.}
Reflection preserves the joint law of the two perfect matchings by
Section~\ref{subsec:reflection}, so the estimates also hold for the
reflected cuts.
\end{proof}

\section{Finite-puncture observables}
\label{sec:punctures}

To identify the law of the loop homotopy classes, we study the
expected product of normalized monodromy traces. Kenyon's determinant
identity reduces this observable to a matrix indexed by edges crossed
by the cuts. For unipotent monodromies the same-cut blocks drop out
exactly, so only transfers between distinct cuts enter the determinant
in Section~\ref{sec:interior-boundary}.

Use the liquid region $L$ and its intrinsic coordinate $\Phi$ from
Section~\ref{subsec:liquid-coordinates}.
Fix distinct punctures $z_1,\ldots,z_m\in L$ and put
$\zeta_i=\Phi(z_i)$. Choose deterministic representatives
$z_{i,N}^0$ off the scaled embedded graph such that
\[
 \max_i|z_{i,N}^0-z_i|\le C/N,
\]
taking $z_{i,N}^0=z_i$ whenever $z_i$ is off the graph.
Isotopies supported in pairwise disjoint neighborhoods in $L$
identify these punctures and the associated representation bases with
those at $z_1,\ldots,z_m$.

For a loop collection $\omega$, write
$\Lam_{p_1,\ldots,p_m}(\omega)$ for the free homotopy classes of its
loops enclosing at least two of the marked points. The reduction
retains nesting and multiplicities. In particular,
\[
 \Gamma_N=\Lam_{z_{1,N}^0,\ldots,z_{m,N}^0}(\mathcal L_N)
\]
takes values in a fixed space of laminations. The map $\Phi$
identifies lamination classes in $L$ with those in $\DD$.

For sufficiently large $N$, choose pairwise disjoint dual paths,
or \emph{zippers}, in the scaled embedding $N^{-1}G_N$.
They run from the faces containing $z_{i,N}^0$ to distinct points
on the outer boundary. Orient each zipper toward the
unbounded face and write $g_i$ for the corresponding counterclockwise
generator about $z_i$. In the oriented $z$-plane a crossing is
positive if it passes from the right side of the oriented zipper to
its left. Let $X_{\rm unip}$ denote the variety of
$\mathrm{SL}_2(\mathbb C)$ representations with unipotent monodromy
around each puncture. Thus for $\rho\in X_{\rm unip}$,
\begin{equation}\label{eq:unipotent}
 M_i=\rho(g_i)=I_2+A_i,\qquad A_i^2=0,
 \qquad 1\le i\le m.
\end{equation}
Different $A_i$ need not commute. Since $M_i^{-1}=I_2-A_i$, a crossing
of zipper $i$ with sign $\sigma=\pm1$ contributes $I_2+\sigma A_i$.

The monodromy of a discrete loop $\gamma$ is the ordered product of
these matrices along its crossings. Set
\begin{equation}\label{eq:topological-correlator}
 F_N(\rho)=\EE\prod_{\gamma\in\mathcal L_N}
                  \frac{\Tr\rho(\gamma)}2.
\end{equation}
Contractible loops and loops around one puncture have factor $1$.
Also $\Tr M^{-1}=\Tr M$ for $M\in\mathrm{SL}_2$, so loop orientations
do not affect the observable. The integrand is therefore a function
of $\Gamma_N$. The proof of Theorem~\ref{thm:aztec-finite-puncture}
recovers the lamination law by varying $\rho$.

\subsection{Compression to zipper edges}
\label{subsec:zipper-compression}

Let $K_N$ be the scalar Kasteleyn matrix, with rows
indexed by white vertices and columns by black vertices.
In the unscaled lattice coordinates $(X,Y)$, set
\begin{equation}\label{eq:kasteleyn-definition}
 K_N(w,b)=
 \begin{cases}
  -\bk,& b-w=(0,1),\\
  1,& b-w\in\{(1,0),(-1,0),(0,-1)\},\\
  0,& \text{otherwise}.
 \end{cases}
\end{equation}
Since $b-w=(0,1)$ on a down edge, these moduli are precisely the
weights of the measure proportional to $\bk^{\dn(M)}$ in
Section~\ref{subsec:aztec-model}. The Kasteleyn signs give
$|\det K_N|=Z_N$.
The unscaled diagonal lattice coordinates are
\begin{equation}\label{eq:diagonal-coordinates}
 (u,v)=(X+Y+N,X-Y+N).
\end{equation}
For zipper $i$, let $E_i$ be the set of crossed primal edges.
Write $e=(w_e,b_e)$ and $i(e)=i$ for $e\in E_i$.
Set $\sigma_e=+1$ if $w_e$ lies to the right of the
oriented zipper in the $z$-plane,
and $\sigma_e=-1$ if it lies to the left.
Define the block matrix $K_{N,\rho}$ by
\begin{equation}\label{eq:twisted-perturbation}
 K_{N,\rho}
 =K_N\otimes I_2
 +\sum_{e\in E}|w_e\rangle
       \sigma_e K_N(w_e,b_e)A_{i(e)}\langle b_e|,
 \qquad E=\bigsqcup_{i=1}^m E_i.
\end{equation}
Here $|w\rangle C\langle b|$ denotes the block matrix
whose $(w,b)$ block is $C$ and whose other blocks are zero.
For two independent perfect matchings with these weights, Kenyon's
double-dimer identity in bipartite block form gives
\begin{equation}\label{eq:kenyon-determinant}
 F_N(\rho)=\frac{\det K_{N,\rho}}{(\det K_N)^2}.
\end{equation}
See \cite{Kenyon} and \cite[Lemma 2.10]{BI}.

Factor the perturbation in~\eqref{eq:twisted-perturbation} as $UV$,
where
\[
 U_{w,e}=\one_{w=w_e}\,\sigma_e K_N(w_e,b_e)A_{i(e)},
 \qquad V_{e,b}=\one_{b=b_e}I_2.
\]
Sylvester's identity $\det(I+PQ)=\det(I+QP)$ now yields
\begin{equation}\label{eq:edge-determinant}
 F_N(\rho)=\det(I+B_N),\qquad
 (B_N)_{ef}=K_N^{-1}(b_e,w_f)\,
       \sigma_f K_N(w_f,b_f)A_{i(f)}.
\end{equation}
Define scalar blocks
\begin{equation}\label{eq:cross-kernel}
 C_{ij,N}(e,f)=K_N^{-1}(b_e,w_f)\,
                     \sigma_f K_N(w_f,b_f),
 \qquad e\in E_i,\ f\in E_j.
\end{equation}
Thus $(B_N)_{ij}=C_{ij,N}\otimes A_j$: the matrix factor is indexed
by the column zipper $j$.

\paragraph{Conversion of scalar Kasteleyn conventions.}
Let $K_N^{\rm CJY}$ denote the matrix
in~\cite[equation~(2.11)]{CJY}, with bias $a=\sqrt\bk$ and coordinates
identified by~\eqref{eq:diagonal-coordinates}. It has black rows and
white columns, and its entries have modulus $a$ on the vertical edges
and modulus $1$ on the horizontal edges. Writing $(X_v,Y_v)$ for
the unscaled coordinates of a vertex $v$, put
\[
 \vartheta(v)=\ii^{X_v+1/2},\qquad
 h_b=(-1)^{X_b+N-1/2}.
\]
The four edge directions give
\begin{equation}\label{eq:CJY-phase-conversion}
 \begin{aligned}
 K_N(w,b)
 &=-\ii\,\vartheta(w)h_b\vartheta(b)^{-1}a^{\,Y_b-Y_w}
       K_N^{\rm CJY}(b,w),\\
 K_N^{-1}(b,w)
 &=\ii\,\vartheta(b)h_b\vartheta(w)^{-1}a^{\,Y_w-Y_b}
       (K_N^{\rm CJY})^{-1}(w,b).
 \end{aligned}
\end{equation}
Here the transpose acts on the scalar matrix and leaves the
transport matrices $A_i$ untouched. The positive gauge factor
$a^{Y_b-Y_w}$ transfers the vertical weight $a$ to a down-edge
weight $a^2=\bk$. Formula~\eqref{eq:CJY-phase-conversion} identifies
the conventions in \cite{CJY}. Lemma~\ref{lem:scalar-exact-inverse}
verifies the scalar inverse used in this paper directly.

Let $D_{\rm W}$ and $D_{\rm B}$ be diagonal matrices
indexed by white and black vertices, respectively,
with diagonal entries of modulus one. Set
$K_N'=D_{\rm W}K_ND_{\rm B}$.
Let $B_N'$ be defined by~\eqref{eq:edge-determinant}
with $K_N$ replaced by $K_N'$ and the transport
matrices unchanged. Then
\begin{equation}\label{eq:Kasteleyn-gauge-covariance}
 B_N'=D_E^{-1}B_ND_E,
 \qquad
 (D_E)_{ef}
 =\one_{\{e=f\}}D_{\rm B}(b_e,b_e)I_2.
\end{equation}
The matrix $D_E$ is unitary.

\subsection{Exact removal of same-zipper blocks}

\begin{proposition}\label{prop:nilpotent-reduction}
Let $B_N$ be the matrix in~\eqref{eq:edge-determinant}, and set
\[
 D_N=\operatorname{diag}((B_N)_{11},\ldots,(B_N)_{mm}),
 \qquad O_N=B_N-D_N.
\]
For every $N$ for which the zippers have been chosen and every
$\rho\in X_{\rm unip}$,
\begin{equation}\label{eq:nilpotent-reduction}
 F_N(\rho)=\det(I+O_N)=\det{}_2(I+O_N).
\end{equation}
Here $\det{}_2$ denotes the regularized Fredholm determinant;
in finite dimension,
\[
 \det{}_2(I+O_N)=\det(I+O_N)e^{-\Tr O_N}.
\]
\end{proposition}

\begin{proof}
Since $A_i^2=0$, each diagonal block of $D_N^2$ is zero. Moreover,
for $i\ne j$,
\[
 (O_ND_N)_{ij}
 =C_{ij,N}C_{jj,N}\otimes A_j^2=0.
\]
Consequently
\[
 I+B_N=(I+O_N)(I+D_N),\qquad D_N^2=0.
\]
All eigenvalues of $D_N$ vanish, hence $\det(I+D_N)=1$.
This proves the first identity. Since $O_N$ has zero diagonal blocks,
$\Tr O_N=0$, which proves the second.
\end{proof}

In particular, consecutive zipper labels must differ in every surviving
cyclic product, including across the last and first factors. For
$k\ge2$ this gives
\begin{equation}\label{eq:cyclic-traces}
 \Tr O_N^k
 =
 \sum_{\substack{
   i_1,\ldots,i_k\in\{1,\ldots,m\}\\
   i_a\ne i_{a+1},\;1\le a<k\\
   i_k\ne i_1
 }}
 \Tr\!\left(C_{i_1i_2,N}\cdots C_{i_ki_1,N}\right)
 \Tr\!\left(A_{i_2}\cdots A_{i_k}A_{i_1}\right).
\end{equation}

\section{Word length and uniform integrability}
\label{sec:word-length}

The trace product in~\eqref{eq:topological-correlator} can grow
exponentially in the total word length of the loops. We bound this
length by the number of traversals of finitely many fixed slabs.
Proposition~\ref{prop:traversal-moments} then gives uniform moments of
all orders, as needed to remove the frozen parts of the zippers and
identify the lamination law.

Fix a finite set $P=\{p_1,\ldots,p_m\}$ of distinct points in a
bounded simply connected planar domain and choose free generators for
the punctured fundamental group. For a finite family $\omega$ of
disjoint simple loops, set
\[
 \ell_P(\omega)=\sum_{\gamma\text{ retained}}|[\gamma]|_{\rm cyc}.
\]
The sum runs over loops enclosing at least two points of $P$, counted
with multiplicity. Here $|[\gamma]|_{\rm cyc}$ is the shortest word length
in the conjugacy class of $\gamma$ in the chosen generators and their
inverses.

For a slab $\mathcal R$ between two lines of constant first or second
coordinate, let $\NN_\omega(\mathcal R)$ count the subarcs joining the
two boundary lines with interior in the open slab. Each traversal is
counted separately, including those belonging to the same loop.

\begin{lemma}\label{lem:word-length}
There are finitely many fixed slabs
$\mathcal R_1,\ldots,\mathcal R_J$ and a constant $C_P$
such that, for every finite collection $\omega$ of
disjoint simple polygonal loops in the punctured domain,
\begin{equation}\label{eq:word-traversal}
 \ell_P(\omega)
 \le C_P\sum_{j=1}^J\NN_\omega(\mathcal R_j).
\end{equation}
The same slabs and constant work when each $p_i$ moves
in a sufficiently small fixed neighborhood of its
initial location, with the generators transported
continuously.
\end{lemma}

\begin{proof}
\emph{Coordinates and slabs.}
Choose rotated coordinates $(x,y)$ in which the horizontal
coordinates of the punctures differ by more than $4\delta$, for some
$\delta>0$. This separation persists on a small interval of rotation
angles. In the original coordinates take a grid of spacing
$h<\delta/(4\sqrt2)$. Take the finitely many slabs of this grid that meet the domain.

\emph{The word of a loop.}
For the given polygonal collection choose an angle in the prescribed
interval such that no vertex lies on a downward ray from a puncture
or on a line $x=x_i\pm\delta$, where $x_i$ is the rotated horizontal
coordinate of $p_i$. This excludes only finitely many angles.
Orient the retained loops and read a letter $a_i$ or $a_i^{-1}$
at each crossing of the downward ray from $p_i$, according as the
loop crosses from left to right or from right to left. These rays
give the usual free generators in the punctured plane. The inclusion
of the punctured simply connected domain into that plane identifies
their fundamental groups.

\emph{One letter for each arc in a band.}
The bands $B_i=\{|x-x_i|<\delta\}$ are disjoint, and a loop
contained in one band encloses at most one puncture. Let $\alpha$ be the closure
of a component of a retained loop inside $B_i$. It is an arc with ends
on $\partial B_i$, and its crossings reduce to
$a_i^{n(\alpha)}$, where $n(\alpha)$ is the signed crossing count
of the one ray in this band. To see that $|n(\alpha)|\le1$, close
$\alpha$ along its boundary line if both endpoints are on the same
line. If they are on opposite lines, continue vertically up the two
boundary lines and join them above both $\alpha$ and $p_i$. In either
case the added path meets the arc only at its endpoints and misses
the ray. The resulting simple closed curve has winding number
$n(\alpha)$ about $p_i$.

\emph{One traversal for each letter.}
Whenever $n(\alpha)\ne0$, a subarc from an endpoint of $\alpha$ to
the ray spans $\delta$ in the rotated $x$ coordinate. It thus spans
at least $\delta/\sqrt2>4h$ in one original coordinate and contains
a traversal between two consecutive grid lines. Distinct contributing
arcs have disjoint interiors, so their traversals are distinct.
The sum of the traversal counts consequently bounds the total cyclic
word length in the ray generators. Writing each $a_i$ in the original
free basis costs at most a fixed factor $C_P$.

\emph{Moving punctures.}
Throughout the allowed interval of angles no ray passes through another
puncture, so the change of generators is fixed. Shrinking the puncture
neighborhoods if necessary preserves the $4\delta$ separation for the
same interval of angles. The same grid and $C_P$ then work in the
transported basis.
\end{proof}

\begin{proposition}[Uniform moments of topological observables]
\label{prop:topological-ui}
For every compact set $\mathcal C\subset X_{\rm unip}$
and every $1\le p<\infty$,
\begin{equation}\label{eq:topological-ui}
 \sup_N\EE\!\left[
   \sup_{\rho\in\mathcal C}
   \left|\prod_{\gamma\in\mathcal L_N}
                 \frac{\Tr\rho(\gamma)}2\right|^p\right]<\infty.
\end{equation}
This remains true for finitely many moving punctures in the fixed
separated neighborhoods of Lemma~\ref{lem:word-length}, provided their
representations remain in a compact set in the associated fixed bases.
\end{proposition}

\begin{proof}
Let $P_N$ be the discrete puncture set. Compactness gives a single
bound $R\ge1$ for the operator norms of the images of the generators
and their inverses. Loops enclosing at most one puncture have
normalized trace $1$. For a retained loop, evaluate the trace on a
shortest representative of its conjugacy class to obtain
\[
 \left|\frac{\Tr\rho(\gamma)}2\right|
 \le R^{|[\gamma]|_{\rm cyc}}.
\]
The product is therefore at most $R^{\ell_{P_N}(\mathcal L_N)}$.
Apply~\eqref{eq:word-traversal} in the $(s,y)$ coordinates and
discard grid slabs whose traversal counts are eventually zero.
For each remaining slab $\mathcal R_j$, choose a fixed strictly
smaller slab $\mathcal R'_j$ with boundary coordinates in $(0,1)$.
For all sufficiently large $N$, each traversal of $\mathcal R_j$
contains a distinct traversal between the discrete boundary sections
of $\mathcal R'_j$. Hence
\[
 \NN_{\mathcal L_N}(\mathcal R_j)\le\NN_N(\mathcal R'_j).
\]
Apply~\eqref{eq:traversal-exponential-moment} with
$A=p C_P\log R$. For the finitely many smaller $N$ the expectation is
finite, because there are finitely many configurations. The slabs
and $C_P$ are uniform on the prescribed puncture neighborhoods.
\end{proof}

\section{Interior limits and the arctic boundary reduction}
\label{sec:interior-boundary}
\label{sec:arctic-tail}

This section computes the limit of the correlators $F_N(\rho)$ and
proves Theorem~\ref{thm:aztec-finite-puncture}. We truncate the
zippers at distance $\varepsilon$ inside the arctic boundary, compute
the limit on the truncated cuts, and then remove the truncation. The
steps are as follows.
\begin{enumerate}[label=\textnormal{(\arabic*)},leftmargin=*]
\item Section~\ref{subsec:scalar-inverse} gives an exact formula for
the inverse Kasteleyn matrix, as a double contour integral
(Lemma~\ref{lem:scalar-exact-inverse}).
\item Section~\ref{subsec:interior-traces} contains the saddle-point
asymptotics in the liquid region. The four saddle pairs give the limit
$D_\varepsilon$ of the determinant on the truncated cuts
(Proposition~\ref{prop:interior-determinant}).
\item Sections~\ref{subsec:arctic-edge}--\ref{subsec:boundary-reduction}
remove the truncation. Restoring the arctic collars costs
$O(\varepsilon^{1/4})$ in the correlator, and the probability of a
nondoubled edge on the remaining frozen suffixes tends to zero
(Theorem~\ref{thm:boundary-reduction}).
\item In Section~\ref{subsec:canonical-determinant}, as
$\varepsilon\downarrow0$ the determinants $D_\varepsilon$ converge to
a determinant $D_0$, which is identified with the correlator of nested
$\CLE_4$ by comparison with the half-plane model
of~\cite{Dubedat,BI}, followed by an analytic continuation in the
monodromy
(Lemma~\ref{lem:canonical-reference-identification}).
\item Section~\ref{subsec:frozen-regions} shows that the perfect matchings are
frozen outside the arctic ellipse
(Proposition~\ref{prop:frozen-regions}).
\end{enumerate}

\subsection{An exact scalar inverse}
\label{subsec:scalar-inverse}

In the coordinates of~\eqref{eq:diagonal-coordinates}, write
\[
 \begin{aligned}
 v^{\rm W}_{m,n}&=(2m+1,2n+2),&&0\le m<N,\ -1\le n<N,\\
 v^{\rm B}_{m,n}&=(2m,2n+1),&&0\le m\le N,\ 0\le n<N.
 \end{aligned}
\]
These are the white and black vertices, respectively. To a vertex with indices $(m,n)$ we assign the macroscopic
coordinates
\begin{equation}\label{eq:vertex-macroscopic-coordinates}
 \boldsymbol x_{v^{\rm W}_{m,n}}
 =\boldsymbol x_{v^{\rm B}_{m,n}}
 =\left(\frac{2m}{N}-1,\frac{2n}{N}-1\right).
\end{equation}
The nonzero entries of~\eqref{eq:kasteleyn-definition} become
\begin{equation}\label{eq:scalar-kasteleyn-entries}
\begin{aligned}
 K_N(v^{\rm W}_{m,n},v^{\rm B}_{m,n})&=1,&
 K_N(v^{\rm W}_{m,n},v^{\rm B}_{m,n+1})&=1,\\
 K_N(v^{\rm W}_{m,n},v^{\rm B}_{m+1,n})&=-\bk,&
 K_N(v^{\rm W}_{m,n},v^{\rm B}_{m+1,n+1})&=1,
\end{aligned}
\end{equation}
whenever the indicated vertices exist. In the notation of
Sections~\ref{subsec:seams} and~\ref{subsec:layer}, one has
$b_{r,i}=v^{\rm B}_{r,N-1-i}$ and $w_{r,x}=v^{\rm W}_{r-1,N-1-x}$.

\begin{lemma}[Exact inverse Kasteleyn matrix]
\label{lem:scalar-exact-inverse}
Let $m,n,M,\nu$ be integers satisfying
\[
 0\le m\le N,\qquad 0\le n<N,\qquad
 0\le M<N,\qquad -1\le\nu<N.
\]
In the contour variables $z,t$, let $\Gamma_0=\{|z|=r_0\}$ and $\Gamma_1=\{|t-\bk^{-1}|=r_1\}$
be positively oriented, where $r_0,r_1>0$ and $r_0+r_1<\bk^{-1}$.
Then
\begin{equation}\label{eq:scalar-exact-inverse}
\begin{aligned}
 K_N^{-1}(v^{\rm B}_{m,n},v^{\rm W}_{M,\nu})
 &=\frac1{(2\pi\ii)^2}
 \int_{\Gamma_0}\int_{\Gamma_1}
 \frac{(z+1)^m(\bk z-1)^{N-m}}{z^{n+1}}
 \frac{t^{\nu+1}}{(t+1)^{M+1}(\bk t-1)^{N-M}}
 \frac{dt\,dz}{z-t}\\
 &\quad-\frac{\one_{\{m>M\}}}{2\pi\ii}
 \int_{\Gamma_0}
 \frac{z^{\nu-n}(z+1)^{m-M-1}}{(\bk z-1)^{m-M}}\,dz.
\end{aligned}
\end{equation}
\end{lemma}

\begin{proof}
The biased inverse in \cite[Theorem~2.3]{CJY} is related to
\eqref{eq:scalar-exact-inverse} by exchanging axes, applying the gauge
in~\eqref{eq:CJY-phase-conversion}, rescaling the contours and rewriting
the single integral by residues. We verify the formula in our coordinates.

The matrix $K_N$ is square and $|\det K_N|=Z_N\ne0$, so it suffices
to show that the right-hand side $G(b,w)$
of~\eqref{eq:scalar-exact-inverse} satisfies
\begin{equation}\label{eq:inverse-check}
 \sum_bK_N(w',b)\,G(b,w)=\one_{\{w'=w\}}
\end{equation}
for all white vertices $w'=v^{\rm W}_{\mu,\lambda}$ and
$w=v^{\rm W}_{M,\nu}$. Write $G=\mathcal D-\mathcal E$ for the double
and the single integral, and put
\[
 \begin{gathered}
 f_{m,n}(z)=\frac{(z+1)^m(\bk z-1)^{N-m}}{z^{n+1}},
 \qquad
 g(t)=\frac{t^{\nu+1}}{(t+1)^{M+1}(\bk t-1)^{N-M}},\\
 C(z)=(z+1)^\mu(\bk z-1)^{N-1-\mu}.
 \end{gathered}
\]
Then $\mathcal D(v^{\rm B}_{m,n},w)$ is the double integral of
$f_{m,n}(z)g(t)/(z-t)$, and
\[
 \mathcal E(v^{\rm B}_{m,n},w)
 =\frac{\one_{\{m>M\}}}{2\pi\ii}\int_{\Gamma_0}f_{m,n}(z)g(z)\,dz.
\]
By~\eqref{eq:scalar-kasteleyn-entries}, the left-hand side
of~\eqref{eq:inverse-check} involves the combination
\begin{equation}\label{eq:four-term-identity}
 f_{\mu,\lambda}+f_{\mu,\lambda+1}
 -\bk f_{\mu+1,\lambda}+f_{\mu+1,\lambda+1}
 =\frac{C(z)}{z^{\lambda+2}}
 \bigl[(\bk z-1)z+(\bk z-1)-\bk(z+1)z+(z+1)\bigr]=0
\end{equation}
when $0\le\lambda\le N-2$, so that all four black neighbors exist.
For $\lambda=-1$ only $v^{\rm B}_{\mu,0}$ and $v^{\rm B}_{\mu+1,0}$
exist, and
\begin{equation}\label{eq:two-term-low}
 f_{\mu,0}+f_{\mu+1,0}=(1+\bk)C(z).
\end{equation}
For $\lambda=N-1$ only $v^{\rm B}_{\mu,N-1}$ and
$v^{\rm B}_{\mu+1,N-1}$ exist, and
\begin{equation}\label{eq:two-term-high}
 f_{\mu,N-1}-\bk f_{\mu+1,N-1}=-(1+\bk)\frac{C(z)}{z^N}.
\end{equation}

\emph{The double integral.}
By~\eqref{eq:four-term-identity}, the contribution of $\mathcal D$
to~\eqref{eq:inverse-check} vanishes for $0\le\lambda\le N-2$.
For $\lambda=-1$, \eqref{eq:two-term-low} shows that the integrand
is holomorphic in $z$ inside $\Gamma_0$, because $t\in\Gamma_1$ lies
outside $\Gamma_0$; the contribution vanishes again.
For $\lambda=N-1$, the function $C(z)z^{-N}(z-t)^{-1}$ is
$O(|z|^{-2})$ at infinity and its only poles are $0$ and $t$, so
\[
 \frac1{2\pi\ii}\int_{\Gamma_0}\frac{C(z)}{z^N}\frac{dz}{z-t}
 =-\frac{C(t)}{t^N}.
\]
By~\eqref{eq:two-term-high}, the contribution of $\mathcal D$ is
therefore
\begin{equation}\label{eq:double-boundary-term}
 \frac{1+\bk}{2\pi\ii}\int_{\Gamma_1}h(t)\,dt,
 \qquad
 h(t)=\frac{C(t)g(t)}{t^N}
 =t^{\nu+1-N}(t+1)^{\mu-M-1}(\bk t-1)^{M-\mu-1}.
\end{equation}
It vanishes for $\mu<M$, since $h$ is then holomorphic inside
$\Gamma_1$.

\emph{The case $\mu<M$.}
By~\eqref{eq:scalar-kasteleyn-entries}, every black neighbor of $w'$
has first index $m\in\{\mu,\mu+1\}$, so $m\le M$ and the indicator
in $\mathcal E$ vanishes. The double-integral contribution also
vanishes, as shown above, proving~\eqref{eq:inverse-check} in this case.

\emph{The case $\mu>M$.}
All neighbors of $w'$ carry a single integral with the common factor
$g(z)$, and~\eqref{eq:four-term-identity}--\eqref{eq:two-term-high}
apply to them. For $0\le\lambda\le N-2$ the contribution of
$\mathcal E$ vanishes. For $\lambda=-1$ the integrand
$(1+\bk)C(z)g(z)$ is holomorphic inside $\Gamma_0$, because
$\nu+1\ge0$. For $\lambda=N-1$ the contribution of $\mathcal E$ is
$-(1+\bk)(2\pi\ii)^{-1}\int_{\Gamma_0}h(z)\,dz$. Together
with~\eqref{eq:double-boundary-term}, the left-hand side
of~\eqref{eq:inverse-check} becomes
\[
 \frac{1+\bk}{2\pi\ii}\left(\int_{\Gamma_0}+\int_{\Gamma_1}\right)h(t)\,dt.
\]
Since $\mu>M$, all finite poles of $h$ lie in $\{0,\bk^{-1}\}$, and
$h(t)=O(|t|^{\nu-N-1})=O(|t|^{-2})$ at infinity. The sum of the two
integrals is therefore zero.

\emph{The case $\mu=M$.}
A neighbor carries a single integral precisely when its first index is $M+1$.
Expanding
$(\bk z-1)^{-1}=-\sum_{j\ge0}\bk^jz^j$ gives
\[
 \mathcal E(v^{\rm B}_{M+1,n},w)
 =\frac1{2\pi\ii}\int_{\Gamma_0}\frac{z^{\nu-n}}{\bk z-1}\,dz
 =-\bk^{\,n-\nu-1}\one_{\{n>\nu\}}.
\]
For $0\le\lambda\le N-2$ the contribution of $\mathcal D$ vanishes,
and that of $-\mathcal E$ is
\[
 -\bigl[-\bk\,\mathcal E(v^{\rm B}_{M+1,\lambda},w)
        +\mathcal E(v^{\rm B}_{M+1,\lambda+1},w)\bigr]
 =\bk^{\lambda-\nu}
   \bigl[\one_{\{\lambda\ge\nu\}}-\one_{\{\lambda>\nu\}}\bigr]
 =\one_{\{\lambda=\nu\}}.
\]
For $\lambda=-1$ the contribution of $\mathcal D$ vanishes and that of
$-\mathcal E$ is $-\mathcal E(v^{\rm B}_{M+1,0},w)
=\bk^{-\nu-1}\one_{\{\nu<0\}}=\one_{\{\nu=-1\}}$.
For $\lambda=N-1$ the contribution of $-\mathcal E$ is
$\bk\,\mathcal E(v^{\rm B}_{M+1,N-1},w)
=-\bk^{N-1-\nu}\one_{\{\nu<N-1\}}$, while
\eqref{eq:double-boundary-term} equals
\[
 (1+\bk)\operatorname{Res}_{t=1/\bk}
 \frac{t^{\nu+1-N}}{(t+1)(\bk t-1)}
 =\bk^{N-1-\nu}.
\]
Their sum is $\bk^{N-1-\nu}\one_{\{\nu=N-1\}}=\one_{\{\lambda=\nu\}}$.
In all cases the left-hand side of~\eqref{eq:inverse-check} equals
$\one_{\{\mu=M,\ \lambda=\nu\}}$.
\end{proof}

\subsection{Interior cyclic traces and determinants}
\label{subsec:interior-traces}

For a point $z$ in the rescaled physical plane, write
\begin{equation}\label{eq:physical-diagonal-coordinates}
 \boldsymbol x=(\xi,\eta)
 =(\Re z+\Im z,\Re z-\Im z).
\end{equation}
At a lattice vertex these coordinates are $(u/N-1,v/N-1)$
by~\eqref{eq:diagonal-coordinates}. Subtracting the indexed coordinates
in~\eqref{eq:vertex-macroscopic-coordinates} from this pair gives
$(1,2)/N$ at white vertices and $(0,1)/N$ at black vertices.
Join each $z_i$ to a distinct point of $\partial\Diamond$
by pairwise disjoint paths of finitely many horizontal and vertical
segments in the $(\xi,\eta)$-plane.
Each path crosses the arctic boundary once, transversely on a vertical
segment, away from its contacts with $\partial\Diamond$.

In the scaled embedding, approximate the paths by disjoint
lattice zippers starting in the faces containing $z_{i,N}^0$, as specified before~\eqref{eq:unipotent}
in Section~\ref{sec:punctures}. In the lattice coordinates $(u,v)$,
an $f$-type face has center $f_{m,n}=(2m,2n+2)$ and an $f^*$-type
face has center $f^*_{m,n}=(2m+1,2n+3)$. The two-step periods for the
positive horizontal and vertical directions are, respectively,
\begin{equation}\label{eq:zipper-periods}
 f_{m,n}\to f^*_{m,n}\to f_{m+1,n},\qquad
 f_{m,n}\to f^*_{m,n}\to f_{m,n+1}.
\end{equation}
From an initial $f^*$-type face, first take one dual edge to an
adjacent $f$-type face. Join successive segments at $f$-type face
centers, reverse the periods for negative directions, and remove any
immediate backtracking at the joins.

For sufficiently small $\varepsilon>0$, stop each path at vertical
distance $\varepsilon$ inside the arctic boundary and its lattice
zipper at the nearest face center along the zipper.

Let $E_{N,\varepsilon}\subset E$ be the edges crossed
by the corresponding truncated lattice zippers.
Define the diagonal projection
\[
 (P_{N,\varepsilon})_{ef}
 =\one_{\{e=f\}}\one_{\{e\in E_{N,\varepsilon}\}}I_2,
 \qquad e,f\in E,
\]
and, with $O_N(\rho)$ as in Proposition~\ref{prop:nilpotent-reduction}, set
\begin{equation}\label{eq:cutoff-operator}
 O_{N,\varepsilon}=P_{N,\varepsilon}O_NP_{N,\varepsilon},
 \qquad
 F_{N,\varepsilon}(\rho)
 =\det{}_2(I+O_{N,\varepsilon}(\rho)).
\end{equation}
For fixed $\varepsilon$, the truncated macroscopic paths
stay a positive distance from one another and from the
arctic boundary. Each truncated lattice zipper crosses
$O(N)$ edges.

Write $L_{\xi,\eta}$ for the image of $L$ in the $(\xi,\eta)$-plane. Set
\begin{equation}\label{eq:xi-hat}
 \hat\xi(\xi,\eta)
 =\frac{(1+\bk)\xi+(\bk-1)\eta}{2\sqrt\bk},
 \qquad
 L_{\xi,\eta}=\bigl\{(\xi,\eta)\in\mathbb R^2:
                  \hat\xi(\xi,\eta)^2+\eta^2<1\bigr\}.
\end{equation}
For $\bk=1$ we have $\hat\xi=\xi$, and $L_{\xi,\eta}$ is the unit disk.
In general
\begin{equation}\label{eq:liquid-symmetric-form}
 4\bk\bigl(1-\hat\xi(\xi,\eta)^2-\eta^2\bigr)
 =4\bk-(1+\bk)^2(\xi^2+\eta^2)-2(\bk^2-1)\xi\eta,
\end{equation}
so $L_{\xi,\eta}$ is an open ellipse which is symmetric under
$(\xi,\eta)\mapsto(\eta,\xi)$ and inscribed in the square
$[-1,1]^2$. In terms of $z$,
\[
 \hat\xi=\sqrt\bk\,\Re z+\frac{\Im z}{\sqrt\bk},
 \qquad
 \hat\xi^2+\eta^2
 =(1+\bk)\Bigl((\Re z)^2+\frac{(\Im z)^2}{\bk}\Bigr),
\]
so $z\in L$ if and only if $\boldsymbol x\in L_{\xi,\eta}$,
by~\eqref{eq:physical-liquid}. We write $\hat\xi(\boldsymbol x)$ for
$\hat\xi(\xi,\eta)$.
Define $q:L\to\mathbb H$, where
$\mathbb H=\{w\in\mathbb C:\Im w>0\}$, by
\begin{equation}\label{eq:spectral-coordinate}
 q(z)=\frac{\hat\xi+\ii\sqrt{1-\hat\xi^2-\eta^2}}{\sqrt\bk\,(1-\eta)}
 =\frac{\ii}{\sqrt\bk}\,
   \frac{1+e^{-\ii\theta_\bk}\overline{\Phi(z)}}
        {1-e^{-\ii\theta_\bk}\overline{\Phi(z)}},
 \qquad
 \theta_\bk=\arctan\sqrt\bk\in(0,\pi/2).
\end{equation}
For $\bk=1$ the angle is $\theta_1=\pi/4$.
To verify the second equality, write $w=\Phi(z)$ and
$\omega=e^{-\ii\theta_\bk}\overline w$. Since
$\cos\theta_\bk=(1+\bk)^{-1/2}$ and
$\sin\theta_\bk=\sqrt\bk\,(1+\bk)^{-1/2}$, the formula
$z=\Psi(w)$ in~\eqref{eq:intrinsic-map} gives
\[
 \eta=\frac{2\Re\omega}{1+|\omega|^2},\qquad
 \hat\xi=-\frac{2\Im\omega}{1+|\omega|^2},\qquad
 1-\eta=\frac{|1-\omega|^2}{1+|\omega|^2},\qquad
 \sqrt{1-\hat\xi^2-\eta^2}=\frac{1-|\omega|^2}{1+|\omega|^2},
\]
and therefore
\[
 \frac{\hat\xi+\ii\sqrt{1-\hat\xi^2-\eta^2}}{1-\eta}
 =\frac{-2\Im\omega+\ii(1-|\omega|^2)}{|1-\omega|^2}
 =\ii\,\frac{1+\omega}{1-\omega}.
\]
The map $w\mapsto\ii\bk^{-1/2}(1+e^{-\ii\theta_\bk}\overline w)
/(1-e^{-\ii\theta_\bk}\overline w)$ is an anticonformal bijection
from $\DD$ onto $\mathbb H$.
Thus $q$ is a bijection from $L$ onto $\mathbb H$ and
reverses orientation. Let $\gamma_i^\varepsilon$ be the image under $q$ of
the truncated $i$th zipper, oriented from the puncture
toward the real axis. Define the oriented chain
\begin{equation}\label{eq:cutoff-contours}
 C_i^\varepsilon
 =\gamma_i^\varepsilon-\overline{\gamma_i^\varepsilon},
\end{equation}
where the bar denotes reflection across the real axis
with the inherited parametrization, and the minus sign
reverses its orientation.
On $\mathcal H_\varepsilon=\bigoplus_iL^2(C_i^\varepsilon,|dq|;
\mathbb C^2)$ define the off-diagonal Cauchy operator by
\begin{equation}\label{eq:canonical-cutoff-kernel}
 (\mathcal O_\varepsilon f)_i(q)
 =\sum_{j\ne i}\int_{C_j^\varepsilon}
       \frac{A_j f_j(r)}{2\pi\ii(q-r)}\,dr,
 \qquad D_\varepsilon(\rho)=\det{}_2(I+\mathcal O_\varepsilon(\rho)).
\end{equation}
Transport the representation by $q$, so that
$\rho_q(\beta)=\rho(q^{-1}\beta)$.
Let $g_i$ be the counterclockwise generator associated
with zipper $i$, with monodromy $M_i$.
Since $q$ reverses orientation, the counterclockwise
generator associated with the image zipper is the
reversed path $(q\circ g_i)^{-1}$.
Its monodromy under the transported connection is
$M_i^{-1}=I_2-A_i$.

Distinct contours are separated and have finite length, so
$\mathcal O_\varepsilon$ is Hilbert--Schmidt.

For $\boldsymbol x=(\xi,\eta)$, define the phase function
in the complex variable $\zeta$ by
\begin{equation}\label{eq:uniform-action}
 F(\zeta;\boldsymbol x)
 =\frac{1+\eta}{2}\log\zeta
  -\frac{1+\xi}{2}\log(\zeta+1)
  -\frac{1-\xi}{2}\log(\bk\zeta-1).
\end{equation}
Its critical points are the roots of
$\bk(1-\eta)\zeta^2-2\sqrt\bk\,\hat\xi(\boldsymbol x)\zeta
+(1+\eta)=0$. For $\boldsymbol x\in L_{\xi,\eta}$ they are
\[
 \zeta_\pm(\boldsymbol x)
 =\frac{\hat\xi\pm\ii\sqrt{1-\hat\xi^2-\eta^2}}
        {\sqrt\bk\,(1-\eta)},
 \qquad \hat\xi=\hat\xi(\boldsymbol x).
\]
For $z\in L$ and $\boldsymbol x$ as
in~\eqref{eq:physical-diagonal-coordinates},
we have $q(z)=\zeta_+(\boldsymbol x)$.

Fix $\varepsilon>0$ as above. All constants below may depend on
$\varepsilon$ and on the fixed paths.
We also write $w(\zeta)=(\zeta+1)/(\bk\zeta-1)$.
For lattice vertices we use~\eqref{eq:vertex-macroscopic-coordinates}.
A full horizontal or vertical period changes the corresponding
coordinate by $2/N$.
Choose a compact set $\mathcal K$ strictly inside
$L_{\xi,\eta}$ containing a neighborhood of all truncated paths.
The estimates for retained lattice vertices below are for
$N\ge N_0$, with $N_0$ sufficiently large.

Differentiating~\eqref{eq:uniform-action} gives
\begin{equation}\label{eq:bulk-saddle-equation}
 F_\zeta(\zeta;\boldsymbol x)
 =-\frac{\bk(1-\eta)(\zeta-\zeta_+(\boldsymbol x))
                     (\zeta-\zeta_-(\boldsymbol x))}
          {2\zeta(\zeta+1)(\bk\zeta-1)}.
\end{equation}
Thus the two saddles are distinct, and
$|F_{\zeta\zeta}(\zeta_\pm(\boldsymbol x);\boldsymbol x)|$
is bounded above and below by positive constants on $\mathcal K$.
Use the principal logarithms on the upper half-plane and their
conjugate branches on the lower half-plane, and set
\[
 \mathcal F_\pm(\boldsymbol x)
 =F(\zeta_\pm(\boldsymbol x);\boldsymbol x),\qquad
 \theta_{\rm osc}(\boldsymbol x)
 =2\Im\mathcal F_+(\boldsymbol x).
\]
These functions have bounded derivatives of every fixed order on
$\mathcal K$, and $\mathcal F_- =\overline{\mathcal F_+}$.

\begin{lemma}\label{lem:bulk-four-saddle}
Let $\mathcal K$ be a compact subset of the set $L_{\xi,\eta}$
of~\eqref{eq:xi-hat}
and fix $\delta>0$. Let $b=v^{\rm B}_{m,n}$ and $w=v^{\rm W}_{M,\nu}$ be
black and white vertices, respectively, with macroscopic
coordinates $\boldsymbol x_b,\boldsymbol x_w$ as
in~\eqref{eq:vertex-macroscopic-coordinates}.
Assume
\[
 \boldsymbol x_b,\boldsymbol x_w\in\mathcal K,
 \qquad |\boldsymbol x_b-\boldsymbol x_w|\ge\delta.
\]
Use the phase function $F$ defined
in~\eqref{eq:uniform-action}.
Define
\begin{equation}\label{eq:real-saddle-value}
 \mathcal A(\boldsymbol x)
 =\Re\mathcal F_+(\boldsymbol x)
 =\Re\mathcal F_-(\boldsymbol x).
\end{equation}
Then the inverse of the Kasteleyn matrix
defined in~\eqref{eq:kasteleyn-definition} satisfies
\begin{equation}\label{eq:four-saddle-entry}
\begin{aligned}
 K_N^{-1}(b,w)
 &=\frac1N\sum_{\alpha,\beta\in\{+,-\}}
 c_{\alpha\beta}(\boldsymbol x_b,\boldsymbol x_w)
 e^{N[\mathcal F_\beta(\boldsymbol x_w)
       -\mathcal F_\alpha(\boldsymbol x_b)]}
\\&\quad+e^{N[\mathcal A(\boldsymbol x_w)-\mathcal A(\boldsymbol x_b)]}
 O(N^{-3/2}).
\end{aligned}
\end{equation}
Set
\[
 d_b(\zeta)=\ii\zeta,\qquad
 d_w(\zeta)=\frac{\ii}{1+\bk}(\zeta+1)(\bk\zeta-1).
\]
For $\alpha,\beta\in\{+,-\}$, define
\[
 \begin{aligned}
 H_{b,\alpha}
 &=F_{\zeta\zeta}\bigl(\zeta_\alpha(\boldsymbol x_b);
                       \boldsymbol x_b\bigr)
   d_b\bigl(\zeta_\alpha(\boldsymbol x_b)\bigr)^2,\\
 H_{w,\beta}
 &=-F_{\zeta\zeta}\bigl(\zeta_\beta(\boldsymbol x_w);
                        \boldsymbol x_w\bigr)
   d_w\bigl(\zeta_\beta(\boldsymbol x_w)\bigr)^2.
 \end{aligned}
\]
Both $H_{b,\alpha}$ and $H_{w,\beta}$ have positive real part.
With the principal square roots, set
\[
 \mathfrak c_{b,\alpha}
 =\frac{d_b\bigl(\zeta_\alpha(\boldsymbol x_b)\bigr)}
        {\sqrt{H_{b,\alpha}}},\qquad
 \mathfrak c_{w,\beta}
 =\frac{d_w\bigl(\zeta_\beta(\boldsymbol x_w)\bigr)}
        {\sqrt{H_{w,\beta}}}.
\]
The coefficients in~\eqref{eq:four-saddle-entry} are
\begin{equation}\label{eq:uniform-saddle-coefficient}
 c_{\alpha\beta}(\boldsymbol x_b,\boldsymbol x_w)
 =-\frac{\zeta_\beta(\boldsymbol x_w)
          \mathfrak c_{b,\alpha}\mathfrak c_{w,\beta}}
 {2\pi\zeta_\alpha(\boldsymbol x_b)
  \bigl(\zeta_\beta(\boldsymbol x_w)+1\bigr)
  \bigl(\zeta_\alpha(\boldsymbol x_b)
        -\zeta_\beta(\boldsymbol x_w)\bigr)}.
\end{equation}
These coefficients have uniformly bounded spatial derivatives of every
fixed order. In particular,
\begin{equation}\label{eq:uniform-saddle-coefficient-bound}
 |c_{\alpha\beta}|+|\nabla_{\boldsymbol x_b}c_{\alpha\beta}|
 +|\nabla_{\boldsymbol x_w}c_{\alpha\beta}|\le C.
\end{equation}
\end{lemma}

\begin{proof}
We deform the two contours of the exact
formula~\eqref{eq:scalar-exact-inverse} to contours on which the real
part of the phase is extremal at the saddles, bound the residue term
and the parts of the integral away from the saddles, and evaluate the
contributions of the four pairs of saddles.

\emph{Saddles and contours.}
Write $\boldsymbol x_b=(\xi_b,\eta_b)$ and
$\boldsymbol x_w=(\xi_w,\eta_w)$, and put
$\hat\xi_b=\hat\xi(\boldsymbol x_b)$,
$\hat\xi_w=\hat\xi(\boldsymbol x_w)$. Direct substitution gives
\[
 \Re H_{b,\alpha}
 =\frac{2\bk\,(1-\hat\xi_b^2-\eta_b^2)}
       {(1+\bk)^2(1-\xi_b^2)}>0,
 \qquad
 \Re H_{w,\beta}
 =\frac{2\bk\,(1-\hat\xi_w^2-\eta_w^2)}
       {(1+\bk)^2(1-\eta_w^2)}>0.
\]

Put $\mathsf w(z)=(z+1)/(\bk z-1)$. This M\"obius transformation
is an involution of the Riemann sphere; it exchanges $-1$ with $0$
and $\bk^{-1}$ with $\infty$. A direct computation gives the duality
\begin{equation}\label{eq:action-duality}
 F\bigl(t;(\xi,\eta)\bigr)
 =-F\bigl(\mathsf w(t);(\eta,\xi)\bigr)-\log(1+\bk),
\end{equation}
up to the choice of branches, which does not affect real parts.
Define
\begin{equation}\label{eq:contour-radii}
 r_b=\sqrt{\frac{1+\eta_b}{\bk(1-\eta_b)}},\qquad
 r_w=\sqrt{\frac{1+\xi_w}{\bk(1-\xi_w)}}.
\end{equation}
Let $C_b$ be the positively oriented circle $|z|=r_b$.
Define
\begin{equation}\label{eq:white-contour}
 W=\{t:|\mathsf w(t)|>r_w\},\qquad C_w=\partial W.
\end{equation}
Orient $C_w$ with $W$ on its left.
This region $W$ contains $\bk^{-1}$ and excludes $-1$; it is
considered on the Riemann sphere when its boundary passes through
infinity.

\emph{Real-part bounds.}
On $C_b$, writing $u=\cos\theta$ gives
\[
 \Re F(r_be^{\ii\theta};\boldsymbol x_b)
 =\frac{1+\eta_b}{2}\log r_b
 -\frac{1+\xi_b}{4}\log(r_b^2+1+2r_bu)
 -\frac{1-\xi_b}{4}\log(\bk^2r_b^2+1-2\bk r_bu).
\]
This is strictly convex in $u$. Its derivative in $u$ vanishes when
$(1+\xi_b)(\bk^2r_b^2+1-2\bk r_bu)
=\bk(1-\xi_b)(r_b^2+1+2r_bu)$; since
$\bk r_b^2+1=2/(1-\eta_b)$ and $\bk r_b^2-1=2\eta_b/(1-\eta_b)$,
this is the point $u=\hat\xi_b/\sqrt{1-\eta_b^2}$, which lies in
$(-1,1)$ because $\boldsymbol x_b\in L_{\xi,\eta}$.
The two minimizing points on $C_b$ are $\zeta_\pm(\boldsymbol x_b)$:
indeed $|\zeta_\pm(\boldsymbol x_b)|=r_b$ and
$\Re\zeta_\pm(\boldsymbol x_b)/r_b=\hat\xi_b/\sqrt{1-\eta_b^2}$.
For $v=\mathsf w(t)$, the duality~\eqref{eq:action-duality} reads
\[
 \Re F(t;\boldsymbol x_w)
 =-\log(1+\bk)+\frac{1+\eta_w}{2}\log|v+1|
 +\frac{1-\eta_w}{2}\log|\bk v-1|
 -\frac{1+\xi_w}{2}\log|v|.
\]
On $|v|=r_w$ this is strictly concave in $\cos\arg v$,
with its two maxima at the images of $\zeta_\pm(\boldsymbol x_w)$.
Indeed, by~\eqref{eq:action-duality} the right-hand side is
$-\log(1+\bk)-\Re F(v;(\eta_w,\xi_w))$; the point $(\eta_w,\xi_w)$
belongs to $L_{\xi,\eta}$ by the symmetry of $L_{\xi,\eta}$, the circle
$|v|=r_w$ is its black circle~\eqref{eq:contour-radii}, and
$\mathsf w$ maps the critical points of $F(\,\cdot\,;\boldsymbol x_w)$
to those of $F(\,\cdot\,;(\eta_w,\xi_w))$. The claim therefore follows
from the convexity statement just proved.
Consequently, with $\mathcal A$ defined
in~\eqref{eq:real-saddle-value},
\[
 \Re F(z;\boldsymbol x_b)\ge\mathcal A(\boldsymbol x_b)
 \quad(z\in C_b),
 \qquad
 \Re F(t;\boldsymbol x_w)\le\mathcal A(\boldsymbol x_w)
 \quad(t\in C_w).
\]
Both inequalities have uniform quadratic gaps near their equality
points and uniform positive gaps away from fixed neighborhoods of them.

\emph{Contour deformation.}
We record the residue left by the contour deformation. Put
$d=m-M$ and
\[
 R(z)=\frac{e^{N[F(z;\boldsymbol x_w)-F(z;\boldsymbol x_b)]}}{z+1}
 =\frac{z^{\nu-n}(z+1)^{d-1}}{(\bk z-1)^d},
\]
where the last identity follows from~\eqref{eq:uniform-action}
and~\eqref{eq:vertex-macroscopic-coordinates}.
The exact formula~\eqref{eq:scalar-exact-inverse} becomes
\begin{equation}\label{eq:uniform-deformed-inverse}
 K_N^{-1}(b,w)
 =\frac1{(2\pi\ii)^2}\int_{C_b}\int_{C_w}
 e^{N[F(t;\boldsymbol x_w)-F(z;\boldsymbol x_b)]}
 \frac{t\,dt\,dz}{(t+1)z(z-t)}+\mathcal R_N,
\end{equation}
where
\[
 \mathcal R_N=
 \begin{cases}
 -\dfrac1{2\pi\ii}\displaystyle\int_{C_b\setminus W}R(z)\,dz,
     &d>0,\\[2mm]
 \dfrac1{2\pi\ii}\displaystyle\int_{C_b\cap W}R(z)\,dz,
     &d\le0.
 \end{cases}
\]
All subarcs inherit the orientation of $C_b$.

We now prove~\eqref{eq:uniform-deformed-inverse}.
We first assume that $r_b\ne\bk^{-1}$ and that
$C_b$ and $C_w$ are not tangent to each other. Let $\mathcal D(C,C')$ denote the double-integral term
in~\eqref{eq:uniform-deformed-inverse}, with $z$ integrated
along $C$ and $t$ along $C'$. Take $\Gamma_0$ inside $C_b$ and choose $\Gamma_1$ so small
that its closed disk lies in $W$ and is disjoint from $C_b$.
The pole $z=t$ has residues $R(t)$ in $z$ and $-R(z)$ in $t$.
The residue theorem gives
\[
 \mathcal D(\Gamma_0,\Gamma_1)
 =\mathcal D(C_b,C_w)
  +\frac1{2\pi\ii}\int_{C_b\cap W}R(z)\,dz
  -\one_{\{r_b>1/\bk\}}\operatorname{Res}_{z=1/\bk}R(z).
\]
Indeed, moving the $z$-contour from $\Gamma_0$ to $C_b$ crosses the
pole $z=t$, $t\in\Gamma_1$, exactly when $\Gamma_1$ lies inside $C_b$,
and $(2\pi\ii)^{-1}\int_{\Gamma_1}R(t)\,dt
=\operatorname{Res}_{z=1/\bk}R(z)$; moving the $t$-contour from
$\Gamma_1$ to $C_w$ inside $W$ crosses the pole $t=z$ exactly for
$z\in C_b\cap W$.
No additional residues arise at $t=0$ or at infinity: the
$t$-integrand is regular at $0$ because $\nu+1\ge0$, and it is
$O(|t|^{\nu-N-1})=O(|t|^{-2})$ at infinity.
For $d>0$, $R$ is regular at $-1$, and
\[
 \frac1{2\pi\ii}\int_{\Gamma_0}R(z)\,dz
 +\one_{\{r_b>1/\bk\}}\operatorname{Res}_{z=1/\bk}R(z)
 =\frac1{2\pi\ii}\int_{C_b}R(z)\,dz.
\]
For $d\le0$, $R$ is regular at $\bk^{-1}$, and the original single
integral is absent. Substitution
in~\eqref{eq:scalar-exact-inverse} gives the stated two cases.
The uniform integrability estimate below extends the identity
to $r_b=\bk^{-1}$ and contour tangencies.

\emph{The residue.}
An intersection $z\in C_b\cap C_w$ cannot simultaneously
minimize $\Re F(\,\cdot\,;\boldsymbol x_b)$ on $C_b$ and
maximize $\Re F(\,\cdot\,;\boldsymbol x_w)$ on $C_w$.
Otherwise, it would be a common saddle. The injectivity of
$\zeta_+$, together with conjugation, would then imply
$\boldsymbol x_b=\boldsymbol x_w$, contradicting
$|\boldsymbol x_b-\boldsymbol x_w|\ge\delta$.
Thus
\[
 \bigl[\Re F(z;\boldsymbol x_b)-\mathcal A(\boldsymbol x_b)\bigr]
 +\bigl[\mathcal A(\boldsymbol x_w)-\Re F(z;\boldsymbol x_w)\bigr]
 >0.
\]
The triples $(\boldsymbol x_b,\boldsymbol x_w,z)$ satisfying
$\boldsymbol x_b,\boldsymbol x_w\in\mathcal K$,
$|\boldsymbol x_b-\boldsymbol x_w|\ge\delta$, and
$z\in C_b\cap C_w$ form a compact set. The displayed sum is continuous on this set, so it attains
a positive minimum $c$. Consequently,
\[
 \Re[F(z;\boldsymbol x_w)-F(z;\boldsymbol x_b)]
 \le\mathcal A(\boldsymbol x_w)-\mathcal A(\boldsymbol x_b)-c
\]
uniformly at all intersections.
For $z\in C_b$,
\[
 \Re[F(z;\boldsymbol x_w)-F(z;\boldsymbol x_b)]
 =\frac{\eta_w-\eta_b}{2}\log r_b
  +\frac dN\log|\mathsf w(z)|.
\]
For $d>0$, the residue integral is over $C_b\setminus W$,
where $|\mathsf w(z)|\le r_w$.
For $d\le0$, it is over $C_b\cap W$,
where $|\mathsf w(z)|>r_w$.
Thus, on the closure of the corresponding integration arc,
\[
 \Re[F(z;\boldsymbol x_w)-F(z;\boldsymbol x_b)]
 \le \frac{\eta_w-\eta_b}{2}\log r_b
     +\frac dN\log r_w.
\]
When the contours intersect, this bound is attained at
the arc endpoints, where $|\mathsf w(z)|=r_w$. For $z=r_be^{\ii\theta}\in C_b$,
\eqref{eq:contour-radii} and~\eqref{eq:white-contour} give
\[
 z\in W
 \quad\Longleftrightarrow\quad
 \cos\theta>\frac{\hat\xi(\xi_w,\eta_b)}{\sqrt{1-\eta_b^2}}.
\]
To see this, note that $|\mathsf w(z)|^2>r_w^2$ means
$\bk(1-\xi_w)(r_b^2+1+2r_b\cos\theta)
>(1+\xi_w)(\bk^2r_b^2+1-2\bk r_b\cos\theta)$, and use again
$\bk r_b^2+1=2/(1-\eta_b)$ and
$\bk r_b^2-1=2\eta_b/(1-\eta_b)$.
Also, since $d=m-M$, \eqref{eq:vertex-macroscopic-coordinates}
gives $d=N(\xi_b-\xi_w)/2$; the function $\hat\xi(\xi,\eta)$ is
strictly increasing in $\xi$; and
$|\hat\xi(\xi_b,\eta_b)|<\sqrt{1-\eta_b^2}$ by the liquid-region
condition $\boldsymbol x_b\in L_{\xi,\eta}$.

If $C_b$ and $C_w$ do not intersect, then either
$C_b\cap W=\varnothing$ or $C_b\subset W$.
If $C_b\cap W=\varnothing$, the preceding criterion gives
$\hat\xi(\xi_w,\eta_b)\ge\sqrt{1-\eta_b^2}>\hat\xi(\xi_b,\eta_b)$,
hence $\xi_w>\xi_b$ and $d<0$;
the residue integral is therefore over
$C_b\cap W=\varnothing$.
If $C_b\subset W$, it gives
$\hat\xi(\xi_w,\eta_b)\le-\sqrt{1-\eta_b^2}<\hat\xi(\xi_b,\eta_b)$,
hence $\xi_w<\xi_b$ and $d>0$;
the residue integral is therefore over
$C_b\setminus W=\varnothing$.
Thus $\mathcal R_N=0$ whenever the contours do not intersect.
For $d\ge1$, use
$|R|=r_b^{\nu-n}|\mathsf w|^{d-1}/|\bk z-1|$ on
$|\mathsf w|\le r_w$; for $d\le0$, use the original expression
on $|\mathsf w|\ge r_w$. Their remaining denominators are bounded
away from zero, because
$|\bk z-1|+\bk|z+1|\ge1+\bk$ and
$|z+1|=|\mathsf w(z)|\,|\bk z-1|$. It follows that
\[
 |\mathcal R_N|
 \le C e^{N[\mathcal A(\boldsymbol x_w)-\mathcal A(\boldsymbol x_b)-c]}.
\]

\emph{Integrability at the pole and at infinity.}
The two contours may intersect. We check uniform absolute
integrability after removing the real exponential factor.
In the parameters
$z=r_be^{\ii\theta}$ and $v=r_we^{\ii\phi}$, with
$t=\mathsf w(v)$, we have $z-t=D/(\bk v-1)$, where the Cauchy
denominator is
\[
 D=(\bk z-1)v-(z+1).
\]
Near $D=0$, $z$ stays away from $\bk^{-1}$. For each $\theta$, the
set of $\phi$ satisfying $|D|<s$ has length at most $Cs$ and is
empty unless $|\,|\mathsf w(z)|-r_w|<Cs$. Moreover,
\[
 |\mathsf w(r_be^{\ii\theta})|^2
 =\frac{r_b^2+1+2r_b\cos\theta}
        {\bk^2r_b^2+1-2\bk r_b\cos\theta}.
\]
Near its intersection values, the last condition implies
$|\cos\theta-c(r_b,r_w)|<Cs$. Since
$|\{\theta:|\cos\theta-c|<s\}|\le C\sqrt{s}$ uniformly in $c$,
\[
 |\{(\theta,\phi):|D|<s\}|\le Cs^{3/2},\qquad
 \int_{|D|<s}\frac{d\theta\,d\phi}{|D|}\le C\sqrt{s}.
\]
Since $t/(t+1)=(v+1)/((1+\bk)v)$ and
$dt=-(1+\bk)(\bk v-1)^{-2}dv$, the remaining amplitude,
including the contour differentials, has absolute value
\[
 \frac{|v+1|}{|\bk v-1|\,|D|}\,d\theta\,d\phi.
\]
Only $v=\bk^{-1}$, corresponding to $t=\infty$, needs a further
check. There $D=-(1+\bk)/\bk$ is bounded away from zero, and the
transformed phase gives
\[
 \Re F(t;\boldsymbol x_w)-\mathcal A(\boldsymbol x_w)
 \le c_0\log|\bk v-1|+C_0,
 \qquad c_0=\tfrac12\inf_{\mathcal K}(1-\eta)>0.
\]
On a sufficiently small fixed neighborhood of $v=\bk^{-1}$, this
is at most $-c_1<0$, and, for all sufficiently large $N$,
\[
 \frac{
 e^{\frac N2[\Re F(t;\boldsymbol x_w)-\mathcal A(\boldsymbol x_w)]}}
 {|\bk v-1|}\le C.
\]
Use this half of the exponent to bound the amplitude and the other
half for exponential decay. Away from that neighborhood the amplitude
is bounded by $C/|D|$. The preceding integrability estimate and phase
gaps therefore make the part outside the four saddle neighborhoods
exponentially small, uniformly in the positions. They also justify
the contour limits above, including at tangencies.

\emph{The four Gaussian contributions.}
Parametrize $C_b$ counterclockwise and the $v$-circle clockwise.
The tangent vectors at $\zeta_\alpha(\boldsymbol x_b)$ and
$\zeta_\beta(\boldsymbol x_w)$ are
$d_b\bigl(\zeta_\alpha(\boldsymbol x_b)\bigr)$ and
$d_w\bigl(\zeta_\beta(\boldsymbol x_w)\bigr)$, respectively;
for the latter, $t=\mathsf w(v)$ and $v=r_we^{-\ii\phi}$ give
$dt/d\phi=\ii(1+\bk)v(\bk v-1)^{-2}
=\ii(t+1)(\bk t-1)/(1+\bk)$.
Taylor expansion and the elementary Gaussian integral
therefore give, for each pair of saddles,
\[
 \begin{aligned}
 &-\frac{e^{N[\mathcal F_\beta(\boldsymbol x_w)
                    -\mathcal F_\alpha(\boldsymbol x_b)]}}{2\pi N}
 \\
 &\quad\times
 \frac{\zeta_\beta(\boldsymbol x_w)}
 {\zeta_\alpha(\boldsymbol x_b)
  \bigl(\zeta_\beta(\boldsymbol x_w)+1\bigr)
  \bigl(\zeta_\alpha(\boldsymbol x_b)
        -\zeta_\beta(\boldsymbol x_w)\bigr)}
 \frac{d_b\bigl(\zeta_\alpha(\boldsymbol x_b)\bigr)
       d_w\bigl(\zeta_\beta(\boldsymbol x_w)\bigr)}
      {\sqrt{H_{b,\alpha}}\sqrt{H_{w,\beta}}}.
 \end{aligned}
\]
The Taylor remainder is uniformly $O(N^{-3/2})$ after the real
exponential factor is removed. All four saddle pairs stay a positive
distance from the pole $z=t$, so their amplitudes and derivatives
are uniformly bounded. Summing proves the expansion.
The coefficient derivative bounds follow from its explicit formula,
the nonvanishing of $H_{b,\alpha}$ and $H_{w,\beta}$,
and separation of the saddle images.
\end{proof}

\begin{lemma}\label{lem:bulk-oscillatory-sum}
Let $\boldsymbol d\in\{(1,0),(-1,0),(0,1),(0,-1)\}$, and let
$\boldsymbol x_n=\boldsymbol x_0+2n\boldsymbol d/N$,
$0\le n\le L$, lie on a segment in $\mathcal K$, with $L\le C_0N$.
Suppose that $a_N$ is continuously differentiable on that segment and
\[
 \sup\bigl(|a_N|+|\partial_{\boldsymbol d}a_N|\bigr)\le C_1.
\]
Then, for either sign,
\begin{equation}\label{eq:bulk-oscillatory-sum}
 \left|\frac1N\sum_{n=0}^{L-1}
 a_N(\boldsymbol x_n)
 e^{\pm\ii N\theta_{\rm osc}(\boldsymbol x_n)}\right|
 \le\frac C N.
\end{equation}
The constant $C$ depends only on $\mathcal K$, $C_0$, and $C_1$.
\end{lemma}

\begin{proof}
Since $F_\zeta(\zeta_+;\boldsymbol x)=0$, differentiation gives
\begin{equation}\label{eq:bulk-phase-derivatives}
 \partial_\eta\theta_{\rm osc}=\arg\zeta_+,
 \qquad
 \partial_\xi\theta_{\rm osc}=-\arg w(\zeta_+).
\end{equation}
Here $0<\arg\zeta_+<\pi$ and $-\pi<\arg w(\zeta_+)<0$, because
\[
 \Im w(\zeta_+)=-\frac{(1+\bk)\Im\zeta_+}{|\bk\zeta_+-1|^2}<0.
\]
Compactness therefore gives
\[
 \inf_{\boldsymbol x\in\mathcal K}
 \left|e^{2\ii\partial_{\boldsymbol d}
                         \theta_{\rm osc}(\boldsymbol x)}-1\right|>0.
\]
Set $u_n=e^{\pm\ii N\theta_{\rm osc}(\boldsymbol x_n)}$ and
$r_n=u_{n+1}/u_n$, for $0\le n<L$.
Taylor's formula and the derivative bounds on $\theta_{\rm osc}$ imply
\[
 r_n=e^{\pm2\ii\partial_{\boldsymbol d}
                         \theta_{\rm osc}(\boldsymbol x_n)}+O(N^{-1}),
 \qquad |r_n-1|\ge c>0,
 \qquad |r_n-r_{n-1}|\le C/N
\]
for all sufficiently large $N$.
Writing $b_n=a_N(\boldsymbol x_n)/(r_n-1)$, we obtain
$|b_n|\le C$ and $|b_n-b_{n-1}|\le C/N$.
Summation by parts gives
\[
 \begin{aligned}
 \sum_{n=0}^{L-1}a_N(\boldsymbol x_n)u_n
 &=\sum_{n=0}^{L-1}b_n(u_{n+1}-u_n)\\
 &=b_{L-1}u_L-b_0u_0
   +\sum_{n=1}^{L-1}(b_{n-1}-b_n)u_n.
 \end{aligned}
\]
Its absolute value is at most $C(1+L/N)\le C$.
Dividing by $N$ proves~\eqref{eq:bulk-oscillatory-sum}.
The empty sum is zero, and the remaining bounded values of $N$
are covered by increasing $C$.
\end{proof}

\begin{lemma}\label{lem:scalar-cutoff-cycle}
With $C_i^\varepsilon$ as defined in~\eqref{eq:cutoff-contours},
for every fixed $k\ge2$ and labels $i_1,\ldots,i_k\in\{1,\ldots,m\}$
with $i_s\ne i_{s+1}$, where $i_{k+1}=i_1$, one has
\begin{equation}\label{eq:scalar-cutoff-cycle}
 \begin{aligned}
 &\sum_{\substack{e_s\in E_{i_s}\cap E_{N,\varepsilon}\\1\le s\le k}}
   \prod_{s=1}^k C_{i_si_{s+1},N}(e_s,e_{s+1})\\
 &\qquad\longrightarrow
 \frac1{(2\pi\ii)^k}
 \int_{C_{i_1}^\varepsilon\times\cdots\times C_{i_k}^\varepsilon}
       \prod_{s=1}^k\frac{dq_s}{q_s-q_{s+1}},
 \qquad e_{k+1}=e_1,\quad q_{k+1}=q_1.
 \end{aligned}
\end{equation}
\end{lemma}

\begin{proof}
By Lemma~\ref{lem:bulk-four-saddle} every factor of a cyclic product is
a sum of four oscillating terms. The terms in which some edge carries
two different saddles cancel by oscillation, and the remaining terms
are Riemann sums for the integral in~\eqref{eq:scalar-cutoff-cycle}.

\emph{Expansion of the cyclic product.}
By~\eqref{eq:cross-kernel}, each summand on the left is
\[
 \prod_{s=1}^k
 \left[K_N^{-1}(b_{e_s},w_{e_{s+1}})
       \sigma_{e_{s+1}}K_N(w_{e_{s+1}},b_{e_{s+1}})\right].
\]
Adjacent labels differ, so Lemma~\ref{lem:bulk-four-saddle}
applies to every inverse entry in this product.

The real exponential factors in a cyclic product satisfy
\[
 \begin{aligned}
 &\sum_{s=1}^k
 \bigl[\Re\mathcal F_+(\boldsymbol x_{w_{e_{s+1}}})
       -\Re\mathcal F_+(\boldsymbol x_{b_{e_s}})\bigr]\\
 &\qquad=
 \sum_{s=1}^k
 \bigl[\Re\mathcal F_+(\boldsymbol x_{w_{e_s}})
       -\Re\mathcal F_+(\boldsymbol x_{b_{e_s}})\bigr]
 =O(N^{-1}).
 \end{aligned}
\]
Thus each cyclic product is $O(N^{-k})$, and replacing all inverse
entries by the four leading terms in~\eqref{eq:four-saddle-entry}
changes its sum by $O(N^{-1/2})$: the error per tuple is
$O(N^{-k-1/2})$, and there are $O(N^k)$ tuples.
Tuples with a repeated edge at nonadjacent positions number
$O(N^{k-1})$ and contribute $O(N^{-1})$.

For each edge $e_s$, use $\boldsymbol x_{b_{e_s}}$ as its reference
position. On each fixed edge type, the white endpoint has coordinates
$\boldsymbol x_{b_{e_s}}+\boldsymbol\delta/N$ for a fixed bounded
offset $\boldsymbol\delta$. The displacement factor is controlled by
\[
 N\bigl[\mathcal F_\nu(\boldsymbol x+\boldsymbol\delta/N)
           -\mathcal F_\nu(\boldsymbol x)\bigr]
 =\int_0^1\nabla\mathcal F_\nu
      (\boldsymbol x+t\boldsymbol\delta/N)
      \cdot\boldsymbol\delta\,dt,
 \qquad \nu\in\{+,-\}.
\]
This expression and its first derivatives are uniformly bounded.
Its exponential can therefore be included in the coefficient.

Let $\alpha_s$ and $\beta_s$ be the saddle choices at the black and
white endpoints of $e_s$. After this inclusion, each choice of branches
contributes a sum of the form
\begin{equation}\label{eq:expanded-branch-cycle}
 \frac1{N^k}
 \sum_{\substack{e_s\in E_{i_s}\cap E_{N,\varepsilon}\\1\le s\le k}}
 a_N(e_1,\ldots,e_k)
 \prod_{s\in\mathcal D}
 e^{\tau_s\ii N\theta_{\rm osc}(\boldsymbol x_{b_{e_s}})},
 \qquad
 \mathcal D=\{s:\alpha_s\ne\beta_s\},\quad
 \tau_s\in\{1,-1\}.
\end{equation}
Here $a_N$ and its first derivatives are uniformly bounded after
fixing the finitely many edge types. Only adjacent positions enter
each inverse entry; nonadjacent positions may coincide without
creating a pole.

\emph{Cancellation of mixed branches.}
Suppose $\mathcal D\ne\varnothing$, and fix $s\in\mathcal D$.
Split the $s$th zipper into its finitely many straight segments and
edge types. Within each type, successive full periods are spaced by
$2/N$ in one coordinate. Lemma~\ref{lem:bulk-oscillatory-sum} gives
\begin{equation}\label{eq:cycle-inner-oscillation}
 \sup_{(e_j)_{j\ne s}}
 \left|\frac1N
 \sum_{e_s\in E_{i_s}\cap E_{N,\varepsilon}}
 a_N(e_1,\ldots,e_k)
 e^{\tau_s\ii N\theta_{\rm osc}(\boldsymbol x_{b_{e_s}})}\right|
 \le\frac C N.
\end{equation}
The bounded number of corner and incomplete-period edges contributes
$O(N^{-1})$ to this normalized sum. If repeated-edge tuples were
removed, adding back at most $k-1$ excluded sites has the same bound.
The other exponential factors in~\eqref{eq:expanded-branch-cycle}
have modulus $1$ and are independent of $e_s$. Hence that entire sum
has absolute value at most
\[
 \frac{O(N^{k-1})}{N^{k-1}}\frac C N=O(N^{-1}).
\]
There are finitely many branch choices for fixed $k$, so all choices
with $\mathcal D\ne\varnothing$ together contribute $o(1)$.

\emph{One period of a zipper.}
It remains to sum the choices with $\alpha_s=\beta_s$ at every
position. The map~\eqref{eq:diagonal-coordinates} reverses orientation,
so $\sigma_e=+1$ when the white endpoint is to the left of the
oriented path in $(u,v)$.

Each positive period in~\eqref{eq:zipper-periods} crosses two edges with $K_N(w,b)=1$ and $\sigma_e=-1$;
by~\eqref{eq:scalar-kasteleyn-entries}, the down edges
$v^{\rm W}_{m,n}v^{\rm B}_{m+1,n}$, which carry the entry $-\bk$,
are not crossed by these periods.
Their white-minus-black offsets
$\boldsymbol\delta=N(\boldsymbol x_w-\boldsymbol x_b)$ are
\[
\begin{array}{c|c|c}
\text{period}&\text{crossed edge}&\boldsymbol\delta\\\hline
\xi&v^{\rm W}_{m,n}v^{\rm B}_{m,n+1}&(0,-2)\\
\xi&v^{\rm W}_{m,n}v^{\rm B}_{m+1,n+1}&(-2,-2)\\
\eta&v^{\rm W}_{m,n}v^{\rm B}_{m,n+1}&(0,-2)\\
\eta&v^{\rm W}_{m,n+1}v^{\rm B}_{m,n+1}&(0,0)
\end{array}
\]
At a saddle $\zeta$, stationarity gives
$\partial_\xi\mathcal F=-\tfrac12\log w(\zeta)$ and
$\partial_\eta\mathcal F=\tfrac12\log\zeta$.
The corresponding displacement factors
$\exp(\boldsymbol\delta\cdot\nabla\mathcal F)$ are therefore
$\zeta^{-1}$, $w(\zeta)\zeta^{-1}$, and $1$, respectively.
Thus the signed period sums are
\[
 E_\xi(\zeta)=-\frac{1+w(\zeta)}\zeta,
 \qquad E_\eta(\zeta)=-\left(1+\frac1\zeta\right).
\]

Put $H_\nu(\boldsymbol x)=
F_{\zeta\zeta}(\zeta_\nu(\boldsymbol x);\boldsymbol x)$.
In the exact scalar inverse, the rational prefactor is
$t/[(t+1)z(z-t)]$. In a consistent cyclic product, the factors
$t/(t+1)$ and $1/z$ at each position combine to $1/(\zeta+1)$.
The oriented Gaussian factors in Lemma~\ref{lem:bulk-four-saddle},
evaluated at $\boldsymbol x_b=\boldsymbol x_w=\boldsymbol x$
with $\alpha=\beta=\nu$, satisfy
\[
 \mathfrak c_{b,\nu}\mathfrak c_{w,\nu}
 =\frac{\epsilon_\nu\ii}{H_\nu(\boldsymbol x)},
 \qquad \epsilon_+=1,\quad\epsilon_-=-1.
\]
Both sides vary continuously in the liquid region $L_{\xi,\eta}$ and
their quotient has square $1$, because
$H_{b,\nu}H_{w,\nu}=-H_\nu^2d_b^2d_w^2$.
At $(\xi,\eta)=(0,0)$, the saddles are $\pm\ii\bk^{-1/2}$.
Put $\varsigma=\sqrt\bk-1/\sqrt\bk$. For $\nu=+$,
\[
 H_+=\frac{\bk}{2-\ii\varsigma},\qquad
 d_b(\zeta_+)=-\frac1{\sqrt\bk},\qquad
 d_w(\zeta_+)=-\frac{\ii(2-\ii\varsigma)}{1+\bk},
\]
so that
$H_{b,+}=(2-\ii\varsigma)^{-1}$,
$H_{w,+}=\bk(2-\ii\varsigma)/(1+\bk)^2$,
$\sqrt{H_{b,+}}\sqrt{H_{w,+}}=\sqrt\bk/(1+\bk)$, and
$\mathfrak c_{b,+}\mathfrak c_{w,+}
=\ii(2-\ii\varsigma)/\bk=\ii/H_+$.
The case $\nu=-$ follows by complex conjugation, since
$\zeta_-=\overline{\zeta_+}$,
$d_b(\zeta_-)d_w(\zeta_-)=\overline{d_b(\zeta_+)d_w(\zeta_+)}$,
and $H_{b,-}=\overline{H_{b,+}}$, $H_{w,-}=\overline{H_{w,+}}$.
This fixes the two signs throughout $L_{\xi,\eta}$.
For $\bk=1$ these values are $\zeta_\pm=\pm\ii$, $H_\pm=1/2$,
black tangents $-1,+1$ and both white tangents $-\ii$.

\emph{The limiting differential.}
For a positive period in direction $d\in\{\xi,\eta\}$,
the contribution at one position is
\[
 \frac{\epsilon_\nu}{2\pi\mathrm i N H_\nu}
 \frac{E_d(\zeta_\nu)}{\zeta_\nu+1}.
\]
Since $w(\zeta)=(\zeta+1)/(\bk\zeta-1)$, we have
$1+w(\zeta)=(1+\bk)\zeta/(\bk\zeta-1)$ and
\[
 \frac{E_\xi(\zeta)}{\zeta+1}
 =-\frac{1+\bk}{(\zeta+1)(\bk\zeta-1)},
 \qquad
 \frac{E_\eta(\zeta)}{\zeta+1}=-\frac1\zeta.
\]
Implicit differentiation of the saddle equation
$F_\zeta(\zeta_\nu(\boldsymbol x);\boldsymbol x)=0$, with
$F_{\zeta\xi}=(1+\bk)/[2(\zeta+1)(\bk\zeta-1)]$ and
$F_{\zeta\eta}=1/(2\zeta)$, gives
\[
 \partial_\xi\zeta_\nu
 =-\frac{1+\bk}{2(\zeta_\nu+1)(\bk\zeta_\nu-1)H_\nu},\qquad
 \partial_\eta\zeta_\nu=-\frac1{2\zeta_\nu H_\nu}.
\]
For the period length $h_N=2/N$, the summed factor is consequently
$\epsilon_\nu h_N\partial_d\zeta_\nu/(2\pi\mathrm i)$.
Reversing a period changes both $\sigma_e$ and its increment's sign.
Thus the two branches give the differentials
$d\zeta_+/(2\pi\mathrm i)$ and
$-d\overline{\zeta_+}/(2\pi\mathrm i)$ along the outward path.

\emph{Summation over periods.}
Fix one full period on each zipper. Let $\boldsymbol x_s$ be
its initial macroscopic position and let
$\boldsymbol d_s\in\{(1,0),(-1,0),(0,1),(0,-1)\}$ be its direction.
For a consistent branch choice $\nu_s\in\{+,-\}$, put
\[
 \zeta_{s,\nu}=\zeta_\nu(\boldsymbol x_s),\qquad
 \Delta\zeta_{s,+}
 =\zeta_+(\boldsymbol x_s+h_N\boldsymbol d_s)
   -\zeta_+(\boldsymbol x_s),\qquad
 \Delta\zeta_{s,-}=-\overline{\Delta\zeta_{s,+}}.
\]
Use cyclic indices $\boldsymbol x_{k+1}=\boldsymbol x_1$ and
$\nu_{k+1}=\nu_1$. Let $\mathfrak p_s$ be the two crossed edges,
and let $a_N^{\boldsymbol\nu}$ be the coefficient in
\eqref{eq:expanded-branch-cycle} with $\alpha_s=\beta_s=\nu_s$.
Replacing positions within a period by $\boldsymbol x_s$ changes
these bounded coefficients by $O(N^{-1})$. Therefore
\begin{equation}\label{eq:period-cycle-sum}
 \frac1{N^k}\sum_{\substack{e_s\in\mathfrak p_s\\1\le s\le k}}
 a_N^{\boldsymbol\nu}(e_1,\ldots,e_k)
 =
 \frac1{(2\pi\mathrm i)^k}
 \prod_{s=1}^k
 \frac{\Delta\zeta_{s,\nu_s}}
      {\zeta_{s,\nu_s}-\zeta_{s+1,\nu_{s+1}}}
 +O(N^{-k-1}).
\end{equation}
There are $O(N^k)$ tuples of full periods, so the errors sum to
$O(N^{-1})$. Incomplete periods contribute $O(N^{-1})$ as well.
The leading terms are Riemann sums. Summing over the $2^k$ branch
choices gives the integral over the chains
$C_i^\varepsilon=\gamma_i^\varepsilon-\overline{\gamma_i^\varepsilon}$
in~\eqref{eq:scalar-cutoff-cycle}.
\end{proof}

\begin{proposition}[Fixed-cutoff limit]\label{prop:interior-determinant}
For every fixed sufficiently small $\varepsilon>0$ and every compact
$\mathcal C\subset X_{\rm unip}$,
\begin{equation}\label{eq:interior-determinant-limit}
 \sup_{\rho\in\mathcal C}
 |F_{N,\varepsilon}(\rho)-D_\varepsilon(\rho)|\longrightarrow0.
\end{equation}
\end{proposition}

\begin{proof}
By~\eqref{eq:cyclic-traces} and Lemma~\ref{lem:scalar-cutoff-cycle},
for every fixed $k\ge2$,
\begin{align}
 \Tr O_{N,\varepsilon}^k\longrightarrow
 \Tr\mathcal O_\varepsilon^k
 &=\frac1{(2\pi\ii)^k}
 \sum_{\substack{i_1,\ldots,i_k\in\{1,\ldots,m\}\\
                 i_s\ne i_{s+1},\ 1\le s\le k}}
 \Tr(A_{i_2}\cdots A_{i_k}A_{i_1})\notag\\[-2mm]
 &\quad\times
 \int_{C_{i_1}^\varepsilon\times\cdots\times C_{i_k}^\varepsilon}
       \prod_{s=1}^k\frac{dq_s}{q_s-q_{s+1}},
 \qquad i_{k+1}=i_1,\quad q_{k+1}=q_1.
 \label{eq:interior-cycle-limit}
\end{align}
The integral is the trace of the composed kernel
\eqref{eq:canonical-cutoff-kernel} and converges absolutely by the
separation of adjacent cuts. For each label sequence the ordered
matrix factor is independent of the edge positions. There are
finitely many such factors, uniformly bounded on $\mathcal C$, so
\eqref{eq:interior-cycle-limit} holds uniformly there.

We next give a uniform Hilbert--Schmidt bound after a scalar
similarity. Put
\[
 \mathfrak s(\boldsymbol x)=\Re\mathcal F_+(\boldsymbol x),\qquad
 (\mathsf G_{N,\varepsilon})_{ee}
 =e^{N\mathfrak s(\boldsymbol x_{b_e})}I_2,
 \qquad e\in E_{N,\varepsilon},
\]
and, on the retained edge coordinates, set
\[
 \widetilde O_{N,\varepsilon}
 =\mathsf G_{N,\varepsilon}O_{N,\varepsilon}
                        \mathsf G_{N,\varepsilon}^{-1}.
\]
Deleting the zero rows and columns outside $E_{N,\varepsilon}$
does not change the traces or determinants.
Equation~\eqref{eq:four-saddle-entry} and~\eqref{eq:cross-kernel} give
\[
 |(O_{N,\varepsilon})_{e\alpha,f\beta}|
 \le\frac{C_{\varepsilon,\mathcal C}}N
 e^{N[\mathfrak s(\boldsymbol x_{w_f})
             -\mathfrak s(\boldsymbol x_{b_e})]}.
\]
The same-label blocks are zero. Since $w_f,b_f$ are adjacent,
\[
 |\mathfrak s(\boldsymbol x_{w_f})
       -\mathfrak s(\boldsymbol x_{b_f})|\le C_\varepsilon/N,
 \qquad
 |(\widetilde O_{N,\varepsilon})_{e\alpha,f\beta}|
 \le C_{\varepsilon,\mathcal C}/N.
\]
There are $O(N)$ retained edges and two matrix coordinates per edge;
summing the squares of all entries therefore gives
\[
 \sup_{N\ge N_0,\rho\in\mathcal C}
 \|\widetilde O_{N,\varepsilon}(\rho)\|_2<\infty.
\]
The continuum bound is equally explicit:
\[
 \|\mathcal O_\varepsilon(\rho)\|_2^2
 =\frac1{4\pi^2}\sum_{i\ne j}
 \left(\sum_{\alpha,\beta=1}^2|(A_j)_{\alpha\beta}|^2\right)
 \int_{C_i^\varepsilon}\int_{C_j^\varepsilon}
       \frac{|dq|\,|dr|}{|q-r|^2}.
\]
The lengths are finite and the distinct contours are separated, so
this is uniformly bounded on $\mathcal C$.
Choose $H<\infty$ bounding both Hilbert--Schmidt norms.

Write
\[
 \det{}_2(I+t\widetilde O_{N,\varepsilon}(\rho))
 =\sum_{n=0}^\infty d_{N,n}(\rho)t^n,
 \qquad
 \det{}_2(I+t\mathcal O_\varepsilon(\rho))
 =\sum_{n=0}^\infty d_n(\rho)t^n.
\]
For any of these operators $T$, its powers $T^j$, $j\ge2$, are
trace class and satisfy
\[
 |\Tr T^j|\le\|T^2\|_1\|T\|^{j-2}\le H^j.
\]
Thus, for $|t|<(1+H)^{-1}$,
\[
 \log\det{}_2(I+tT)
 =\sum_{j=2}^\infty\frac{(-1)^{j+1}}j\,t^j\Tr T^j.
\]
Differentiating and comparing coefficients gives
\begin{equation}\label{eq:det-two-coefficient-recursion}
 d_0=1,\qquad d_1=0,\qquad
 n d_n=\sum_{j=2}^n(-1)^{j+1}\Tr(T^j)d_{n-j},\qquad n\ge2.
\end{equation}
The same recurrence holds for $d_{N,n}$.
Similarity preserves traces, so~\eqref{eq:interior-cycle-limit}
and induction in~\eqref{eq:det-two-coefficient-recursion} give
\[
 \sup_{\rho\in\mathcal C}|d_{N,n}(\rho)-d_n(\rho)|\longrightarrow0
 \qquad\text{for every fixed }n.
\]

For $R>1$, the bound
$|\det{}_2(I+tT)|\le\exp(|t|^2\|T\|_2^2/2)$ and Cauchy's estimate
bound the $n$th coefficient of either determinant by
$e^{R^2H^2/2}R^{-n}$. This summable majorant gives uniform convergence
at $t=1$. Similarity preserves $\det{}_2$, proving
\eqref{eq:interior-determinant-limit}.
\end{proof}

\subsection{Direct estimates at the arctic boundary}
\label{subsec:arctic-edge}

The arctic boundary $\partial L_{\xi,\eta}$ touches the boundary of the
square $[-1,1]^2$ at the four tangency points
\begin{equation}\label{eq:tangency-points}
 \Bigl(\pm1,\ \mp\frac{\bk-1}{\bk+1}\Bigr),
 \qquad
 \Bigl(\mp\frac{\bk-1}{\bk+1},\ \pm1\Bigr),
\end{equation}
which for $\bk=1$ are the four points with $\xi\eta=0$. All other
points of $\partial L_{\xi,\eta}$ are called regular arctic points; at a
regular arctic point the tangent of $\partial L_{\xi,\eta}$ is neither
horizontal nor vertical in the coordinates $(\xi,\eta)$.
Choose the terminal crossings at distinct regular arctic points.
Continue each terminal
segment through its frozen region to the actual boundary.
Put $\delta_N=N^{-2/3+1/100}$ and stop each zipper at a face
$r_{i,N}$ at outward transverse distance $\delta_N+O(N^{-1})$.
Here transverse distance is measured in $\eta$ along its terminal
vertical segment. The chosen crossings are transverse, so this
distance is comparable to Euclidean distance near the boundary.
Let $E_N^{\rm edge}$ be the retained crossed edges and
$O_N^{\rm edge}$ the corresponding compression of $O_N$.
Bounded lattice displacements do not change the estimates below.

\paragraph{The common scalar gauge.}
Write $\boldsymbol x=(\xi,\eta)$ and set
\begin{equation}\label{eq:explicit-entropy-gauge}
 \begin{aligned}
 \mathcal A(\xi,\eta)
 ={}&\frac14\bigl[(1+\eta)\log(1+\eta)+(1-\eta)\log(1-\eta)
 -(1+\xi)\log(1+\xi)-(1-\xi)\log(1-\xi)\bigr]\\
 &-\frac12\log(1+\bk)+\frac{\xi-\eta}{4}\log\bk.
 \end{aligned}
\end{equation}
In the liquid region, substitution of
$\zeta_+(\boldsymbol x)$ into the phase gives
$\mathcal A(\boldsymbol x)
 =\Re F(\zeta_+(\boldsymbol x);\boldsymbol x)$.
The formula extends smoothly to $(-1,1)^2$, with
\[
 \partial_\xi\mathcal A=-\frac14\log\frac{1+\xi}{\bk(1-\xi)},
 \qquad
 \partial_\eta\mathcal A=\frac14\log\frac{1+\eta}{\bk(1-\eta)}.
\]
The retained paths stay away from the four tangency points, so these
derivatives are bounded in their neighborhoods. Consequently
\begin{equation}\label{eq:adjacent-action}
 |\mathcal A(\boldsymbol x_{w_e})-
   \mathcal A(\boldsymbol x_{b_e})|\le C/N.
\end{equation}
Use the scalar gauge
\begin{equation}\label{eq:edge-gauge}
 H_N(e,e)=e^{N\mathcal A(\boldsymbol x_{b_e})}I_2,
 \qquad\widetilde O_N^{\rm edge}=H_NO_N^{\rm edge}H_N^{-1}.
\end{equation}

\begin{lemma}\label{lem:scalar-edge-bound}
In each terminal arctic neighborhood, let $d$ be the signed inward
distance along the transverse coordinate. On the retained paths define
\begin{equation}\label{eq:edge-envelope}
 a_N(\boldsymbol x)=
 \begin{cases}
 N^{-1/2}d^{-1/4},&d>2N^{-2/3}\text{ in a terminal neighborhood},\\
 N^{-1/3},&-2\delta_N\le d\le 2N^{-2/3}\text{ in that neighborhood},\\
 N^{-1/2},&\text{elsewhere on the retained paths}.
 \end{cases}
\end{equation}
For an edge $e$, write $d_e=d(\boldsymbol x_{b_e})$ and
$a_N(e)=a_N(\boldsymbol x_{b_e})$.
Then
\begin{equation}\label{eq:factorized-edge-bound}
 \|\widetilde O_N^{\rm edge}(e,f)\|_{\rm HS}
 \le C\|A_{i(f)}\|_{\rm HS}a_N(e)a_N(f),
 \qquad i(e)\ne i(f).
\end{equation}
The constant is uniform in $N$ and in the retained edges.
\end{lemma}

\begin{proof}
We first bound the oscillatory integral over one contour, by van der
Corput estimates near the coalescing saddles, and then combine the two
contours in the formula~\eqref{eq:uniform-deformed-inverse}.

\emph{The black contour.}
Consider a small fixed neighborhood of a regular arctic point
$(\xi_0,\eta_0)$, that is, a point with
$\hat\xi(\xi_0,\eta_0)^2+\eta_0^2=1$ which is not one of the four
points~\eqref{eq:tangency-points}. In particular $|\eta_0|<1$ and
$\hat\xi(\xi_0,\eta_0)=\pm\sqrt{1-\eta_0^2}\ne0$.
Let $d$ be signed inward transverse distance to
$\partial L_{\xi,\eta}$, positive on the liquid side. Put
\[
 \rho^2=\frac{1+\eta}{\bk(1-\eta)},\qquad
 u_0=\frac{\hat\xi(\xi,\eta)}{\sqrt{1-\eta^2}},\qquad
 z=\rho e^{\ii\theta}.
\]
On this circle, $g(u)=\Re F(\rho e^{\ii\theta};\boldsymbol x)$,
$u=\cos\theta$, satisfies
\[
 g(u)=\frac{1+\eta}{2}\log\rho
 -\frac{1+\xi}{4}\log(1+\rho^2+2\rho u)
 -\frac{1-\xi}{4}\log(1+\bk^2\rho^2-2\bk\rho u).
\]
Its second derivative is positive, its unique critical point is $u_0$,
and $g(u_0)=\mathcal A(\xi,\eta)$. These statements remain valid when
$u_0$ is slightly outside $[-1,1]$: the logarithm arguments at $u_0$
are $(1+\bk)(1+\xi)/[\bk(1-\eta)]$ and
$(1+\bk)(1-\xi)/(1-\eta)$, both positive.
Thus
\[
 h_N(\theta):=\exp\{N[\mathcal A-g(\cos\theta)]\}\le1.
\]
Convexity of $g$ and monotonicity of cosine show that $h_N$ is monotone
on at most four intervals. Its total variation is therefore uniformly
bounded.

\emph{Coalescing saddles.}
Put $p(\theta)=\Im[-F(\rho e^{\ii\theta};\boldsymbol x)]$.
Direct differentiation gives
\begin{equation}\label{eq:circle-phase-derivative}
 p'(\theta)=
 \frac{\bk\rho^2(\cos\theta-u_0)
       \bigl[2\eta\cos\theta+(\bk-1)\rho(1-\eta)\bigr]}
      {|\rho e^{\ii\theta}+1|^2\,|\bk\rho e^{\ii\theta}-1|^2}.
\end{equation}
Indeed $p'(\theta)=-\Re[zF_\zeta(z;\boldsymbol x)]$ is a rational
function of $u=\cos\theta$ with the displayed denominator and a
quadratic numerator, which vanishes at $u_0$ and has leading
coefficient $2\bk\rho^2\eta$. For $\bk=1$ the bracket is
$2\eta\cos\theta$.
The coalescing saddles occur near $\theta_*=0$ if
$\hat\xi(\xi_0,\eta_0)>0$ and
near $\theta_*=\pi$ if $\hat\xi(\xi_0,\eta_0)<0$.
At $(\xi_0,\eta_0)$ and $\theta=\theta_*$, the bracket
in~\eqref{eq:circle-phase-derivative} equals
$\pm2\eta_0+(\bk-1)\sqrt{(1-\eta_0^2)/\bk}$, with the sign of
$\hat\xi(\xi_0,\eta_0)$. A direct computation shows that on
$\partial L_{\xi,\eta}$ this quantity vanishes only at the two tangency
points with $\xi_0=\pm1$, where the double saddle is $\bk^{-1}$ or
$-1$; at every regular arctic point the double saddle
$\pm\rho$ is different from $-1$ and $\bk^{-1}$, so the denominator
in~\eqref{eq:circle-phase-derivative} is bounded away from zero near
$\theta_*$.
Consequently, a fixed sufficiently small neighborhood of
$\theta_*$ has $|p'''|\ge c>0$ for all parameters sufficiently close
to $(\xi_0,\eta_0)$. The third-derivative van der Corput estimate,
applied with the bounded-variation amplitude $h_N$ times a smooth
factor, therefore gives an $O(N^{-1/3})$ bound.

\emph{Separated saddles on the liquid side.}
For $d\ge N^{-2/3}$ on the liquid side, the two saddle angles have
separation comparable to $\sqrt d$. On neighborhoods of radius
$c_1\sqrt d$ around them, \eqref{eq:circle-phase-derivative} gives
$|p''|\ge c\sqrt d$. The second-derivative estimate gives
$O(N^{-1/2}d^{-1/4})$. On the rest of the fixed neighborhood of
$\theta_*$, $|p'|\ge cd$ and $p'$ is monotone on a uniformly bounded
number of intervals. Integration by parts gives $O((Nd)^{-1})$,
which is bounded by $O(N^{-1/2}d^{-1/4})$ because $d\ge N^{-2/3}$.
The parts outside this fixed neighborhood are exponentially small:
there $g(\cos\theta)-\mathcal A\ge c>0$.

Here the versions of van der Corput used above require only a
uniform bound on the supremum and total variation of the amplitude.
They follow by splitting where $|p'|\le t$, using the sublevel bound
$|\{|p'|\le t\}|\le C(t/\lambda)^{1/(r-1)}$ when
$|p^{(r)}|\ge\lambda$, and integrating by parts on the complement.
Optimizing $t$ gives $C(N\lambda)^{-1/r}$, for $r=2,3$.
The number of monotonicity intervals here is bounded by
\eqref{eq:circle-phase-derivative}.

\emph{The white contour.}
The white contour satisfies the same estimates. Indeed, under
$v=w(t)=(t+1)/(\bk t-1)$, by~\eqref{eq:action-duality},
\[
 F(t;\boldsymbol x)=-F(v;(\eta,\xi))-\log(1+\bk),
 \qquad
 \mathcal A(\xi,\eta)=-\mathcal A(\eta,\xi)-\log(1+\bk),
\]
for real parts and locally compatible logarithms. Apply the preceding
argument on $|v|^2=(1+\xi)/[\bk(1-\xi)]$, at the point
$(\eta,\xi)$; the set $L_{\xi,\eta}$, its four tangency points, and
hence the set of regular arctic points are invariant under
$(\xi,\eta)\mapsto(\eta,\xi)$. All changes of variable are
uniformly regular near the selected arctic point.

\emph{The double integral.}
Use the black circles and white Apollonius circles in the exact
deformed-contour identity~\eqref{eq:uniform-deformed-inverse}.
For endpoints on distinct retained zippers, their saddle neighborhoods
are separated, including at the distinct terminal arctic points.
In a product of two such neighborhoods, the Cauchy denominator and
its derivatives are bounded. Apply the one-dimensional estimate first
to one variable, also applying it after differentiating the smooth
factor in the other variable. The resulting function and its first
derivative are bounded by $Ca_N(\boldsymbol x_w)$. A second application in
the other variable bounds the local part of the gauged inverse by
$Ca_N(\boldsymbol x_b)a_N(\boldsymbol x_w)$.

\emph{Residues and remote contours.}
Write
$F_v(\zeta)=F(\zeta;\boldsymbol x_v)$ and
$\mathcal A_v=\mathcal A(\boldsymbol x_v)$ for $v\in\{b,w\}$.
At a contour intersection,
\[
 \Re F_w\le\mathcal A_w,\qquad
 \Re F_b\ge\mathcal A_b.
\]
Equality in both would require a common saddle, which is excluded by
separation of the zipper spectral images. Compactness gives a fixed
strict gap. When the contours intersect, the selected residue arc has
its maximal real exponent at an intersection, as proved in the contour
identity.

A full-circle residue can also occur in the shrinking outward collar.
Write $r_b^2=(1+\eta_b)/[\bk(1-\eta_b)]$,
$R_w^2=(1+\xi_w)/[\bk(1-\xi_w)]$, and
$t_b=\sqrt{1-\eta_b^2}$. If the selected residue is the whole
black circle, then, by the criterion for $z\in W$ in the proof of
Lemma~\ref{lem:bulk-four-saddle}, either
\[
 \xi_b>\xi_w\ \text{ and }\ \hat\xi(\xi_w,\eta_b)\ge t_b,
 \qquad\text{or}\qquad
 \xi_b\le\xi_w\ \text{ and }\ \hat\xi(\xi_w,\eta_b)\le-t_b.
\]
Since $\hat\xi(\xi,\eta)$ is increasing in $\xi$, in both cases
$|\hat\xi(\xi_b,\eta_b)|\ge t_b$, so the black vertex lies in the
outward collar.
Along any convergent sequence of such pairs, the outward collar
shrinks to zero, so the limiting coordinates satisfy
$\xi_b=\xi_w$ and $\hat\xi(\xi_b,\eta_b)=\pm t_b$.
The limiting contours are tangent at
the black double saddle. The white saddle is different, because the
two zipper images are separated, so the strict gap holds at this
tangency point. On the black circle,
\[
 \Re(F_w-F_b)
 =\frac{\eta_w-\eta_b}{2}\log r_b
  +\frac{\xi_b-\xi_w}{2}\log|w(z)|.
\]
On the selected circle the last term is at most
$\tfrac12(\xi_b-\xi_w)\log R_w$ in both cases: for $\xi_b>\xi_w$ the
circle lies in the complement of $W$, where $|w(z)|\le R_w$, and for
$\xi_b\le\xi_w$ it lies in $W$, where $|w(z)|>R_w$.
The radius $R_w$ is bounded above and below, and $|w|$ is bounded
above and below at the black double saddle, which at a regular arctic
point is different from $-1$ and $\bk^{-1}$. Since
$\xi_b-\xi_w\to0$, the supremum of $\Re(F_w-F_b)$ over the selected
circle exceeds its value at the tangency point by $o(1)$.
The same strict gap consequently
holds on the entire selected circle for all sufficiently large $N$.
A compactness argument makes this bound uniform over all retained
pairs. Thus these full-circle residues are exponentially small too.

The remote double-integral pieces have the same strict gap; the
locally integrable Cauchy singularity has a uniformly bounded integral.
All these contributions are exponentially small.

The edge weights and the displacement factor in~\eqref{eq:adjacent-action}
are bounded. This proves~\eqref{eq:factorized-edge-bound}.
\end{proof}

\begin{proposition}[Arctic collar estimate]\label{prop:edge-collar}
Fix $N_0$ large enough that the truncated zippers are defined for every
$N\ge N_0$. For every compact $\mathcal C\subset X_{\rm unip}$,
\begin{equation}\label{eq:edge-hs-bounded}
 \sup_{N\ge N_0}\sup_{\rho\in\mathcal C}
       \|\widetilde O_N^{\rm edge}(\rho)\|_2<\infty,
\end{equation}
and
\begin{equation}\label{eq:edge-hs-tail}
 \limsup_{N\to\infty}\sup_{\rho\in\mathcal C}
 \|\widetilde O_N^{\rm edge}(\rho)
           -\widetilde O_{N,\varepsilon}(\rho)\|_2^2
 \le C_{\mathcal C}\sqrt\varepsilon.
\end{equation}
Here $\widetilde O_{N,\varepsilon}$ uses the same gauge and is extended
by zero on the omitted edges.
\end{proposition}

\begin{proof}
Transversality of the terminal segments gives at most
$C(N\Delta+1)$ edges whose transverse distance lies in any interval
of length $\Delta$. Hence
\begin{align*}
 \sum_{-2\delta_N\le d_e\le 2N^{-2/3}}a_N(e)^2
   &\le C N^{1/3+1/100}N^{-2/3}=o(1),\\
 \sum_{\substack{e\text{ liquid}\\2N^{-2/3}<d_e<\varepsilon}}a_N(e)^2
   &\le C\int_0^{C\varepsilon}s^{-1/2}\,ds+o(1)
    \le C\sqrt\varepsilon+o(1),\\
 \sum_{e\in E_N^{\rm edge}}a_N(e)^2&\le C.
\end{align*}
Square~\eqref{eq:factorized-edge-bound} and sum. Pairs with at least
one omitted endpoint are bounded by twice the product of the tail
sum and the full sum, giving~\eqref{eq:edge-hs-tail}.
\end{proof}

\subsection{Deleting the frozen suffixes}
\label{subsec:frozen-suffixes}

For $s\in[0,1]$ write
\begin{equation}\label{eq:section-arctic-endpoints}
 y_\pm(s)=\frac{(1-s)+\bk s\pm2\sqrt{\bk s(1-s)}}{1+\bk},
 \qquad
 s_{\rm top}=\frac{\bk}{1+\bk},\qquad
 s_{\rm bot}=\frac1{1+\bk}.
\end{equation}
By~\eqref{eq:liquid}, the section $\{s\}\times[0,1]$ meets $\Omega$ in
the interval $(y_-(s),y_+(s))$. Since
$(1+\bk)y_\pm(s)=(\sqrt{1-s}\pm\sqrt{\bk s})^2$, the function
$y_+$ is concave, increases from $1/(1+\bk)$ to $1$ on
$[0,s_{\rm top}]$ and decreases to $\bk/(1+\bk)$ on
$[s_{\rm top},1]$, while $y_-$ is convex, decreases from
$1/(1+\bk)$ to $0$ on $[0,s_{\rm bot}]$ and increases to
$\bk/(1+\bk)$ on $[s_{\rm bot},1]$. Moreover
$y_\pm(1-s)=1-y_\mp(s)$.
For $\bk=1$ we have $s_{\rm top}=s_{\rm bot}=1/2$ and
$y_\pm(s)=1/2\pm\sqrt{s(1-s)}$.
When the dependence on the bias matters, we write
$\psi_j^{(\bk)}$, $y_\pm^{(\bk)}$, $s_{\rm top}^{(\bk)}$ and
$s_{\rm bot}^{(\bk)}$.

\begin{lemma}\label{lem:forbidden-krawtchouk}
Let $A_1\Subset(0,s_{\rm top})$ and $A_2\Subset(0,s_{\rm bot})$.
There are $c,d_0>0$ such that, for $0<d\le d_0$, the following four
bounds hold.
\begin{enumerate}[label=\textnormal{(\alph*)},leftmargin=*]
\item If $a\in A_1$ and $x/N\ge y_+(a)+d$, then
\begin{equation}\label{eq:forbidden-krawtchouk}
 \sum_{j=0}^{\lfloor Na\rfloor}\psi_j(x)^2
 \le (N+1)^2e^{-cNd^{3/2}}.
\end{equation}
\item If $a\in A_2$ and $x/N\le y_-(a)-d$,
then~\eqref{eq:forbidden-krawtchouk} holds.
\item If $a\in A_2$ and $x/N\ge y_+(1-a)+d$, then
\begin{equation}\label{eq:forbidden-krawtchouk-high}
 \sum_{j=N-\lfloor Na\rfloor}^{N}\psi_j(x)^2
 \le (N+1)^2e^{-cNd^{3/2}}.
\end{equation}
\item If $a\in A_1$ and $x/N\le y_-(1-a)-d$,
then~\eqref{eq:forbidden-krawtchouk-high} holds.
\end{enumerate}
\end{lemma}

For $\bk=1$, parts \textup{(a)} and \textup{(b)} are the two
symmetric bounds beyond $1/2\pm\sqrt{a(1-a)}$ for the low modes, and
\textup{(c)}, \textup{(d)} are the same bounds for the high modes.

\begin{proof}
\emph{Part \textup{(a)}.}
Put $\alpha=j/N$ and $\beta=x/N$. For
$0<\alpha<s_{\rm top}$ and $\beta>y_+(\alpha)$, set
\[
 v=(1+\bk)\beta-1-(\bk-1)\alpha,\qquad
 \rho=\frac{v-\sqrt{v^2-4\bk\alpha(1-\alpha)}}{2(1-\alpha)},
 \qquad
 \rho_0(\alpha)=\sqrt{\frac{\bk\alpha}{1-\alpha}}.
\]
The condition $\beta>y_+(\alpha)$ is equivalent to
$v>2\sqrt{\bk\alpha(1-\alpha)}$, and $\alpha<s_{\rm top}$ is
equivalent to $\rho_0(\alpha)<\bk$. Then
$0<\rho<\rho_0(\alpha)<\bk$ and
$v=(1-\alpha)\rho+\bk\alpha/\rho$, that is,
\begin{equation}\label{eq:forbidden-beta}
 \beta=\frac{(1+\rho)\bigl(\rho+\alpha(\bk-\rho)\bigr)}
            {(1+\bk)\rho}.
\end{equation}
By~\eqref{eq:krawtchouk}, the logarithm of the
generating-function modulus on $|z|=\rho$, divided by $N$, is
\[
 g(u)=\frac{1-\beta}{2}\log(\bk^2+\rho^2+2\bk\rho u)
      +\frac\beta2\log(1+\rho^2-2\rho u),\qquad u=\cos\theta.
\]
Its derivative has the sign of
$(1-\beta)\bk(1+\rho^2-2\rho u)-\beta(\bk^2+\rho^2+2\bk\rho u)$,
which is decreasing in $u$ and equals
$(1-\beta)\bk(1+\rho)^2-\beta(\bk-\rho)^2$ at $u=-1$.
Substituting~\eqref{eq:forbidden-beta} gives
\begin{equation}\label{eq:forbidden-sign}
 \beta(\bk-\rho)^2-(1-\beta)\bk(1+\rho)^2
 =\frac{(1+\rho)(\bk-\rho)\bigl(\bk\alpha-(1-\alpha)\rho^2\bigr)}
       {\rho}\ge0.
\end{equation}
Thus $g'\le0$ on $[-1,1]$, and the maximum of $g$ occurs at
$z=-\rho$. Cauchy's coefficient formula gives
\[
 |K_j(x)|\le(\bk-\rho)^{N-x}(1+\rho)^x\rho^{-j}.
\]
Let $\mathsf H(u)=-u\log u-(1-u)\log(1-u)$ and put
\[
 \begin{aligned}
 h(\alpha,\beta)&=(1-\beta)\log(\bk-\rho)
             +\beta\log(1+\rho)-\alpha\log\rho,\\
 I(\alpha,\beta)&=\log(1+\bk)+(1-\alpha-\beta)\log\bk
       +\mathsf H(\alpha)-\mathsf H(\beta)-2h(\alpha,\beta).
 \end{aligned}
\]
The elementary binomial bounds
\[
 (N+1)^{-1}e^{N\mathsf H(j/N)}\le\binom Nj\le e^{N\mathsf H(j/N)}
\]
and~\eqref{eq:psi} imply $\psi_j(x)^2\le(N+1)e^{-NI(\alpha,\beta)}$.
Since $\rho$ is a stationary point of $h$ as a function of
$\rho$, differentiation gives
\[
 \partial_\alpha I
 =\log\frac{(1-\alpha)\rho^2}{\bk\alpha}\le0,\qquad
 \partial_\beta I
 =\log\frac{\beta}{\bk(1-\beta)}
       -2\log\frac{1+\rho}{\bk-\rho}\ge0,
\]
where the second inequality is~\eqref{eq:forbidden-sign}.
Since $y_+$ is increasing on $(0,s_{\rm top})$, it follows that
$I(\alpha,\beta)\ge I(a,\beta)$ when
$\alpha\le a<s_{\rm top}$ and $\beta>y_+(a)$.
For $j=0$, the same conclusion follows by continuity at $\alpha=0$,
or directly from $K_0(x)=\bk^{N-x}$ and
$I(0,\beta)=\log(1+\bk)-(1-\beta)\log\bk-\mathsf H(\beta)$.

At $\beta=y_+(a)$ one has $\rho=\rho_0(a)$,
\[
 (1+\bk)\,y_+(a)=(1-a)(1+\rho_0)^2,
 \qquad
 \bk(1+\bk)\bigl(1-y_+(a)\bigr)=(1-a)(\bk-\rho_0)^2,
\]
and substituting these two identities into $h$ gives
$I(a,y_+(a))=0$. If $d=\beta-y_+(a)$ and
$v_*=2\sqrt{\bk a(1-a)}$, then $v=v_*+(1+\bk)d$ and
\[
 \rho_0-\rho
 =\frac{\sqrt{(1+\bk)d\,\bigl(2v_*+(1+\bk)d\bigr)}-(1+\bk)d}
       {2(1-a)}\ge c_1\sqrt d
 \qquad(0<d\le d_0),
\]
uniformly for $a\in A_1$. Moreover, the two identities above give
\[
 \log\frac{y_+(a)}{\bk(1-y_+(a))}
 =2\log\frac{1+\rho_0}{\bk-\rho_0}.
\]
The function $\beta\mapsto\log[\beta/(\bk(1-\beta))]$ is increasing,
and the derivative of $\rho\mapsto\log[(1+\rho)/(\bk-\rho)]$ is
$(1+\bk)/[(1+\rho)(\bk-\rho)]\ge1/\bk$ on $(0,\bk)$.
It follows that
\[
 \partial_\beta I(a,y_+(a)+d)
 \ge 2\log\frac{1+\rho_0}{\bk-\rho_0}
      -2\log\frac{1+\rho}{\bk-\rho}
 \ge\frac2\bk(\rho_0-\rho)\ge\frac{2c_1}{\bk}\sqrt d.
\]
Integrating proves $I(a,y_+(a)+d)\ge c_2d^{3/2}$.
The function $I(a,\beta)$ increases beyond this point, so the
same lower bound holds for all larger $\beta$. Summing over at
most $N+1$ indices proves \textup{(a)}.

\emph{Parts \textup{(b)}--\textup{(d)}.}
For $0\le j,x\le N$,
\begin{equation}\label{eq:krawtchouk-dualities}
 \psi_j^{(\bk)}(N-x)^2=\psi_j^{(1/\bk)}(x)^2,
 \qquad
 \psi_{N-j}^{(\bk)}(x)^2=\psi_j^{(1/\bk)}(x)^2,
 \qquad
 \psi_{N-j}^{(\bk)}(N-x)^2=\psi_j^{(\bk)}(x)^2.
\end{equation}
Indeed, writing $K_j^{(\bk)}$ for the
polynomials~\eqref{eq:krawtchouk}, the substitution $z=-\bk w$ in
the generating function gives
$K_j^{(\bk)}(N-x)=(-\bk)^{-j}\bk^NK_j^{(1/\bk)}(x)$,
and the substitution $z\mapsto1/z$ gives
$K_{N-j}^{(\bk)}(x)=(-1)^x\bk^{N-x}K_j^{(1/\bk)}(x)$.
Inserting these relations into~\eqref{eq:psi} proves the first two
identities, and the third is their composition. For $\bk=1$ they
reduce to $\psi_j(N-x)^2=\psi_j(x)^2=\psi_{N-j}(x)^2$.
We also note that
\[
 y_+^{(1/\bk)}(a)=1-y_-^{(\bk)}(a)=y_+^{(\bk)}(1-a),
 \qquad
 s_{\rm top}^{(1/\bk)}=s_{\rm bot}^{(\bk)}.
\]
For \textup{(b)}, the first identity
in~\eqref{eq:krawtchouk-dualities} gives
$\sum_{j\le Na}\psi_j^{(\bk)}(x)^2
=\sum_{j\le Na}\psi_j^{(1/\bk)}(N-x)^2$, and
$(N-x)/N\ge1-y_-^{(\bk)}(a)+d=y_+^{(1/\bk)}(a)+d$; apply
\textup{(a)} with bias $1/\bk$ and the set
$A_2\Subset(0,s_{\rm top}^{(1/\bk)})$.
For \textup{(c)}, the second identity gives
$\sum_{j\ge N-\lfloor Na\rfloor}\psi_j^{(\bk)}(x)^2
=\sum_{i\le Na}\psi_i^{(1/\bk)}(x)^2$, and
$x/N\ge y_+^{(\bk)}(1-a)+d=y_+^{(1/\bk)}(a)+d$; apply
\textup{(a)} with bias $1/\bk$ again.
For \textup{(d)}, the third identity gives
$\sum_{j\ge N-\lfloor Na\rfloor}\psi_j^{(\bk)}(x)^2
=\sum_{i\le Na}\psi_i^{(\bk)}(N-x)^2$, and
$(N-x)/N\ge1-y_-(1-a)+d=y_+(a)+d$; apply \textup{(a)}.
Decreasing $c$ and $d_0$ makes the constants common to the four parts.
\end{proof}

Let $X_N(\rho)$ be the random product in
\eqref{eq:topological-correlator}, and let $X_N^{\rm edge}(\rho)$ be
the same product with transport only across $E_N^{\rm edge}$.
The latter connection has $2m$ punctures: $z_{i,N}^0$ and the
new endpoints $r_{i,N}$. Around a new endpoint its monodromy is conjugate
to $M_i^{-1}$, and is therefore still unipotent. In particular, the
small loops around these new endpoints have factor $1$.
Disjointness and the fixed isotopy of the cuts put these representations
in a compact family whenever $\rho$ ranges over a fixed compact set.
Proposition~\ref{prop:topological-ui} applies to both observables.

Let $D_N$ be the event that a nondoubled edge crosses a deleted suffix.
Then
\begin{equation}\label{eq:tail-agreement}
 X_N(\rho)=X_N^{\rm edge}(\rho)\quad\text{on }D_N^c,
 \qquad\text{simultaneously for all }\rho.
\end{equation}
We estimate $\PP(D_N)$ by the probability of an exceptional particle
or hole on a section meeting a deleted suffix.

The section-state law follows from~\eqref{eq:cap-state}: since
$\PP\{S_r=S\}=L_r(S)R_r(S)/Z_N$ and the gauge factors cancel,
\[
 \PP\{S_r=S\}
 =\left|\det[\psi_j(x)]_{x\in S,\ 0\le j<r}\right|^2.
\]
Cauchy--Binet and orthonormality therefore give
\begin{equation}\label{eq:section-occupation-projection}
 \Pi_r(x,x):=\PP\{x\in S_r\}
 =\sum_{j=0}^{r-1}\psi_j(x)^2.
\end{equation}
By~\eqref{eq:liquid}, the transverse arctic endpoints on section $r$
are $y_-(r/N)$ and $y_+(r/N)$. Along a white section,
$\eta=1-2y+O(N^{-1})$, so transverse $y$-distance and
$\eta$-distance differ only by the factor $2$ and a lattice error.
At the upper endpoint, for $r/N<s_{\rm top}$ apply
Lemma~\ref{lem:forbidden-krawtchouk}\textup{(a)}
with $a=(r-1)/N$; for $r/N>s_{\rm top}$ use
\[
 1-\Pi_r(x,x)=\sum_{j=r}^N\psi_j(x)^2
\]
and Lemma~\ref{lem:forbidden-krawtchouk}\textup{(c)} with
$a=(N-r)/N$. At the lower endpoint, for $r/N<s_{\rm bot}$ apply
Lemma~\ref{lem:forbidden-krawtchouk}\textup{(b)} with $a=(r-1)/N$,
and for $r/N>s_{\rm bot}$ apply
Lemma~\ref{lem:forbidden-krawtchouk}\textup{(d)} with $a=(N-r)/N$.
Accordingly, the frozen state beyond the upper endpoint is empty for
$r/N<s_{\rm top}$ and full for $r/N>s_{\rm top}$, and the frozen
state beyond the lower endpoint is empty for $r/N<s_{\rm bot}$ and
full for $r/N>s_{\rm bot}$.
Thus, for $\delta_N\le d\le d_0$, the wrong frozen occupation
at a site of outward distance at least $d$ has probability at most
\[
 C(N+1)^2 e^{-cNd^{3/2}}.
\]
The $O(N^{-1})$ changes of arctic location between adjacent sections
are absorbed by replacing $d$ by $d/2$. The points of $\partial\Omega$
with $s\in\{0,1\}$, the point $(s_{\rm top},1)$ and the point
$(s_{\rm bot},0)$ are the four tangency
points~\eqref{eq:tangency-points}. Regularity of the terminal
crossings therefore keeps $r/N$ in a compact subset of
$(0,1)\setminus\{s_{\rm top}\}$ at a crossing of the upper arc, and
in a compact subset of $(0,1)\setminus\{s_{\rm bot}\}$ at a crossing
of the lower arc.

The layer formula also shows that frozen section states force the
original matching edges. In layer $r$, put $S=S_r$, $T=S_{r+1}$.
If both states are empty on $\{L,\ldots,N\}$, the selected-white
prefix count is $c_i=i$ for $i\ge L$, forcing
$b_{r,i}$ to match $w_{r+1,i+1}$.
If both are full there, it instead matches $w_{r,i+1}$.
On a lower suffix $\{0,\ldots,L\}$, the prefix count is
$c_i=i+1$; empty or full states force respectively
$w_{r+1,i}$ or $w_{r,i}$. These follow from the counting argument
in Lemma~\ref{lem:layer-matching}.
Consequently an edge on a frozen terminal suffix can be nondoubled
only if one of the two perfect matchings has an exceptional particle or hole
in a neighboring section suffix.

Fix a sufficiently small regular arctic neighborhood, of transverse
width $d_0>0$, at each terminal crossing. Sum the preceding bound over
section suffixes, the two neighboring sections, the two perfect matchings,
and the $O(N)$ crossed edges. The defect probability in the portion
between outward distances $\delta_N$ and $d_0$ is at most
\[
 C(N+1)^4e^{-cN\delta_N^{3/2}}
 =C(N+1)^4e^{-cN^{3/200}}\longrightarrow0.
\]

Choose the remaining suffixes to lie at distance at least $d_0/2$
outside the arctic boundary, within the corresponding frozen regions,
all the way to the actual graph boundary. Each of them runs along
its terminal section, so the same parameter $a$ serves all its sites,
and the bounds of Lemma~\ref{lem:forbidden-krawtchouk} hold at every
site beyond the distance $d$, up to $x=0$ or $x=N$. Applying the
lemma with the fixed distance $d=d_0/4$ to the $O(N)$ sites of these
suffixes and of the neighboring section suffixes, and to both
independent perfect matchings, shows that these entire
suffixes are in the prescribed frozen tiling with probability at
least $1-C(N+1)^3e^{-cN}$. Combining the two suffix estimates gives
\begin{equation}\label{eq:frozen-defect}
 \PP(D_N)\longrightarrow0.
\end{equation}

\begin{proposition}[Frozen-tail deletion]\label{prop:frozen-deletion}
For every compact $\mathcal C\subset X_{\rm unip}$,
\begin{equation}\label{eq:frozen-deletion}
 \sup_{\rho\in\mathcal C}
 |F_N(\rho)-\det{}_2(I+O_N^{\rm edge}(\rho))|
 \longrightarrow0.
\end{equation}
\end{proposition}

\begin{proof}
The determinant identity and nilpotent elimination apply to the
truncated connection exactly as to the full connection, so
$\EE X_N^{\rm edge}=\det{}_2(I+O_N^{\rm edge})$.
By~\eqref{eq:tail-agreement}, Cauchy--Schwarz, and the uniform $L^2$
bound of Proposition~\ref{prop:topological-ui},
\[
 \sup_{\rho\in\mathcal C}
 |\EE X_N(\rho)-\EE X_N^{\rm edge}(\rho)|
 \le C_{\mathcal C}\PP(D_N)^{1/2}\longrightarrow0.
\]
\end{proof}

\subsection{Boundary reduction for the correlators}
\label{subsec:boundary-reduction}

\begin{theorem}\label{thm:boundary-reduction}
For every compact $\mathcal C\subset X_{\rm unip}$,
\begin{equation}\label{eq:determinant-boundary-reduction}
 \limsup_{N\to\infty}\sup_{\rho\in\mathcal C}
 |F_N(\rho)-\det{}_2(I+O_{N,\varepsilon}(\rho))|
 \le C_{\mathcal C}\varepsilon^{1/4}.
\end{equation}
\end{theorem}

\begin{proof}
Finite-dimensional similarity preserves $\det{}_2$.
For Hilbert--Schmidt operators, its continuity estimate has the form
\[
 |\det{}_2(I+A)-\det{}_2(I+B)|
 \le\|A-B\|_2
    \exp\{C(1+\|A\|_2+\|B\|_2)^2\}.
\]
Apply this with the common gauge from~\eqref{eq:edge-gauge}, use
Proposition~\ref{prop:edge-collar}, and then use
Proposition~\ref{prop:frozen-deletion}.
\end{proof}

\subsection{The canonical determinant and finite-puncture identification}
\label{subsec:canonical-determinant}

Let $F_{\DD}(\rho)$ denote the topological correlator of nested
CLE$_4$ in $\DD$ with punctures $\zeta_1,\ldots,\zeta_m$ and normalized
loop factors $\Tr\rho(\gamma)/2$.

Let $\gamma_i^0$ be the complete spectral cut from the $i$th puncture
to its real terminal point, with the convention of
\eqref{eq:canonical-cutoff-kernel}, and put
$C_i^0=\gamma_i^0-\overline{\gamma_i^0}$.
Define $\mathcal O_0$ by~\eqref{eq:canonical-cutoff-kernel} on these
complete contours and write $D_0=\det{}_2(I+\mathcal O_0)$.
The distinct real endpoints and disjointness of the cuts imply
\[
 \inf\{|q-r|:q\in C_i^0,\ r\in C_j^0,\ i\ne j\}>0.
\]
On a terminal segment let $d$ be the inward vertical distance to its
crossing point. By transversality, $1-\hat\xi^2-\eta^2$ has a simple
zero at $d=0$, and $1-\eta$ stays positive. Differentiating
\eqref{eq:spectral-coordinate} along the segment gives $|q'(d)|=O(d^{-1/2})$. Thus
the contours have finite length and their omitted arcs have total
length $O(\sqrt\varepsilon)$.
Extend $\mathcal O_\varepsilon$ by zero to the full contour space.
For every compact $\mathcal C\subset X_{\rm unip}$,
\begin{equation}\label{eq:canonical-collar-limit}
 \sup_{\rho\in\mathcal C}
 \|\mathcal O_\varepsilon(\rho)-\mathcal O_0(\rho)\|_2^2
 \le C_{\mathcal C}\sqrt\varepsilon,
 \qquad
 D_\varepsilon\longrightarrow D_0
 \quad\hbox{locally uniformly on }X_{\rm unip}.
\end{equation}
Here the cross-cut kernel is bounded uniformly on $\mathcal C$.

\begin{lemma}
\label{lem:canonical-reference-identification}
With the same cut generators and transport signs,
\begin{equation}\label{eq:canonical-reference-identity}
 D_0(\rho)=F_{\DD}(\rho),\qquad \rho\in X_{\rm unip}.
\end{equation}
\end{lemma}

\begin{proof}
Consider unbiased square-lattice double dimers in $\mathbb H$ at
mesh $\delta$, with representation $\rho_q$ and zippers approximating
the spectral cuts up to their distinct real endpoints. We show that
the correlators of this reference model converge to $D_0$,
identify $D_0$ with the $\CLE_4$ correlator of the half-plane, first
for real representations close to the trivial one and then by analytic
continuation, and then pass from the half-plane to the disk.

\emph{Gauge and orientations.}
On the reference graph put
$\vartheta_\delta(v)=\ii^{\Re v/\delta}$.
The reference~\cite{BI} uses $K_\delta^{\rm BI}(w,b)=b-w$, whereas we
use the dimensionless signs of~\eqref{eq:kasteleyn-definition} with
$\bk=1$ in the reference lattice coordinates. On the same reference
domain,
\begin{equation}\label{eq:BI-phase-conversion}
 \begin{aligned}
 K_\delta^{\rm ref}(w,b)
 &=\frac{\ii}{\delta}\,
   \vartheta_\delta(w)\vartheta_\delta(b)^{-1}
   K_\delta^{\rm BI}(w,b),\\
 (K_\delta^{\rm ref})^{-1}(b,w)
 &=-\ii\delta\,
   \vartheta_\delta(b)\vartheta_\delta(w)^{-1}
   (K_\delta^{\rm BI})^{-1}(b,w).
 \end{aligned}
\end{equation}
The factors $\delta$ cancel in the edge kernel. Writing
$(Q_E)_{ee}=\vartheta_\delta(b_e)I_2$, we obtain
$B_\delta^{\rm ref}=Q_EB_\delta^{\rm BI}Q_E^{-1}$.
Thus, by~\eqref{eq:Kasteleyn-gauge-covariance}, the determinant and the
cyclic traces are the same in the two conventions. The determinant
identity for the reference model is \cite[Lemma 2.10]{BI}; the
definition by finite-volume exhaustion is also given
in~\cite[Section 5.4]{Dubedat}. The exact elimination of
Proposition~\ref{prop:nilpotent-reduction} applies to the reference
model.

\emph{The reference kernel.}
Write $\chi_v=(-1)^{\Im v/\delta}$ and
$b_+=b,\ b_-=\bar b,\ w_+=w,\ w_-=\bar w$.
On distinct cuts, \cite[equations~(2.8)--(2.10)]{BI} give
\[
 (K_\delta^{\rm BI})^{-1}(b,w)
 =\frac1{2\pi}(1,-\chi_b)
   \left[\frac1{b_s-w_t}\right]_{s,t\in\{+,-\}}
   \binom{1}{-\chi_w}+O(\delta).
\]
The error is uniform up to the distinct real endpoints. For the
oriented dual edge increment $\Delta r$, the signed weight is
$\ii\Delta r$ and its product with $\chi_b\chi_w$ is
$-\ii\Delta\bar r$ by \cite[equations~(2.14)--(2.15)]{BI}.

Choose the reference zippers with bounded total length and
$O(\delta^{-1/2})$ straight runs, by approximating successive arcs of
length at most $\sqrt\delta$ with coordinate segments and rounding to
the mesh within disjoint tubes. A mismatch of the two branch labels at
an edge in a cyclic product leaves a single factor $\chi_b$ or $\chi_w$.
It alternates on each straight run. The separated Cauchy denominators
are bounded with bounded first derivatives, so pairing successive
edges costs $O(\delta^2)$ per pair after including the edge increment.
The unpaired ends cost $O(\delta)$ per run. Summation in this edge
therefore gives $O(\sqrt\delta)$, uniformly in the other positions.
Only branches with the same label survive in the limit.
With the transported monodromy $I_2-A_i$, the kernel
differentials along $\gamma_i^0$ and its reflection are
\[
 \frac{(-A_i)\ii\,dr}{2\pi(q-r)}
 =\frac{A_i\,dr}{2\pi\ii(q-r)},
 \qquad
 \frac{(-A_i)(-\ii\,d\bar r)}{2\pi(q-\bar r)}
 =-\frac{A_i\,d\bar r}{2\pi\ii(q-\bar r)}.
\]
These are the column factors in~\eqref{eq:canonical-cutoff-kernel}
on the chain $C_i^0=\gamma_i^0-\overline{\gamma_i^0}$. The Riemann
sums on the approximating paths converge to the corresponding cyclic
integrals. Indeed, the paths converge uniformly with bounded variation
and the Cauchy kernels are smooth on the separated cuts. The finitely
many endpoint cells contribute $O(\delta)$. Summing over the label
sequences gives $\Tr\mathcal O_0^k$. The cross entries are
$O(\delta)$ and each cut has $O(\delta^{-1})$ edges, so the
Hilbert--Schmidt norms are bounded uniformly in $\delta$. The
auxiliary-$t$ argument in the proof of
Proposition~\ref{prop:interior-determinant} then shows that the
reference correlators converge to $D_0$, locally uniformly on
$X_{\rm unip}$.

\emph{Identification with the half-plane correlator.}
Let $F_{\mathbb H}(\rho_q)$ be the topological correlator of nested
$\CLE_4$ in $\mathbb H$ at the punctures $q(z_1),\ldots,q(z_m)$, with
loop factors $\Tr\rho_q(\ell)/2$. For real $A_i$ the representation
$\rho_q$ takes values in $\mathrm{SL}_2(\mathbb R)$, with unipotent
local monodromies $I_2-A_i$. By~\cite[Theorem~1]{Dubedat}, which
concerns this half-plane model, for real $A_i$ close enough to $0$ the
reference correlators converge to $F_{\mathbb H}(\rho_q)$; hence, by
the first part, $D_0(\rho)=F_{\mathbb H}(\rho_q)$ for these $\rho$.
Both sides extend to entire functions of the matrices
$(A_1,\ldots,A_m)$: $D_0$ because $\mathcal O_0$ is linear in the
$A_i$, and $F_{\mathbb H}$ because its expansion over laminations
converges uniformly on bounded sets. Indeed, the loop products and the
number of laminations grow at most exponentially in the lamination
complexity, while the lamination probabilities decay faster than every
exponential by the estimate of~\cite{BW}, recorded
in~\cite[Theorem~4.2(1)]{BI}. Every nilpotent $2\times2$ matrix has
the form
$\bigl(\begin{smallmatrix}uv&-u^2\\ v^2&-uv\end{smallmatrix}\bigr)$,
which is real for real $(u,v)$. As entire functions of
$(u_i,v_i)_{i\le m}\in\mathbb C^{2m}$, the two sides agree near the
origin of $\mathbb R^{2m}$, hence everywhere. Thus
$D_0(\rho)=F_{\mathbb H}(\rho_q)$ for every $\rho\in X_{\rm unip}$.

\emph{From the half-plane to the disk.}
By~\eqref{eq:spectral-coordinate}, $q\circ\Psi$ is an
anticonformal bijection from $\DD$ onto $\mathbb H$.
Conformal and reflection invariance of $\CLE_4$, with the transport
of $\rho$ fixed above, identify $F_{\mathbb H}(\rho_q)$ with
$F_{\DD}(\rho)$.
\end{proof}

\begin{theorem}[Finite-puncture identification]
\label{thm:aztec-finite-puncture}
Fix distinct points $z_1,\ldots,z_m$ in $L$ and put $\zeta_i=\Phi(z_i)$.
With the representatives and $\Gamma_N$ of Section~\ref{sec:punctures},
\begin{equation}\label{eq:aztec-tau-limit}
 \sup_{\rho\in\mathcal C}|F_N(\rho)-F_{\DD}(\rho)|
       \longrightarrow0
 \qquad\hbox{for every compact }\mathcal C\subset X_{\rm unip}.
\end{equation}
Consequently,
\begin{equation}\label{eq:aztec-lamination-limit}
 \Gamma_N
 \ \Longrightarrow\
 \Lam_{\zeta_1,\ldots,\zeta_m}
       \bigl(\CLE_4^{\rm nested}(\DD)\bigr).
\end{equation}
The convergence of the reduced-lamination laws is in total variation.
Both conclusions remain valid for deterministic punctures
$z_{i,N}\to z_i$ lying off the scaled embedded graph.
Identify homotopy classes and representation bases by isotopies
supported in disjoint liquid neighborhoods of the distinct
limiting points $z_i$. The limiting punctures remain
$\zeta_i=\Phi(z_i)$.
\end{theorem}

\begin{proof}
Combine Theorem~\ref{thm:boundary-reduction},
Proposition~\ref{prop:interior-determinant},
\eqref{eq:canonical-collar-limit}, and
Lemma~\ref{lem:canonical-reference-identification}. First let
$N\to\infty$ with $\varepsilon$ fixed, then let $\varepsilon\downarrow0$.
This proves~\eqref{eq:aztec-tau-limit} uniformly on every compact
subset of the whole complex unipotent variety.

Let $\Gamma_{\rm CLE}$ be the reduced lamination of nested CLE$_4$
at $\zeta_1,\ldots,\zeta_m$. Fix a triangulation of the punctured disk
and let $|\Gamma|$ be the minimal intersection number of a lamination
with its edges. For a reduced lamination $\Gamma$ put
\[
 p_N(\Gamma)=2^{-\#\Gamma}\PP\{\Gamma_N=\Gamma\},\qquad
 p_{\DD}(\Gamma)=2^{-\#\Gamma}\PP\{\Gamma_{\rm CLE}=\Gamma\},
\]
and set $f_\Gamma(\rho)=\prod_{\gamma\in\Gamma}\Tr\rho(\gamma)$.
Then
\[
 F_N=\sum_\Gamma p_N(\Gamma)f_\Gamma,\qquad
 F_{\DD}=\sum_\Gamma p_{\DD}(\Gamma)f_\Gamma.
\]
The second series converges uniformly on compact subsets of
$X_{\rm unip}$ by the complexity bounds used in the proof of
Lemma~\ref{lem:canonical-reference-identification}; thus both
correlators are invariant entire functions.
The coefficient estimate of~\cite[Theorem 1.4 and Corollary 1.6]{BC}
applied to their difference $F_N-F_{\DD}$ gives, for
each $R>0$, a compact $K_R\subset X_{\rm unip}$ and a constant $C_R$
such that
\[
 |p_N(\Gamma)-p_{\DD}(\Gamma)|
 \le C_RR^{-|\Gamma|}
       \sup_{\rho\in K_R}|F_N(\rho)-F_{\DD}(\rho)|.
\]
Hence every reduced-lamination probability converges. The state
space is countable and the limiting probabilities sum to one;
pointwise convergence therefore implies total-variation convergence.
This proves~\eqref{eq:aztec-lamination-limit}.

For moving punctures, write $F_N^{\mathrm{mov}}$ for the
correlator and put $\Delta_N=\max_i|z_{i,N}-z_{i,N}^0|$.
Choose their zippers and the reference zippers from $z_{i,N}^0$
with identical tails.
Their initial portions lie in disjoint compact liquid
neighborhoods and can be chosen so that their signed
crossed-edge lists differ on at most $C(N\Delta_N+1)$ edges.
Use the same monodromy matrices in the identified bases.

For fixed $\varepsilon$, extend the two cutoff operators
by zero to the union of their edge sets and use the common
gauge of~\eqref{eq:edge-gauge}.
Their entries agree unless a row or column corresponds to
a changed edge. There are $O(N)$ indices, and different
zipper labels stay uniformly separated. The entry bound
in the proof of Proposition~\ref{prop:interior-determinant}
therefore gives
\[
 \sup_{\rho\in\mathcal C}
 \|\widetilde O_{N,\varepsilon}^{\mathrm{mov}}(\rho)
       -\widetilde O_{N,\varepsilon}(\rho)\|_2^2
 \le C_{\varepsilon,\mathcal C}(\Delta_N+N^{-1})
 \longrightarrow0.
\]
The two Hilbert--Schmidt norms are uniformly bounded,
so continuity of $\det_2$ gives the same cutoff limit.
The boundary reduction is uniform for these paths:
their terminal portions are unchanged, and
Proposition~\ref{prop:topological-ui} applies uniformly
to punctures in the chosen neighborhoods. Hence
\[
 \limsup_{N\to\infty}\sup_{\rho\in\mathcal C}
 |F_N^{\mathrm{mov}}(\rho)-F_N(\rho)|
 \le C_{\mathcal C}\varepsilon^{1/4}.
\]
Letting $\varepsilon\downarrow0$ proves the same correlator
limit. The preceding coefficient argument in the identified
bases then gives total-variation convergence for the
moving-puncture lamination laws.
\end{proof}

\subsection{Frozen regions outside the arctic ellipse}
\label{subsec:frozen-regions}

The global curve limit requires control throughout each frozen region.
For $\bk=1$, the Arctic Circle Theorem of Jockusch, Propp and
Shor~\cite[Theorem~1]{JPS} gives freezing with probability tending to
one at fixed positive distance from the liquid region. For biased tilings,
Cohn, Elkies and Propp~\cite[Theorem~26]{CEP} proved that, for each fixed
$\varepsilon>0$, the probability that the frozen boundary deviates from
the arctic ellipse by more than $\varepsilon$ is exponentially small in $N$.
Lemma~\ref{lem:forbidden-krawtchouk} gives the estimates below as well as
the bound at the shrinking distance $\delta_N$ used in
Proposition~\ref{prop:frozen-deletion}.

Recall $\Omega$ from~\eqref{eq:liquid} and $y_\pm$, $s_{\rm top}$,
$s_{\rm bot}$ from~\eqref{eq:section-arctic-endpoints}. The complement
of $\overline\Omega$ in $[0,1]^2$ is the union of the four regions
\begin{equation}\label{eq:frozen-regions}
 \begin{aligned}
 R_{\rm W}&=\{s<s_{\rm top},\ y>y_+(s)\},&
 R_{\rm N}&=\{s>s_{\rm top},\ y>y_+(s)\},\\
 R_{\rm S}&=\{s<s_{\rm bot},\ y<y_-(s)\},&
 R_{\rm E}&=\{s>s_{\rm bot},\ y<y_-(s)\},
 \end{aligned}
\end{equation}
which the map $(s,y)\mapsto(s-y)+\ii(s+y-1)$ sends to the west, north,
south and east corners of $\Diamond\setminus\overline L$.
The \emph{frozen partner} of a black vertex $b_{r,i}$ is
$w_{r+1,i+1}$ if $(r/N,i/N)\in R_{\rm W}$,
$w_{r,i+1}$ if $(r/N,i/N)\in R_{\rm N}$,
$w_{r+1,i}$ if $(r/N,i/N)\in R_{\rm S}$, and
$w_{r,i}$ if $(r/N,i/N)\in R_{\rm E}$.
Thus the frozen edges are vertical in the west and east regions and
horizontal in the north and south regions.

\begin{proposition}[Frozen regions]\label{prop:frozen-regions}
For every $\eta>0$ there are constants $C,c>0$ and $N_0$ such that,
for $N\ge N_0$, the following event has probability at least
$1-C(N+1)^4e^{-cN}$: every black vertex $b_{r,i}$ for which the point
$(r/N,i/N)$ has distance at least $\eta$ from $\Omega$ is matched to
its frozen partner in both perfect matchings $M_N^{(1)}$ and $M_N^{(2)}$.
On this event the edges at these black vertices are doubled, no
nontrivial loop of $\mathcal L_N$ visits them, and every loop of
$\mathcal L_N$ lies within distance $2\eta$ of $\overline L$.
\end{proposition}

\begin{proof}
We measure distances in the Euclidean $(s,y)$-plane.
The map to the $z$-plane in~\eqref{eq:seam-physical-coordinates}
multiplies distances by $\sqrt2$.

\emph{Choice of compact ranges.}
The functions $y_\pm$ are continuous, $y_+(s_{\rm top})=1$,
$y_-(s_{\rm bot})=0$, $y_\pm(0)=1/(1+\bk)$, and
$y_\pm(1)=\bk/(1+\bk)$. We can therefore choose
$0<s_1<s_2<s_{\rm top}<s_3<s_4<1$ with $s_1\le\eta/4$,
$1-s_4\le\eta/4$, and
\[
 1-y_+(s)<\frac\eta2\ \ (s_2\le s\le s_3),\qquad
 y_+(s_1)-y_+(0)\le\frac\eta8,\qquad
 y_+(s_4)-y_+(1)\le\frac\eta8,
\]
and $0<s_1'<s_2'<s_{\rm bot}<s_3'<s_4'<1$ with the corresponding
properties of $y_-$:
$y_-(s)<\eta/2$ on $[s_2',s_3']$,
$y_-(0)-y_-(s_1')\le\eta/8$, $y_-(1)-y_-(s_4')\le\eta/8$,
$s_1'\le\eta/4$ and $1-s_4'\le\eta/4$.

\emph{One-site bounds.}
Let $(r,x)$ be a site with $1\le r\le N$, put $(s,y)=(r/N,x/N)$, and
assume that $(s,y)$ has distance at least $\eta/2$ from $\Omega$.
We claim that, for $N\ge N_0$,
\begin{equation}\label{eq:frozen-one-site}
 \begin{aligned}
 \Pi_r(x,x)&\le(N+1)^2e^{-c_\eta N}
 &&\text{if }(s,y)\in R_{\rm W}\cup R_{\rm S},\\
 1-\Pi_r(x,x)&\le(N+1)^2e^{-c_\eta N}
 &&\text{if }(s,y)\in R_{\rm N}\cup R_{\rm E},
 \end{aligned}
\end{equation}
where $c_\eta=c\,d_\eta^{3/2}$, $d_\eta=\min\{\eta/8,d_0\}$, and
$c,d_0$ are the constants of
Lemma~\ref{lem:forbidden-krawtchouk} for the compact sets chosen
below, which depend only on $\eta$. In each application of the lemma
the hypothesis will hold with a distance at least $\eta/8$, hence
with $d=d_\eta$. We use that,
by~\eqref{eq:section-occupation-projection},
$\Pi_r(x,x)$ is nondecreasing in $r$ for fixed $x$.

Let $(s,y)\in R_{\rm W}$. The vertical distance $y-y_+(s)$ is at
least the distance to $\Omega$, so $y\ge y_+(s)+\eta/2$; in
particular $s<s_2$. If $r\ge r_1:=\lceil Ns_1\rceil$, apply
Lemma~\ref{lem:forbidden-krawtchouk}\textup{(a)} with
$A_1=[s_1/2,s_2]$ and $a=(r-1)/N$: since $y_+$ is
increasing on $[0,s_{\rm top}]$, we have
$x/N\ge y_+(a)+\eta/2$, and
$\Pi_r(x,x)=\sum_{j\le r-1}\psi_j(x)^2$ satisfies the claim.
If $r<r_1$, then $s<s_1\le\eta/4$, and the distance from $(s,y)$
to the point $(0,y_+(0))\in\overline\Omega$ is at least $\eta/2$;
hence $y-y_+(0)\ge\eta/4$ and $y\ge y_+(s_1)+\eta/8$.
By monotonicity and part \textup{(a)} with $a=(r_1-1)/N<s_1$,
$\Pi_r(x,x)\le\Pi_{r_1}(x,x)\le(N+1)^2e^{-c_\eta N}$.

Let $(s,y)\in R_{\rm N}$. Again $y\ge y_+(s)+\eta/2$ and $s>s_3$.
Put $r_4=\lceil Ns_4\rceil$. If $r\le r_4$, apply
Lemma~\ref{lem:forbidden-krawtchouk}\textup{(c)} with
$A_2=[(1-s_4)/2,1-s_3]\Subset(0,s_{\rm bot})$ and $a=(N-r)/N=1-s$:
the sum in~\eqref{eq:forbidden-krawtchouk-high} is
$\sum_{j\ge r}\psi_j(x)^2=1-\Pi_r(x,x)$.
If $r>r_4$, then $1-s<\eta/4$, the distance from $(s,y)$ to
$(1,y_+(1))\in\overline\Omega$ is at least $\eta/2$, and therefore
$y\ge y_+(1)+\eta/4\ge y_+(r_4/N)+\eta/8$, because $y_+$ is
decreasing on $[s_{\rm top},1]$. Monotonicity and part
\textup{(c)} with $a=(N-r_4)/N$ give
$1-\Pi_r(x,x)\le1-\Pi_{r_4}(x,x)\le(N+1)^2e^{-c_\eta N}$.

The regions $R_{\rm S}$ and $R_{\rm E}$ are treated in the same way,
with $y_-$, the numbers $s_1',\ldots,s_4'$, and parts
\textup{(b)} and \textup{(d)} of
Lemma~\ref{lem:forbidden-krawtchouk} in place of $y_+$,
$s_1,\ldots,s_4$, and parts \textup{(a)} and \textup{(c)}.
This proves~\eqref{eq:frozen-one-site}.

\emph{The frozen event.}
Let $G_N$ be the event that, in both perfect matchings, every site
$(r,x)$ at distance at least $\eta/2$ from $\Omega$ has the frozen
occupation: $x\notin S_r$ if $(r/N,x/N)\in R_{\rm W}\cup R_{\rm S}$,
and $x\in S_r$ if $(r/N,x/N)\in R_{\rm N}\cup R_{\rm E}$.
There are $N(N+1)$ sites and two perfect matchings,
so~\eqref{eq:frozen-one-site} gives
$\PP(G_N^c)\le2(N+1)^4e^{-c_\eta N}$.

\emph{Frozen states force frozen edges.}
Assume $G_N$, and let $b_{r,i}$ be a black vertex such that
$(r/N,i/N)$ has distance at least $\eta$ from $\Omega$, say in
$R_{\rm W}$. The set $\Omega$ is convex, so for fixed $s$ the
distance from $(s,y)$ to $\Omega$ is a convex function of $y$; it
vanishes at $y=y_+(s)$ and is therefore nondecreasing for
$y\ge y_+(s)$. Hence, for $N\ge N_0$, all sites $(r,x)$ and
$(r+1,x)$ with $x\ge i$ lie in $R_{\rm W}$ at distance at least
$\eta/2$ from $\Omega$; here we use that $r/N<s_2<s_{\rm top}$, so
that the two sections lie on the same side of $s_{\rm top}$.
On $G_N$, the states $S_r$ and $S_{r+1}$ are thus empty on
$\{i,\ldots,N\}$. By the counting argument of
Lemma~\ref{lem:layer-matching}, recalled in
Section~\ref{subsec:frozen-suffixes}, this forces $b_{r,i}$ to be
matched to $w_{r+1,i+1}$. In $R_{\rm N}$ both states are full on
$\{i,\ldots,N\}$ and $b_{r,i}$ is matched to $w_{r,i+1}$; in
$R_{\rm S}$ and $R_{\rm E}$ both states are empty, respectively
full, on $\{0,\ldots,i+1\}$, and $b_{r,i}$ is matched to
$w_{r+1,i}$, respectively $w_{r,i}$. For the boundary layers
$r=0$ and $r=N$ the same argument applies with the formal states
$S_0=\varnothing$ and $S_{N+1}=E_N$, which have the frozen
occupations of the regions adjacent to $s=0$ and $s=1$.
In every case the partner of $b_{r,i}$ is its frozen partner, in
both perfect matchings.

\emph{Loops.}
On $G_N$, the edge at each such black vertex is doubled, so no
nontrivial loop visits it. Every point of a nontrivial loop lies
on an edge of the loop, within physical distance $1/N$ of a black
vertex of the loop. That vertex is within distance $\eta$ of
$\Omega$ in the coordinates $(s,y)$, hence within physical
distance $\sqrt2\,\eta+O(N^{-1})$ of $L$. For $N\ge N_0$ this is less
than $2\eta$.
\end{proof}

The records of Theorem~\ref{thm:aztec-finite-puncture} retain loops
enclosing all the punctures. The next
section upgrades their convergence to convergence of the full curves,
using the opposite-crossing estimate of Theorem~\ref{thm:global-opposite}
and Proposition~\ref{prop:frozen-regions}.

\section{Proof of the global curve limit}
\label{sec:global-curves}

Write $\mathcal C=\Psi(\CLE_4^{\rm nested}(\DD))$. The traversal bound
gives tightness of the entire collection
(Sections~\ref{subsec:polish}--\ref{subsec:tightness}). Along a
subsequence we couple the loops with their enclosure indicators on a
countable dense set of probes (Sections~\ref{subsec:probes}--\ref{subsec:coupling}).
The finite-puncture limit identifies the sets of probes enclosed by
the loops of $\mathcal C$. We use the opposite-crossing estimate twice.
First it rules out excursions into a disk on which the enclosure
indicator is constant, which identifies each limiting trace
(Sections~\ref{subsec:signatures}--\ref{subsec:traces}). It then rules
out backtracking along that trace, which identifies the curve up to
reparametrization (Section~\ref{subsec:monotone}).

We use the diagonal coordinates $(s,y)$ of Section~\ref{subsec:seams};
their physical realization includes the $O(N^{-1})$ translation
in~\eqref{eq:seam-physical-coordinates}.

\subsection{The space of loop collections}
\label{subsec:polish}

\begin{lemma}\label{lem:loop-collection-polish}
The function $d_{\rm coll}$ is a metric, and
$(\mathcal X_\Diamond,d_{\rm coll})$ is complete and separable.
\end{lemma}

\begin{proof}
\emph{The metric.}
For individual loops,
\[
 |\diam\gamma-\diam\eta|
 \le2d_{\rm curv}(\gamma,\eta).
\]
Composing partial pairings and using this inequality proves the
triangle inequality for $d_{\rm coll}$. To check that zero distance
separates collections, fix a loop of diameter $\varepsilon>0$. Pairings
of error less than $\varepsilon/4$ must pair each copy of this loop with
a loop of diameter at least $\varepsilon/2$ in the other collection.
There are only finitely many such candidates. Letting the pairing
errors tend to zero forces the same loop to occur there with at least
the same multiplicity. Reversing the roles of the collections gives
equality of every multiplicity.

\emph{Single loops.}
The single-loop space with constant maps included is complete:
from a Cauchy sequence, take a subsequence with summable successive
distances and choose successive reparametrizations, allowing
summable additional errors in approximating the infima, so that
its representatives are uniformly Cauchy. Their uniform limit is a
continuous map into the closed set $\Diamond$.

\emph{Completeness.}
Take a $d_{\rm coll}$-Cauchy sequence and a subsequence with
partial pairings of successive errors $e_n$ satisfying
$\sum_n e_n<\infty$. These pairings form chains of loop copies.
Each infinite chain is Cauchy, so it has a limit in the single-loop
space. Retain the nonconstant limits, one copy per chain. At level
$n$, any terminating chain, or chain with
constant limit, has loop diameter at most
$2\sum_{k\ge n}e_k$; a retained chain born after level $n$ has
limiting diameter at most the same bound. Pairing the other
level-$n$ loops with their limits therefore gives collection
distance at most $\sum_{k\ge n}e_k$. For each positive diameter
cutoff, all retained loops above that cutoff correspond, at a
sufficiently late fixed level, to loops above half the cutoff.
There are only finitely many such loops. The resulting limit
belongs to $\mathcal X_\Diamond$, proving completeness.

\emph{Separability.}
Closed polygonal paths with finitely many rational vertices in
$\Diamond$ form a countable dense family for the single-loop
distance, since $\Diamond$ is convex. Finite collections drawn
from this family are dense in $\mathcal X_\Diamond$: discard all
loops below a small diameter cutoff and approximate the finitely
many remaining loops. This proves separability.
\end{proof}

\subsection{Tightness}
\label{subsec:tightness}

For a finite collection $\Gamma$ of polygonal loops and $\delta>0$,
let $n_\delta(\Gamma)$ be the largest number of loop subarcs whose
endpoints are at distance at least $\delta$ from each other and whose
parameter interiors are pairwise disjoint; subarcs of different loops
are always considered disjoint. This number is finite, since each such
subarc has length at least $\delta$.

\begin{lemma}[Few long subarcs]\label{lem:long-subarcs}
For every $0<\delta\le1$,
\[
 \lim_{M\to\infty}\sup_N
 \PP\{n_\delta(\mathcal L_N)>M\}=0 .
\]
\end{lemma}

\begin{proof}
If two points are at distance at least $\delta$, one of their
coordinates $s,y$ differs by at least $\delta/2$,
by~\eqref{eq:seam-physical-coordinates}. Choose a grid of mesh
$h=1/k<\delta/12$ in both coordinate directions. A coordinate interval
of length at least $\delta/2$ contains two consecutive grid levels in
$(0,1)$ at distance at least $h$ from its endpoints. Round these levels
to sections as in Definition~\ref{def:traversals}; for all sufficiently
large $N$ this moves them by less than $h/2$. In the second coordinate
direction use the reflected sections of
Section~\ref{subsec:reflection}.

Take a subarc whose endpoints are at distance at least $\delta$, and
the two sections just found. Starting from the first endpoint of the
subarc, take its first visit to the farther section and its last
earlier visit to the nearer section. The part of the subarc between
these two visits is a traversal of the slab between the two sections.
Subarcs with disjoint parameter interiors give different traversals.
Thus, for all sufficiently large $N$, $n_\delta(\mathcal L_N)$ is at
most the sum of the traversal counts of the finitely many slabs between
consecutive grid levels. Proposition~\ref{prop:traversal-moments}
bounds the exponential moments of this sum uniformly in $N$. The
finitely many remaining graphs have finitely many configurations.
\end{proof}

\begin{lemma}[Finite approximation]\label{lem:finite-approximation}
For every $\delta>0$ and every integer $M\ge1$ there is a finite set
$F_{\delta,M}\subset\mathcal X_\Diamond$ such that
\[
 n_\delta(\Gamma)\le M\quad\Longrightarrow\quad
 d_{\rm coll}(\Gamma,F_{\delta,M})\le3\delta
\]
for every finite collection $\Gamma$ of polygonal loops in $\Diamond$.
\end{lemma}

\begin{proof}
Suppose $n_\delta(\Gamma)\le M$. Leave the loops of diameter at most
$6\delta$ unpaired; this costs at most $3\delta$ in $d_{\rm coll}$.
Each other loop contains a subarc whose endpoints are at distance at
least $\delta$, so at most $M$ loops remain.

Take one of the remaining loops and a starting point on it. Mark the
first time the loop is at distance $\delta$ from the starting point,
then the first time it is at distance $\delta$ from this mark, and continue to the final time. Each completed step is a
subarc counted by $n_\delta$, so there are at most $M+1$ pieces,
including a possibly incomplete last one. Each piece stays in the
closed $\delta$-ball about its initial point.

Fix a finite $\delta$-net in $\Diamond$. Replace every marked position
by a nearest point of the net, chosen by a fixed rule, and join
consecutive points by straight segments. The two new endpoints of a
piece lie in the $2\delta$-ball about its initial point, hence so does
the segment joining them. So the polygon is within $3\delta$ of the
loop, pointwise in a common parametrization. It lies in $\Diamond$ by
convexity, and it is nonconstant because the loop has diameter larger
than $6\delta$. Up to reparametrization it is given by a word of at
most $M+2$ net points. Let $F_{\delta,M}$ be the finite set of all
collections of at most $M$ such nonconstant closed polygons. Pairing
each remaining loop with its polygon gives
$d_{\rm coll}(\Gamma,F_{\delta,M})\le3\delta$.
\end{proof}

\begin{proposition}[Tightness]\label{prop:tightness}
The laws of $\mathcal L_N$, $N\ge1$, are tight in
$(\mathcal X_\Diamond,d_{\rm coll})$.
\end{proposition}

\begin{proof}
Given $\varepsilon>0$, put $\delta_i=2^{-i}$ and use
Lemma~\ref{lem:long-subarcs} to choose $M_i$ with
$\sup_N\PP\{n_{\delta_i}(\mathcal L_N)>M_i\}\le\varepsilon2^{-i}$,
$i\ge1$. The set
\[
 K_\varepsilon=\bigcap_{i\ge1}
 \{\Gamma\in\mathcal X_\Diamond:
       d_{\rm coll}(\Gamma,F_{\delta_i,M_i})\le3\delta_i\}
\]
is closed and has finite covers by balls of arbitrarily small radius.
It is compact by Lemma~\ref{lem:loop-collection-polish}. By
Lemma~\ref{lem:finite-approximation} and the union bound,
$\inf_N\PP\{\mathcal L_N\in K_\varepsilon\}\ge1-\varepsilon$.
\end{proof}

\subsection{Test lines and probes}
\label{subsec:probes}

\paragraph{Test lines.}
Choose a countable dense set of levels in $(0,1)$ for each of the two
coordinates. The lines $s=c$ and $y=c$, with $c$ one of these levels,
are the \emph{test lines}. Choose the levels so that no test line
passes through a vertex of a scaled graph $N^{-1}G_N$, $N\ge1$. A test
line then meets $\partial L$ in at most two points and the union of all
scaled graphs in a countable set.

\paragraph{Partitions and probes.}
Fix a bounded interval $J_0$ that contains $[0,1]$ in its interior.
On each test line $H$ fix partitions
$\mathcal P_1(H),\mathcal P_2(H),\ldots$ of the segment of $H$ with
transverse coordinate in $J_0$ into finitely many intervals, called
\emph{cells}, of length at most $h_m$ for $\mathcal P_m(H)$, with
$h_m\downarrow0$. The two ends of $J_0$ are endpoints of cells in every
partition. The other endpoints of cells are called \emph{probes}. We
choose them off $\partial L$, off $\partial\Diamond$, off all scaled
graphs, and off all test lines of the other direction. Let
$\mathcal Z$ be the set of all probes, for all test lines and all $m$.
It is countable, $\mathcal Z\cap\Diamond$ is dense in $\Diamond$, and
every probe lies on exactly one test line.

Almost surely no probe lies on a loop of $\mathcal C$. For probes in
$L$ this holds because a fixed point almost surely lies on no loop of
$\CLE_4$; probes outside $\overline L$ lie on no loop because all
loops of $\mathcal C$ are contained in $L$.

\paragraph{Discrete realization.}
Let $H$ be the test line $s=c$; for a test line $y=c$ use the reflected
cut construction of Section~\ref{subsec:reflection}. For each $N$ let
$r_c$ be the white section nearest to $c$, split the graph at $r_c$ as
in Section~\ref{subsec:marked-caps}, and let $H_N$ be the cut, the line
$s=r_c/N$. For a probe $z\in H$ with transverse coordinate $y$, let
$z_N$ be the point of $H_N$ whose transverse coordinate is the number
of the form $(x+\tfrac12)/N$, $x\in\mathbb Z$, nearest to $y$, the
smaller one in case of a tie. Then $z_N$ is not on the graph and
$|z_N-z|\le C/N$. For each fixed probe outside $\Diamond$, its realization also lies
outside $\Diamond$ for all sufficiently large $N$. For fixed $m$ and all sufficiently large $N$, the
points $z_N$ are in the same order as the probes $z$ of
$\mathcal P_m(H)$. Together with the two ends of $J_0$ they are the
endpoints of a partition $\mathcal P_{m,N}(H)$ of $J_0$ on the cut
$H_N$, whose cells have length at most $h_m+1/N$. In the notation of
Section~\ref{subsec:all-intervals}, with $C_0=1$,
\begin{equation}\label{eq:realized-partition}
 \mathcal P_{m,N}(H)\in\mathfrak P_{N,h_m}.
\end{equation}

\paragraph{Indicators.}
Order the loops of $\mathcal L_N$ by decreasing diameter, with ties
resolved by a fixed rule. For $j\ge1$ and $z\in\mathcal Z$ put
\[
 \chi_{N,j}(z)=
 \begin{cases}
  1,&\text{if the $j$th loop of $\mathcal L_N$ exists and surrounds }z_N,\\
  0,&\text{otherwise.}
 \end{cases}
\]
On each cut, the transverse coordinate of a realized probe differs
from that of every section vertex by at least $1/(2N)$. The splitting
disks have radius $1/(4N)$, so the deformation avoids the probes on
that cut and preserves the winding numbers that determine their
indicators. Corresponding loops move by at most $C/N$ in a common
parametrization.
Below, an arc of a discrete loop also denotes the corresponding arc in
the split graph.

\paragraph{Probes outside the liquid region.}
By Proposition~\ref{prop:frozen-regions}, for each $\eta>0$ all discrete
loops lie in the $\eta$-neighborhood of $\overline L$ with probability
tending to one. Convexity of this neighborhood puts their interiors
there as well. Thus for every probe
$z\notin\overline L$,
\begin{equation}\label{eq:exterior-probes}
 \PP\{\chi_{N,j}(z)=0\text{ for all }j\}\longrightarrow1 .
\end{equation}

\subsection{Subsequential limits}
\label{subsec:coupling}

\begin{lemma}[Coupling]\label{lem:coupling}
Every subsequence of $(\mathcal L_N)$ has a further subsequence along
which the following holds on a suitable probability space. There are
a random collection $\Gamma\in\mathcal X_\Diamond$, continuous maps
$\gamma_{N,j},\gamma_j:\mathbb T\to\Diamond$ and numbers
$\chi_j(z)\in\{0,1\}$, for $j\ge1$ and $z\in\mathcal Z$, such that:
\begin{enumerate}[label=\textnormal{(\roman*)},leftmargin=*]
\item For each $N$, $(\mathcal L_N,(\chi_{N,j}(z))_{j,z})$ has its
original law, the map $\gamma_{N,j}$ is a parametrization of the $j$th
loop of $\mathcal L_N$ if this loop exists, and $\gamma_{N,j}$ is the
constant map $0$ otherwise.
\item Almost surely, $d_{\rm coll}(\mathcal L_N,\Gamma)\to0$,
\begin{equation}\label{eq:global-representative-limit}
 \gamma_{N,j}\longrightarrow\gamma_j
 \quad\hbox{uniformly on }\mathbb T\ \hbox{for every }j,
\end{equation}
and $\chi_{N,j}(z)=\chi_j(z)$ for all sufficiently large $N$, for every
$j$ and $z$.
\item The nonconstant maps $\gamma_j$, each counted once, form the
collection $\Gamma$.
\item Almost surely, for every $\varepsilon>0$ there is $j_0<\infty$
such that $\diam\gamma_{N,j}<\varepsilon$ for all $j>j_0$ and all $N$
of the subsequence.
\end{enumerate}
\end{lemma}

We call the indices $j\ge1$ \emph{entries}. An entry $j$ is
\emph{nonconstant} if $\gamma_j$ is nonconstant.

\begin{proof}
\emph{Convergence of ranked loops.}
Let $Y$ be the space of single-loop classes with constants included,
with the distance $d_{\rm curv}$. It is complete and separable by the
proof of Lemma~\ref{lem:loop-collection-polish}. Let $u_{N,j}\in Y$ be
the $j$th loop of $\mathcal L_N$, or the constant loop at $0$ if there
is none. For a compact $K\subset\mathcal X_\Diamond$, the loops of all
collections in $K$, together with the constant loops, have compact
closure in $Y$. Indeed, take a sequence of such loops. If their
diameters tend to zero, a subsequence converges to a constant. If
their diameters are at least $\varepsilon>0$, pass to a subsequence
along which their collections converge to some $\Gamma'$; eventually
each loop can be paired with one of the finitely many
loops of $\Gamma'$ of diameter at least $\varepsilon/2$, with pairing error tending to zero.

The indicator array takes values in the compact space
$\{0,1\}^{\mathbb N\times\mathcal Z}$. Proposition~\ref{prop:tightness}
and the preceding argument give tightness of the triples
$(\mathcal L_N,(u_{N,j})_j,(\chi_{N,j}(z))_{j,z})$ in
$\mathcal X_\Diamond\times Y^{\mathbb N}\times
\{0,1\}^{\mathbb N\times\mathcal Z}$. Prokhorov's theorem and the
Skorokhod representation theorem give, along a further subsequence, a
coupling with almost sure limit $(\Gamma,(u_j)_j,(\chi_j(z))_{j,z})$.
Each indicator is therefore eventually equal to its limit.

\emph{Ranks and multiplicities.}
Since $\mathcal L_N\to\Gamma$, for each $\varepsilon>0$ and all
sufficiently large $N$ the loops of $\mathcal L_N$ of diameter at
least $\varepsilon$ can be paired with distinct loop copies in $\Gamma$ of diameter at
least $\varepsilon/2$. There are finitely many of these, and each of
the finitely many earlier collections is finite. So the number of loops
of $\mathcal L_N$ of diameter at least $\varepsilon$ is bounded in $N$.
This gives~(iv) by the ordering of the ranks. To obtain~(iii),
choose a cutoff $\delta>0$ different from every loop diameter of
$\Gamma$ and every $\diam u_j$. By~(iv), only finitely many ranks
can exceed $\delta$, and convergence at each of them identifies the
multiset of their limits. For large $N$, collection convergence gives a bijection
pairing the loops of these ranks with the loop copies in $\Gamma$ of
diameter greater than $\delta$. Let such cutoffs decrease to zero.

\emph{Uniform convergence of parametrizations.}
Pass to a further deterministic subsequence $N_k$ such that
\[
 \PP\!\left\{\max_{j\le k}d_{\rm curv}(u_{N_k,j},u_j)>2^{-k}\right\}
 \le2^{-k}.
\]
By the Borel--Cantelli lemma, for each fixed $j$ the distances
$d_{\rm curv}(u_{N_k,j},u_{N_{k+1},j})$ are almost surely summable.
Each $u_{N_k,j}$ has a polygonal parametrization. Starting with the
first one, reparametrize each next map so that its uniform distance
from its predecessor is less than their distance in $d_{\rm curv}$ plus
$2^{-k}$. To make the choice measurable, take the first admissible
reparametrization in a fixed countable uniformly dense family of
homeomorphisms of $\mathbb T$ of both orientations, for instance those
with rational piecewise-linear lifts. Such a choice exists by uniform
continuity. The resulting maps are almost surely uniformly Cauchy for
every $j$. Write $\gamma_{N_k,j}$ for these maps and $\gamma_j$ for
their uniform limits. Each $\gamma_j$ represents $u_j$.
\end{proof}

From now on we work with a fixed subsequence and the coupling of
Lemma~\ref{lem:coupling}; all limits in $N$ are along this subsequence.

\subsection{Signatures}
\label{subsec:signatures}

The \emph{signature} of an entry $j$ is the set
\[
 S_j=\{z\in\mathcal Z:\chi_j(z)=1\}.
\]
An entry is \emph{visible} if its signature has at least two elements.
If $z,z'\in S_j$, then for large $N$ the $j$th loop surrounds $z_N$ and
$z'_N$, so $\diam\gamma_{N,j}\ge|z_N-z'_N|$ and
$\diam\gamma_j\ge|z-z'|>0$. Hence
\begin{equation}\label{eq:visible-nonconstant}
 \text{every visible entry is nonconstant.}
\end{equation}
For a loop $C$ of $\mathcal C$, the signature
$\mathcal Z\cap\operatorname{Int}(C)$ is dense in its closed Jordan
domain. The boundary of that domain recovers $C$, so distinct target
loops have distinct signatures.

\begin{lemma}[Identification of the visible signatures]
\label{lem:signatures}
On an extension of the probability space of Lemma~\ref{lem:coupling}
there is a random collection with the law of
$\Psi(\CLE_4^{\rm nested}(\DD))$, again denoted $\mathcal C$, and a
bijection $C\mapsto j(C)$ from the loops of $\mathcal C$ onto the set
of visible entries, such that almost surely
\[
 S_{j(C)}=\mathcal Z\cap\operatorname{Int}(C)
 \qquad\text{for every loop $C$ of }\mathcal C .
\]
\end{lemma}

\begin{proof}
\emph{Finite sets of probes.}
For finite $F\subset\mathcal Z$, let $\mathfrak M_N(F)$ be the multiset
of $\{z\in F:\chi_{N,j}(z)=1\}$ over entries for which this set has
at least two elements. Define $\mathfrak M(F)$ using $\chi_j$, and
$\mathfrak M_{\mathcal C}(F)$ using the interiors of the target loops.

If $|F|<2$ all three multisets are empty. Otherwise every discrete
loop contributing to $\mathfrak M_N(F)$ has diameter at least
$\tfrac12\min\{|z-z'|:z\ne z'\in F\}$ for large $N$. By
Lemma~\ref{lem:coupling}(iv) only a bounded set of ranks can contribute.
On these ranks the indicators eventually equal their limits, so
$\mathfrak M_N(F)=\mathfrak M(F)$ eventually, almost surely.

The reduced lamination on $F\cap L$ determines the multiset of
enclosed subsets of size at least two. Applying
Theorem~\ref{thm:aztec-finite-puncture} with the moving punctures
$z_N$ gives $\mathfrak M_N(F\cap L)\Rightarrow
\mathfrak M_{\mathcal C}(F\cap L)$ whenever $|F\cap L|\ge2$.
All other probes of $F$ lie outside $\overline L$. By
\eqref{eq:exterior-probes} they are surrounded by no discrete loop
with probability tending to one, and by no target loop almost surely.
Thus $\mathfrak M(F)$ and $\mathfrak M_{\mathcal C}(F)$ have the same law.

For $F'\subset F$, the multiset $\mathfrak M(F')$ is obtained from
$\mathfrak M(F)$ by intersecting each of its sets with $F'$ and
keeping those with at least two elements; the same holds for
$\mathfrak M_{\mathcal C}$. Enumerate $\mathcal Z=\{z_1,z_2,\ldots\}$
and put $F_k=\{z_1,\ldots,z_k\}$. It follows that the sequences
\[
 R=(\mathfrak M(F_k))_{k\ge1}
 \qquad\text{and}\qquad
 R_{\mathcal C}=(\mathfrak M_{\mathcal C}(F_k))_{k\ge1}
\]
have the same law. We call $R$ and $R_{\mathcal C}$ the
\emph{records}. They take values in a countable product of countable
sets.

\emph{The record determines the visible signatures.}
For fixed $a<b$, only finitely many signatures contain both $z_a$ and
$z_b$. The submultisets of $\mathfrak M(F_k)$ containing this pair,
$k\ge b$, determine the counting measure on these signatures: they
give its values on all cylinder sets in $\{0,1\}^{\mathcal Z}$.
Enumerate the pairs $a<b$ and assign each signature to its first
pair in this enumeration. In this way the record recovers the
multiset of all visible signatures, including multiplicities. The
same construction applies to $R_{\mathcal C}$.

\emph{Coupling with $\mathcal C$.}
On an extension of the probability space, sample $\mathcal C$ from
its conditional law given its record, evaluated at $R$. Then
$R_{\mathcal C}=R$ almost surely. The visible signatures thus agree
with the signatures of the target loops, which are pairwise distinct.
This gives the required bijection.
\end{proof}

\subsection{An almost sure consequence of the opposite-crossing estimate}
\label{subsec:no-witness}

Let $H$ be a test line, say $s=c$, and let $\rho_-,\rho_+>0$ be
rational with $0<c-\rho_-<c+\rho_+<1$. Tested arcs are defined with
the deep regions $\{s\le c-\rho_-\}$ and $\{s\ge c+\rho_+\}$, and
indices with the cut $H_N$, as in
Section~\ref{subsec:opposite-setting}. For $m\ge1$ let
\[
 B_{N,m}(H,\rho_\pm)=B^{\rm all}_N\bigl(\mathcal P_{m,N}(H)\bigr)
\]
be the event that some loop of $\mathcal L_N$ has two tested arcs
$\alpha,\beta$ with $\iota_I(\alpha)\iota_I(\beta)<0$ for some cell
$I$ of $\mathcal P_{m,N}(H)$. For test lines $y=c$ the definitions are
the reflected ones.

\begin{lemma}[No persistent opposite crossings]\label{lem:no-witness}
Almost surely, for every test line $H$ and all rational
$\rho_-,\rho_+$ as above, the following does \emph{not} happen: there
is $m_0$ such that for every $m\ge m_0$ the event
$B_{N,m}(H,\rho_\pm)$ occurs for all sufficiently large $N$.
\end{lemma}

\begin{proof}
Fix $H$ and $\rho_\pm$. By~\eqref{eq:realized-partition} and
Theorem~\ref{thm:global-opposite} with $\alpha=\tfrac12$,
$\limsup_N\PP(B_{N,m})\le Ch_m^{1/2}$ for all large $m$. Fatou's lemma,
applied twice, gives
\[
 \PP\bigl(\liminf_{m\to\infty}\liminf_{N\to\infty}B_{N,m}\bigr)
 \le\liminf_{m\to\infty}\liminf_{N\to\infty}\PP(B_{N,m})=0 .
\]
The event on the left is the one described in the lemma. There are
countably many choices of $H$ and $\rho_\pm$.
\end{proof}

The cell $I$ for which $B_{N,m}$ occurs may depend on $m$ and $N$. We
use the lemma together with the following elementary fact.

\begin{lemma}[Index and winding numbers]\label{lem:winding-index}
In the coordinates of a cut, let $u$ be its longitudinal coordinate
and $v$ its transverse coordinate, so the cut is $u=\mathrm{const}$
and $\iota_I$ counts crossings towards increasing $u$ positively.
Let $I=[p_-,p_+]$ run towards increasing $v$ and let $\lambda$ be a
closed polygonal path transverse to the cut and avoiding $p_-,p_+$.
With winding numbers taken in the oriented $(u,v)$ plane,
\[
 \iota_I(\lambda)=\operatorname{wind}(\lambda,p_+)
                  -\operatorname{wind}(\lambda,p_-).
\]
For a simple loop the index belongs to $\{-1,0,1\}$ and vanishes
whenever the two endpoints are both inside or both outside the loop.
\end{lemma}

\begin{proof}
Along $I$, the winding number jumps by $+1$ at a crossing towards
increasing $u$ and by $-1$ at a crossing towards decreasing $u$.
For an ordinary cut $(u,v)=(s,y)$; a reflected cut uses the oriented
coordinates of the reflected graph. In either convention a simple
loop has winding number zero outside and a constant winding number
of absolute value one inside.
\end{proof}

\subsection{Traces of the limiting loops}
\label{subsec:traces}

\begin{lemma}[Constant signature on a disk]
\label{lem:constant-signatures}
Almost surely the following holds for every entry $j$ and every open
disk $B$ such that $\chi_j$ is constant on $\mathcal Z\cap B$. If
$[u,v]$ is a parameter interval with $\gamma_j([u,v])\subset B$, then
$\gamma_j(u)=\gamma_j(v)$.
\end{lemma}

The disk $B$ may meet $\partial L$ or the frozen regions: the set
$\mathcal Z$ contains probes outside $\overline L$, and
Theorem~\ref{thm:global-opposite} concerns partitions of the entire
section.

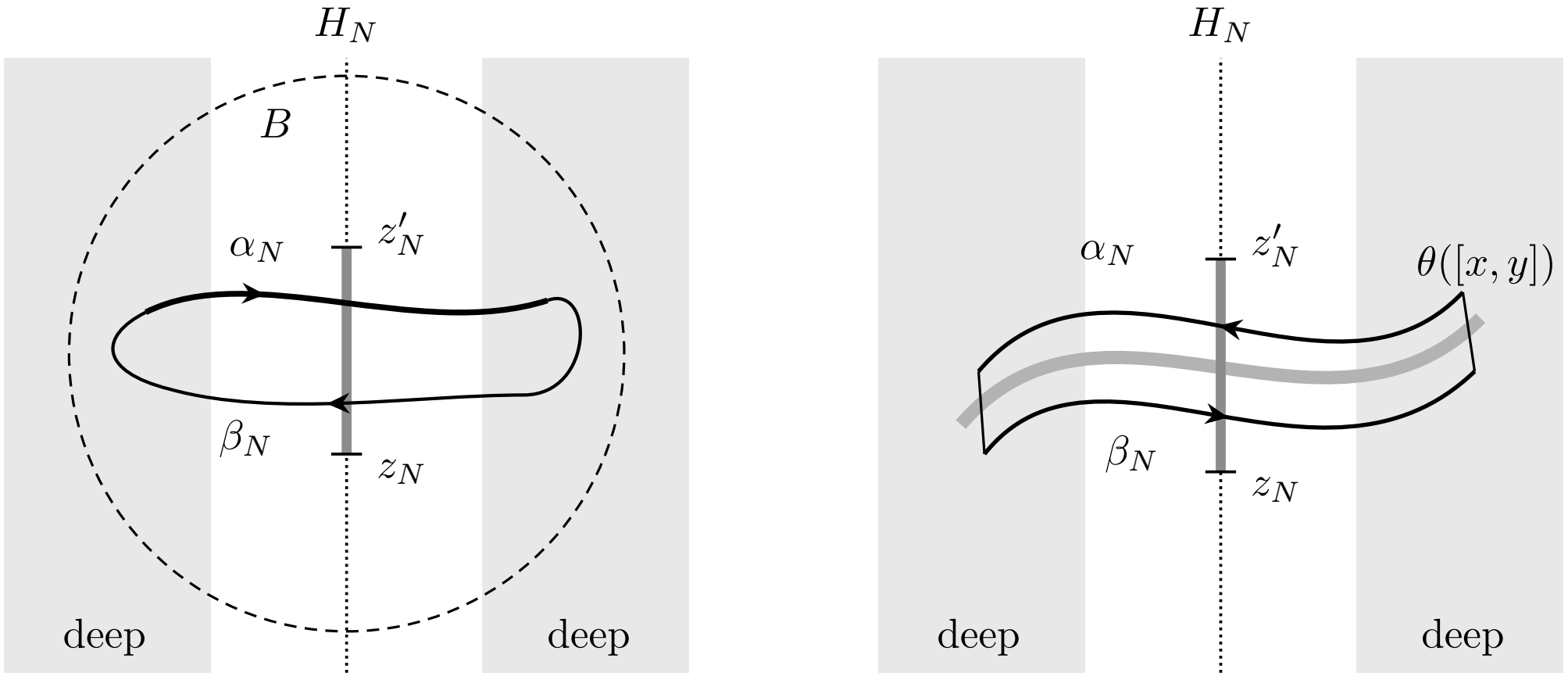
\begin{figure}[htbp]
 \centering
 \begin{tikzpicture}[x=1cm,y=1cm,
   arr/.style={postaction={decorate},
     decoration={markings,mark=at position #1 with
       {\arrow{Stealth[length=5.5pt,width=5pt]}}}},
   pr/.style={line width=0.7pt}]
  \begin{scope}
   \fill[black!9] (-2.9,-2.6) rectangle (-1.15,2.6);
   \fill[black!9] (1.15,-2.6) rectangle (2.9,2.6);
   \draw[densely dashed,line width=0.6pt] (0,0.1) circle (2.35);
   \node at (-0.6,2.05) {\small $B$};
   \draw[densely dotted,line width=0.7pt] (0,-2.6)--(0,2.6);
   \node[above] at (0,2.6) {\small $H_N$};
   \draw[line width=2.4pt,black!45] (0,-0.75)--(0,1.0);
   \draw[pr] (-0.13,1.0)--(0.13,1.0);
   \draw[pr] (-0.13,-0.75)--(0.13,-0.75);
   \node[anchor=west] at (0.12,1.12) {\small $z'_N$};
   \node[anchor=west] at (0.12,-0.9) {\small $z_N$};
   \draw[line width=1.4pt,arr=0.3]
     (-1.7,0.45) .. controls (-0.8,0.9) and (0.6,0.2) .. (1.7,0.55);
   \draw[line width=0.8pt,arr=0.55]
     (1.7,0.55) .. controls (2.1,0.7) and (2.1,-0.25) .. (1.5,-0.25)
       .. controls (0.6,-0.25) and (-0.6,-0.45) .. (-1.5,-0.2)
       .. controls (-2.1,-0.05) and (-2.1,0.25) .. (-1.7,0.45);
   \node at (-0.75,0.98) {\small $\alpha_N$};
   \node at (-0.85,-0.62) {\small $\beta_N$};
   \node at (-2.05,-2.3) {\small deep};
   \node at (2.05,-2.3) {\small deep};
  \end{scope}
  \begin{scope}[xshift=7.4cm]
   \fill[black!9] (-2.9,-2.6) rectangle (-1.15,2.6);
   \fill[black!9] (1.15,-2.6) rectangle (2.9,2.6);
   \draw[densely dotted,line width=0.7pt] (0,-2.6)--(0,2.6);
   \node[above] at (0,2.6) {\small $H_N$};
   \draw[line width=3.2pt,black!30]
     (-2.2,-0.5) .. controls (-1.0,0.9) and (0.9,-0.9) .. (2.2,0.4);
   \node at (2.25,0.85) {\small $\theta([x,y])$};
   \draw[line width=2.4pt,black!45] (0,-0.9)--(0,0.9);
   \draw[pr] (-0.13,0.9)--(0.13,0.9);
   \draw[pr] (-0.13,-0.9)--(0.13,-0.9);
   \node[anchor=west] at (0.12,1.02) {\small $z'_N$};
   \node[anchor=west] at (0.12,-1.05) {\small $z_N$};
   \draw[line width=1.0pt,arr=0.5]
     (2.05,0.62) .. controls (0.95,-0.55) and (-0.95,1.25) .. (-2.05,-0.05);
   \draw[line width=1.0pt,arr=0.5]
     (-2.0,-0.75) .. controls (-1.05,0.45) and (0.85,-1.3) .. (2.15,-0.05);
   \draw[line width=0.6pt] (-2.05,-0.05)--(-2.0,-0.75);
   \draw[line width=0.6pt] (2.15,-0.05)--(2.05,0.62);
   \node at (-0.95,0.95) {\small $\alpha_N$};
   \node at (-0.75,-0.75) {\small $\beta_N$};
   \node at (-2.05,-2.3) {\small deep};
   \node at (2.05,-2.3) {\small deep};
  \end{scope}
 \end{tikzpicture}
 \caption{The two uses of the opposite-crossing estimate. The shaded
 strips are the deep regions and the thick gray segment is a cell
 $I=[z_N,z'_N]$ of the cut $H_N$. Left
 (Lemma~\ref{lem:constant-signatures}): if $z_N$ and $z'_N$ are both
 outside or both inside the loop, the index of the loop on $I$ is zero,
 so an arc $\alpha_N$ with nonzero index forces the opposite index on
 the complementary arc $\beta_N$. Right
 (Lemma~\ref{lem:no-backtracking}): if the limit backtracks on the
 Jordan arc $\theta([x,y])$, two disjoint arcs $\alpha_N,\beta_N$ of
 the same loop run along it in opposite directions; closing them by two
 short segments gives a path with winding number zero around $z_N$ and
 $z'_N$.}
 \label{fig:identification}
\end{figure}

\begin{proof}
Work on the almost sure event of Lemma~\ref{lem:no-witness}, and
suppose that $\gamma_j(u)\ne\gamma_j(v)$. The two points differ in one
coordinate, say $s$; otherwise use the reflected construction. Choose a
test line $H=\{s=c\}$ and rational $\rho_\pm$ such that $\gamma_j(u)$
and $\gamma_j(v)$ lie in the interiors of the two different deep
regions. Let
\[
 \alpha_N=\gamma_{N,j}|_{[u,v]},\qquad
 \beta_N=\gamma_{N,j}|_{[v,u+1]}
\]
be the corresponding arc of the $j$th discrete loop and its
complementary arc, both directed by the parametrization
(Figure~\ref{fig:identification}, left). By the uniform
convergence~\eqref{eq:global-representative-limit}, for large $N$ the
endpoints of $\alpha_N$ lie in the two deep regions, so $\alpha_N$ and
$\beta_N$ are tested arcs of the same loop.

Since the endpoints of $\alpha_N$ are on opposite sides of $H_N$,
its indices sum to $\pm1$ over every partition of $J_0$. The compact
image $\gamma_j([u,v])$ has positive distance from $B^c$ and from
the two ends of $J_0$ on $H$. Choose a closed neighborhood $K\subset B$
of this image disjoint from those ends. For all large $N$, $\alpha_N$
stays a positive distance from $K^c$, uniformly in $N$. Choose $m_0$
so that, for all $m\ge m_0$, every cell of $\mathcal P_m(H)$ meeting
$K$ has both endpoints in $B$ and neither endpoint at an end of $J_0$.

Fix $m\ge m_0$. The realization moves points by at most $C/N$. So for
large $N$ every crossing of $\alpha_N$ with $H_N$ lies in a cell of
$\mathcal P_{m,N}(H)$ that realizes a cell meeting $K$. Hence one of
these cells, $I=[z_N,z'_N]$, has $\iota_I(\alpha_N)\ne0$. There are finitely many such cells, their
endpoints $z,z'$ are probes in $B$, and $\chi_j$ is constant on
$\mathcal Z\cap B$. So for large $N$,
$\chi_{N,j}(z)=\chi_{N,j}(z')$: the points $z_N,z'_N$ are both inside
or both outside the $j$th loop. By Lemma~\ref{lem:winding-index} the
index of the whole loop on $I$ is zero, that is,
\[
 \iota_I(\alpha_N)+\iota_I(\beta_N)=0 .
\]
Thus $\iota_I(\alpha_N)\iota_I(\beta_N)<0$, and $B_{N,m}(H,\rho_\pm)$
occurs for all large $N$. This holds for every $m\ge m_0$ and
contradicts Lemma~\ref{lem:no-witness}.
\end{proof}

\begin{corollary}[Traces]\label{cor:traces}
Almost surely the following holds.
\begin{enumerate}[label=\textnormal{(\roman*)},leftmargin=*]
\item Let $j$ be a nonconstant entry and $F$ a closed set such that
every point outside $F$ has a neighborhood $B$ with $\chi_j$ constant
on $\mathcal Z\cap B$. Then $\gamma_j(\mathbb T)\subset F$.
\item Every nonconstant entry is visible.
\item Let $C$ be a loop of $\mathcal C$ and $j=j(C)$. Then
$\gamma_j(\mathbb T)=C$, and $\gamma_j$, as a map from $\mathbb T$ to
the Jordan curve $C$, has degree $+1$ or $-1$.
\end{enumerate}
\end{corollary}

\begin{proof}
(i) Suppose $\gamma_j(t_0)\notin F$. Choose an open disk $B$ about
$\gamma_j(t_0)$ on whose probes $\chi_j$ is constant. By
Lemma~\ref{lem:constant-signatures}, $\gamma_j$ is constant on every
parameter interval that it maps into $B$. Hence the closed set
$\{t:\gamma_j(t)=\gamma_j(t_0)\}$ is open. As $\mathbb T$ is connected,
$\gamma_j$ is constant, a contradiction.

(ii) Let $j$ be nonconstant. If $S_j$ were empty, (i) would apply with
$F=\varnothing$. If $S_j=\{z\}$, it would apply with $F=\{z\}$ and
$\gamma_j$ would be constant.

(iii) By~\eqref{eq:visible-nonconstant} the entry $j$ is nonconstant.
Its signature $\mathcal Z\cap\operatorname{Int}(C)$ is constant near
every point off $C$, so $\gamma_j(\mathbb T)\subset C$ by~(i). Choose
$z\in S_j$. For large $N$ the $j$th loop surrounds $z_N$, so
$\operatorname{wind}(\gamma_{N,j},z_N)=\pm1$. Since
$z\in\operatorname{Int}(C)$ is at positive distance from
$\gamma_j(\mathbb T)$, uniform convergence gives
$\operatorname{wind}(\gamma_j,z)=\pm1$. This number is, up to sign, the
degree of $\gamma_j:\mathbb T\to C$. A map of nonzero degree to a circle is
surjective.
\end{proof}

\subsection{Order of traversal}
\label{subsec:monotone}

\begin{lemma}[No backtracking]\label{lem:no-backtracking}
Almost surely the following holds for every loop $C$ of $\mathcal C$
and $j=j(C)$. Let $\theta:\mathbb R/\mathbb Z\to C$ be a homeomorphism,
oriented so that $\gamma_j$ has degree $+1$, and let
$\phi:\mathbb R\to\mathbb R$ be any continuous lift satisfying
\[
 \gamma_j(t)=\theta(\phi(t)\bmod1),\qquad \phi(t+1)=\phi(t)+1 .
\]
Then $\phi$ is nondecreasing, and
$d_{\rm curv}(\gamma_j,\theta)=0$.
\end{lemma}

\begin{proof}
Work on the almost sure events of Lemma~\ref{lem:no-witness} and
Corollary~\ref{cor:traces}, and on the event that no probe lies on a
loop of $\mathcal C$. Suppose that $\phi(u)>\phi(v)$ for some
$u<v<u+1$. Choose $x<y$ with
\[
 \phi(v)<x<y<\phi(u),\qquad y-x<1 .
\]

\emph{Two passages in opposite directions.}
Let $v_1$ be the first time after $u$ with $\phi(v_1)=x$, and $u_1$
the last time before $v_1$ with $\phi(u_1)=y$. Then
$u<u_1<v_1<v$ and $\phi([u_1,v_1])=[x,y]$. Since
$\phi(u+1)=\phi(u)+1>y$, we can let $v_2$ be the first time after $v$
with $\phi(v_2)=y$ and $u_2$ the last time before $v_2$ with
$\phi(u_2)=x$. Then $v<u_2<v_2<u+1$ and $\phi([u_2,v_2])=[x,y]$. So
\[
 \alpha=\gamma_j|_{[u_1,v_1]}\quad\text{runs from $\theta(y)$ to
 $\theta(x)$},\qquad
 \beta=\gamma_j|_{[u_2,v_2]}\quad\text{runs from $\theta(x)$ to
 $\theta(y)$},
\]
both inside the Jordan arc $\theta([x,y])$, and their parameter
intervals are disjoint.

\emph{A test line.}
The points $\theta(x)\ne\theta(y)$ differ in one coordinate, say $s$.
Choose a test line $H=\{s=c\}$ and rational $\rho_\pm$ such that
$\theta(x)$ and $\theta(y)$ lie in the interiors of the two different
deep regions. Put $\alpha_N=\gamma_{N,j}|_{[u_1,v_1]}$ and
$\beta_N=\gamma_{N,j}|_{[u_2,v_2]}$. For large $N$ these are two tested
arcs of the $j$th loop, directed by the same orientation.

\emph{A closed path with zero winding numbers.}
Let $\lambda_N$ be the closed path that follows $\alpha_N$, then the
straight segment from the end of $\alpha_N$ to the start of $\beta_N$,
then $\beta_N$, then the straight segment from the end of $\beta_N$ to
the start of $\alpha_N$ (Figure~\ref{fig:identification}, right). For
large $N$ each of the two segments lies in a deep region, which is
convex, so it does not meet $H_N$. As $N\to\infty$, $\lambda_N$
converges uniformly, in a suitable parametrization, to the closed path
$\lambda$ that follows $\alpha$ and then $\beta$. The path $\lambda$
lies in the arc $\theta([x,y])$, which is contractible. Hence
$\operatorname{wind}(\lambda,p)=0$ for every point $p$ off this arc.

Fix $m$. The endpoints of the cells of $\mathcal P_m(H)$ avoid the
limiting arc: the probes miss $C$, and the ends of $J_0$ lie outside
$\Diamond$. Uniform convergence of the paths and convergence of the
realized endpoints preserve these zero winding numbers for all
sufficiently large $N$. By
Lemma~\ref{lem:winding-index}, and since the two segments do not meet
$H_N$,
\[
 \iota_I(\alpha_N)+\iota_I(\beta_N)=\iota_I(\lambda_N)=0
 \qquad\text{for every cell $I$ of }\mathcal P_{m,N}(H).
\]
The endpoints of $\alpha_N$ are on different sides of $H_N$, so
$\sum_I\iota_I(\alpha_N)=\pm1$ and some cell has
$\iota_I(\alpha_N)\ne0$. For this cell
$\iota_I(\alpha_N)\iota_I(\beta_N)<0$. Thus $B_{N,m}(H,\rho_\pm)$
occurs for all large $N$, for every $m$. This contradicts
Lemma~\ref{lem:no-witness}. Hence $\phi$ is nondecreasing.

\emph{The curve distance.}
For $0<\epsilon<1$ the map $\phi_\epsilon(t)=(1-\epsilon)\phi(t)+\epsilon t$
is strictly increasing with $\phi_\epsilon(t+1)=\phi_\epsilon(t)+1$, so
$\theta(\phi_\epsilon\bmod1)$ is a reparametrization of $\theta$. It
converges uniformly to $\gamma_j$ as $\epsilon\downarrow0$. Hence
$d_{\rm curv}(\gamma_j,\theta)=0$.
\end{proof}

\subsection{Proof of Theorem~\ref{thm:global-curves}}
\label{subsec:conclusion}

By Proposition~\ref{prop:tightness} it suffices to identify an
arbitrary subsequential limit $\Gamma$. Use the coupling of
Lemmas~\ref{lem:coupling} and~\ref{lem:signatures}. A nonconstant
entry is visible by Corollary~\ref{cor:traces}(ii), and a visible
entry is nonconstant by~\eqref{eq:visible-nonconstant}. The signature
bijection therefore pairs every loop of $\mathcal C$ with one
nonconstant entry, with multiplicities preserved. By
Corollary~\ref{cor:traces}(iii) and Lemma~\ref{lem:no-backtracking},
each entry represents the same unparametrized curve as its paired loop. The entries exhaust
$\Gamma$ by Lemma~\ref{lem:coupling}(iii). Thus $\Gamma=\mathcal C$
almost surely in $\mathcal X_\Diamond$. \qed

\subsection{Containment in the liquid region}
\label{subsec:containment}

\begin{corollary}[Macroscopic containment and compact windows]
\label{cor:arctic-containment}
For every $\varepsilon>0$, writing
$L_{-\eta}=\{z\in L:\operatorname{dist}(z,\partial L)>\eta\}$,
\begin{equation}\label{eq:macroscopic-arctic-containment}
 \lim_{\eta\downarrow0}\limsup_{N\to\infty}
 \PP\{\exists\gamma\in\mathcal L_N:
          \diam\gamma\ge\varepsilon,\ \gamma\not\subset L_{-\eta}\}=0.
\end{equation}
Fix a compact $K\subset L$ and a cutoff $d>0$. Suppose that almost
surely no limiting loop has diameter $d$ and every limiting loop
of diameter greater than $d$ meeting $K$ enters $\operatorname{int}K$.
Then the collection of whole loops of diameter greater than $d$
that meet $K$ converges in law to the corresponding limiting collection.
\end{corollary}

\begin{proof}
In a coupling realizing collection convergence, every discrete loop
of diameter at least $\varepsilon$ is eventually paired with a target loop
of diameter at least $\varepsilon/2$. This finite family of target
loops lies in $L$. Its distance $\Delta_\varepsilon$ from $L^c$ is
strictly positive almost surely; take $\Delta_\varepsilon=\infty$ for
the empty family. On $\{\Delta_\varepsilon>2\eta\}$ the pairing
errors are eventually smaller than $\eta$, so every such discrete
loop lies in $L_{-\eta}$. The upper limit of the probability
in~\eqref{eq:macroscopic-arctic-containment} is at most
$\PP\{\Delta_\varepsilon\le2\eta\}$, which tends to zero.

For the second assertion, only the finitely many target loops of
diameter at least $d/2$ can be paired with a discrete loop above the cutoff
for large $N$. The two diameters in each pair eventually lie on the same
side of $d$. A target loop above $d$ either enters
$\operatorname{int}K$ or has positive distance from $K$. Its partner
therefore eventually has the same selection status.
\end{proof}

\section*{Acknowledgments}
This work was supported by the National Science Foundation under grant
no.~2512154 and by the Simons Foundation under grant
no.~\mbox{SFI-MPS-TSM-00025603}.

\section*{AI usage}
Generative AI tools assisted in developing the proof and writing the paper. The author takes full responsibility for for its contents.

\end{document}